\documentclass[11pt,a4paper]{article}
\usepackage{amscd}
\usepackage{enumerate}

\usepackage{amsmath, amssymb, latexsym}
\usepackage{bbm}
\usepackage{amsfonts}
\usepackage{amsthm}
\usepackage{relsize}
\usepackage{setspace}
\usepackage{geometry}
\usepackage{url}
\usepackage{xspace}
\usepackage{tocloft}
\usepackage{graphics}
\usepackage{graphicx}
\usepackage{lscape}
\usepackage{microtype}
\usepackage{ulem}

\usepackage[usenames, dvipsnames]{color}
\usepackage[utf8]{inputenc}
\usepackage{tikz-cd}

\usepackage[pagebackref=true]{hyperref}
\usepackage[alphabetic]{amsrefs}
\usepackage[english]{babel}

\usepackage{authblk}

\newtheorem{theorem}{Theorem}[section]
\newtheorem{proposition}[theorem]{Proposition}

\newtheorem{corollary}[theorem]{Corollary}
\newtheorem{lemma}[theorem]{Lemma}

\theoremstyle{definition}
\newtheorem{definition}[theorem]{Definition}

\newtheorem{remark}[theorem]{Remark}

\theoremstyle{problem}

\newcommand{\Aut}{\mathrm{Aut}}

\newcommand{\Opp}{\mathrm{Opp}}
\newcommand{\Ker}{\mathrm{Ker}}

\newcommand{\proj}{\mathrm{proj}}
\newcommand{\RR}{\mathbb{R}}

\newcommand{\QQ}{\mathbb{Q}}
\newcommand{\CC}{\mathbb{C}}
\newcommand{\NN}{\mathbb{N}}

\newcommand{\cat}{$\mathrm{CAT}(0)$\xspace}
\newcommand{\Min}{\mathrm{Min}}
\newcommand{\Ch}{\mathrm{Ch}}

\newcommand{\Stab}{\mathrm{Stab}}

\def\Isom{\operatorname{Isom}}

\newcommand{\id}{\operatorname{id}}

\newcommand{\SL}{\operatorname{SL}}

\newcommand{\dist}{\operatorname{dist}}

\newcommand{\st}{\operatorname{St}}

\newcommand{\limch}{\operatorname{^{Ch}\lim}}

\def\og{\leavevmode\raise.3ex\hbox{$\scriptscriptstyle\langle\!\langle$~}}
\def\fg{\leavevmode\raise.3ex\hbox{~$\!\scriptscriptstyle\,\rangle\!\rangle$}}

\title{Chabauty limits of fixed point groups of $p$-adic involutions}

\author{Corina Ciobotaru\thanks{cociobotaru@math.au.dk}}
\date{September 14, 2025}

\begin{document}
\maketitle
\begin{abstract}
We study Chabauty limits of the fixed-point group of $k$-points \( H_k \)  associated with an involutive $k$-automorphism \( \theta \) of a connected linear reductive group \( G \) defined over a non-Archimedean local field \( k \) of characteristic zero. Leveraging the geometry of the Bruhat--Tits building, the structure of \( (\theta,k) \)-split tori, and the $K\mathcal{B}_kH_k$ decomposition of \( G_k \), we establish that any nontrivial Chabauty limit \( L \) of \( H_k \) is \( G_k \)-conjugate to a subgroup of
\[
U_{\sigma_+}(k) \rtimes (\Ker(\alpha)^0 \cdot (H_k \cap M_{\sigma_-, \sigma_+})) \leq P_{\sigma_+}(k),
\]
where \( \alpha \) is a projection map arising from a Levi factor \( M_{\sigma_-, \sigma_+} \) of a parabolic subgroup \( P_{\sigma_+} \subset G \), and \( \Ker(\alpha)^0 \) denotes the subgroup of elliptic elements in the kernel of \( \alpha \). Our analysis distinguishes between elliptic and hyperbolic elements and constructs explicit unipotent elements in the limit group \( L \) using the Moufang property of \( G_k \). Furthermore, we show that \( L \) acts transitively on the set of ideal simplices opposite to \( \sigma_+ \). These results yield a detailed description of the Chabauty compactification of \( H_k \), and provide new insights into its interaction with the non-Archimedean geometry of \( G_k \).
\end{abstract}

\tableofcontents

\section{Introduction}
Let \( G \) be a connected linear reductive group defined over a field \( k \) of characteristic not $2$. A \textit{symmetric $k$-variety} is the homogeneous space \( G_k/H_k \), where \( \theta \) is an involutive $k$-automorphism of \( G \), and \( H \) is its fixed-point subgroup. The groups \( G_k \) and \( H_k \) denote the \( k \)-rational points of \( G \) and \( H \), respectively. Such symmetric spaces arise naturally in the study of automorphic forms and harmonic analysis, and they play a pivotal role in the representation theory of reductive groups, the Langlands correspondence, and the formulation of Plancherel-type decompositions for symmetric spaces.

In the Archimedean setting, where \( k = \mathbb{R} \) or \( \mathbb{C} \), the homogeneous space \( G_k/H_k \) defines an \textit{affine symmetric space}, extending the classical theory of Riemannian symmetric spaces. Prominent examples include spaces associated with quadratic forms of signature \( (p,q) \), where the fixed-point group \( H_k \) corresponds to the orthogonal group \( \mathrm{O}(p,q) \). These spaces encompass familiar geometries such as spherical, hyperbolic, de Sitter, and anti-de Sitter geometries. Although these models differ in curvature, they are linked through \textit{geometric transitions} -- continuous deformations of geometric structures that yield distinct limiting geometries. Such transitions are of intrinsic interest in both mathematics and theoretical physics, where they model phenomena such as the passage from general to special relativity, or from quantum mechanics to classical mechanics.

A systematic classification of geometric limits arising from affine symmetric spaces over \( \mathbb{R} \) was developed by Cooper--Danciger--Wienhard~\cite{CDW}, who studied degenerations of groups such as \( \mathrm{SO}(p,q) \subset \mathrm{GL}(n, \mathbb{R}) \) via root space decompositions of real Lie algebras.

In line with Klein’s Erlangen Program, geometric structures are determined by their symmetry groups. As a result, understanding limits of geometries reduces to analyzing limits of Lie groups under suitable topologies.

The \textit{Chabauty topology}, introduced by Claude Chabauty in 1950~\cite{Ch}, provides a compact and natural topology on the space \( \mathcal{S}(G) \) of closed subgroups of a locally compact group \( G \). Originally motivated by the study of the set of lattices and its compactness properties, Chabauty’s framework generalizes Mahler’s criterion for lattices in \( \mathbb{R}^n \) and offers a powerful tool for investigating group-theoretic and geometric transitions.

A central feature of the Chabauty topology is the compactness of $\mathcal{S}(G)$, which implies that every sequence of closed subgroups admits a convergent subsequence. The resulting limits, known as \textbf{Chabauty limits}, provide a robust framework for studying degenerations of algebraic and geometric structures. Thus, the notion of limits of groups -- and by extension, limits of geometries -- is both meaningful and rich. The principal challenge lies not in establishing the existence of such limits, but in identifying and characterizing the geometric and algebraic nature of the subgroups that arise.

Despite its foundational importance, the global topology of $\mathcal{S}(G)$ remains elusive and is fully understood only in a few cases. For example, it is known that $\mathcal{S}(\mathbb{R}) \cong [0, \infty]$~\cite{HP}, and $\mathcal{S}(\mathbb{R}^2) \cong \mathbb{S}^4$~\cite{HP}. Kloeckner~\cite{Kloeckner} has shown that while $\mathcal{S}(\mathbb{R}^n)$ fails to be a manifold for $n > 2$, it is a stratified space in the sense of Goresky--MacPherson and is simply connected. Nevertheless, a complete description of $\mathcal{S}(\mathbb{R}^n)$ remains an open problem.

Substantial progress has been made in understanding the Chabauty closures of specific families of subgroups, including abelian subgroups~\cite{Baik1, Baik2, Haettel, Leitner_sl3, Leitner_sln, Htt_2}, connected subgroups~\cite{LL}, and lattices~\cite{BLL, Wang}. More recently, attention has turned to the $p$-adic setting. Bourquin and Valette~\cite{BV} have determined the homeomorphism type of $\mathcal{S}(\mathbb{Q}_p^*)$, and Cornulier~\cite{Cornulier} has investigated structural properties of $\mathcal{S}(G)$ for locally compact abelian groups. Chabauty limits of groups acting on trees have been studied in~\cite{CR, Stulemeijer}, and several open questions are posed in~\cite{CM}.

In joint work with Leitner and Valette, the author has investigated limits of families of subgroups in $\mathrm{SL}(n, \mathbb{Q}_p)$, including parahoric subgroups~\cite{CiLe_p}, Cartan subgroups~\cite{CLV}, and fixed-point groups of involutions in $\mathrm{SL}(2, F)$ for local fields $F$~\cite{CiLe}. Furthermore, compactifications of Bruhat--Tits buildings have been explored in~\cite{GR}, contributing to a deeper understanding of the interplay between group-theoretic limits and non-Archimedean geometry.

This article marks the third contribution in an ongoing research program, following~\cite{CiLe} and~\cite{Cio24}, and establishes a $p$-adic analogue of Chabauty limits for fixed-point subgroups of involutions -- extending the framework developed by Cooper, Danciger, and Wienhard~\cite{CDW} in the real and complex settings. The work in~\cite{CiLe} addresses the specific case of \( \mathrm{SL}(2,F) \), where \( F \) is a local field. In parallel,~\cite{Cio24} explores key dynamical features of certain hyperbolic automorphisms acting on locally finite thick affine buildings, with particular emphasis on the cone topology of the associated spherical building at infinity. It also provides broad topological and algebraic conditions under which a Chabauty limit contains the unipotent radical of a parabolic subgroup. 

While this article builds on the foundational results of \cite{HW93, BeOh, Cio24} and draws conceptual inspiration from the structure of Chabauty limits computed in~\cite{CDW}, the techniques and overall proof strategy developed here are entirely new. A central challenge in the non-Archimedean setting is that, unlike in the real or complex case, the Lie algebra of a \( p \)-adic Lie group does not uniquely determine the group via the exponential map. This discrepancy arises from the fact that non-Archimedean local fields admit small open multiplicative subgroups that share the same Lie algebra but correspond to distinct group-theoretic structures. Consequently, the methods based on Lie algebras and exponential maps employed in~\cite{CDW} are not applicable in the non-Archimedean context.  Instead, the approach developed in this article is entirely geometric, relying heavily on the action of the reductive group on its associated Tits and Bruhat--Tits buildings.

\medskip
Throughout this article, we let \( G \) denote a connected linear reductive group defined over a non-Archimedean local field \( k \) of characteristic zero. The assumption of characteristic zero is essential, as it ensures the applicability of the Inverse Function Theorem in the setting of analytic manifolds over non-Archimedean fields (see Appendix~\ref{appen:B}, and more precisely~\cite{PlaRa}, middle of p.~110).

Let \( \theta \) be an involutive \( k \)-automorphism of \( G \), and define the fixed-point subgroup \( G^{\theta} := \{ h \in G \mid \theta(h) = h \} \). Let \( H := (G^{\theta})^{\circ} \) denote the identity component of \( G^{\theta} \). Since \( \theta \) is defined over \( k \), it follows that \( H \) is also defined over \( k \) (see~\cite[Proposition 1.6]{HW93}). Moreover, by~\cite[Theorem 2.1]{PraYu}, the group \( H \) is reductive. We denote by \( G_k \) and \( H_k \) the groups of \( k \)-rational points of \( G \) and \( H \), respectively.

This article examines the Chabauty topology on the space \( \mathcal{S}(G_k) \) of closed subgroups of \( G_k \), with a particular focus on the Chabauty limits of \( H_k \) under conjugation by sequences in \( G_k \). Our primary objective is to classify all nontrivial Chabauty limits of \( H_k \) (i.e. limits that are not $G_k$-conjugate to $H_k$) and to describe their structure in terms of parabolic subgroups and their unipotent radicals.

\medskip

The central strategy exploits the fact that there exist only finitely many \( H_k \)-conjugacy classes of \( \theta \)-stable maximal \( k \)-split tori in \( G \), and in particular, finitely many $H_k$-conjugacy classes of maximal \( (\theta, k) \)-split tori~\cite{HW93}. Recall that a \( k \)-torus \( A \subset G \) is said to be \((\theta, k)\)-split if it is \( k \)-split and satisfies \( \theta(x) = x^{-1} \) for all \( x \in A \).

Let \( \{A_i \mid 1 \leq i \leq m\} \) be a set of representatives for the \( H_k \)-conjugacy classes of maximal \( (\theta, k) \)-split tori in \( G \), and define
\[
\mathcal{B} := \bigcup_{i=1}^{m} A_i, \quad \mathcal{B}_k := \bigcup_{i=1}^{m} A_i(k),
\]
where \( A_i(k) \) denotes the group of \( k \)-rational points of \( A_i \). By~\cite{BeOh}, we obtain the following polar decomposition: $G_k = K \mathcal{B}_k H_k$, where \( K \subset G_k \) is a compact subset. 

This decomposition shows that, up to \( G_k \)-conjugacy, it suffices to compute nontrivial Chabauty limits of \( H_k \) by conjugating \( H_k \) with sequences drawn from one of the maximal \( (\theta, k) \)-split tori \( A_i(k) \). The main result of this article, informally stated, is as follows:

\begin{theorem}[See Theorem \ref{thm::main_thm}]
\label{thm::main_thm_intro}
Let \( A_k \) denote the set of \( k \)-points of a maximal \( (\theta, k) \)-split torus \( A \subset G \), and let \( \{a_\ell\}_{\ell \in \mathbb{N}} \subset A_k \) be a sequence of hyperbolic elements such that the sequence of conjugates \( \{a_\ell H_k a_\ell^{-1}\}_{\ell \in \mathbb{N}} \) admits a nontrivial Chabauty limit \( L \). Then there exists a parabolic \( k \)-subgroup \( P_{\sigma_+} \subset G \) such that:
\begin{enumerate}
    \item \( L \subseteq U_{\sigma_+}(k) \rtimes \left( \Ker(\alpha)^0 \cdot \left( H_k \cap M_{\sigma_-,\sigma_+} \right) \right) \subseteq P_{\sigma_+} = U_{\sigma_+} \rtimes M_{\sigma_-,\sigma_+} \),  
    where:
    \begin{itemize}
        \item \( \alpha \) is the map defined in Proposition~\ref{prop::res_building_str_tran},
        \item \( \Ker(\alpha)^0 := \{ g \in \Ker(\alpha) \mid g \text{ is elliptic} \} \),
        \item \( M_{\sigma_-,\sigma_+} := P_{\sigma_+} \cap P_{\sigma_-} \) is a Levi factor of \( P_{\sigma_+} \),
        \item \( U_{\sigma_+} \) is the unipotent radical of \( P_{\sigma_+} \).
    \end{itemize}
    
    \item \( H_k \cap M_{\sigma_-,\sigma_+} \subseteq L \).
    
    \item The group \( L \) acts transitively on the set \( \Opp(\sigma_+) \), consisting of all ideal simplices in the Tits building of \( G_k \) that are opposite to \( \sigma_+ \).
    
    \item For every nontrivial unipotent element \( u \in U_{\sigma_+}(k) \), there exists \( m \in \Ker(\alpha)^0 \) such that \( um \in L \).
\end{enumerate}
\end{theorem}

The proof of Theorem~\ref{thm::main_thm_intro} unfolds through a sequence of intermediate results. 

The first is Theorem~\ref{thm::chabauty_in_parabolic}, which asserts that any nontrivial Chabauty limit \( L \) of \( H_k \), arising from conjugation by hyperbolic elements in a \( (\theta,k) \)-split torus \( A_k \subset G \), is contained in a parabolic subgroup \( P_{\sigma_+}(k) \). The proof distinguishes between the elliptic and hyperbolic elements of \( L \), and relies on several technical lemmas from Section~\ref{sec::some_useful_lemmas}, involving projections, midpoints,  fixed points, and convergence in \(\mathrm{CAT}(0)\) spaces, with particular emphasis on the geometry of Bruhat--Tits buildings. 

The second intermediate result investigates the structure of the elements in the Chabauty limit group \( L \) that lie in the Levi factor \( M_{\sigma_-, \sigma_+} \) of the parabolic subgroup \( P_{\sigma_+} \). A key observation is that \( M_{\sigma_-, \sigma_+}(k) \subseteq Z_{G_k}(A) \), the centralizer of the maximal \( (\theta,k) \)-split torus \( A \). In particular, this easily implies that \( H_k \cap M_{\sigma_-, \sigma_+} \subseteq L \). 

Moreover, Corollary~\ref{cor::find_elements_H} shows that for every nontrivial element \( l \in L \), there exists an element \( h \in H_k \cap M_{\sigma_-, \sigma_+} \) such that \( h \in L \) and \( lh^{-1} \in L \) is elliptic in \( G_k \). This, in turn, implies that any hyperbolic element of \( L \) arises from a hyperbolic element in \( H_k \cap M_{\sigma_-, \sigma_+} \); no additional hyperbolic elements are obtained beyond those.

Elliptic elements of \( L \) that lie in \( M_{\sigma_-, \sigma_+} \) are studied in Proposition~\ref{prop::levi_factors}. This result states that if \( l \in L \cap M_{\sigma_-, \sigma_+} \), then
\[
l \in \Ker(\alpha)^0 \cdot \left( H_k \cap M_{\sigma_-, \sigma_+} \right),
\]
where \( \alpha \) is the map defined in Proposition~\ref{prop::res_building_str_tran}, and \( \Ker(\alpha)^0 := \left\{ g \in \Ker(\alpha) \;\middle|\; g \text{ is elliptic} \right\} \).  One of the novel contributions here is the detailed analysis of the action of elements \( l \in L \cap M_{\sigma_-, \sigma_+} \) on the affine subbuilding associated with the Levi factor \( M_{\sigma_-, \sigma_+} \). This perspective leads to the appearance of the subgroup \( \Ker(\alpha) \).

The third set of results concerns the construction of unipotent elements in the limit group \( L \). Building on \cite{Cio24}, which investigates the dynamics of hyperbolic elements with respect to the cone topology on the Tits building of \( G_k \), Proposition~\ref{prop::exist_l_sigma} shows that the set of elliptic elements in the limit group \( L \leq P_{\sigma_+}(k) \) acts transitively on the set \( \Opp(\sigma_+) \) of all ideal simplices in the Tits building of \( G_k \) that are opposite to the ideal simplex \( \sigma_+ \).

Moreover, under the assumption that \( G_k \) satisfies the Moufang property, Proposition~\ref{prop::find_unipotent} establishes that every nontrivial unipotent element \( u \in U_{\sigma_+}(k) \) can be realized as part of an element \( um \in L \), with \( m \in \Ker(\alpha)^0 \). This construction again relies on a detailed analysis of the action of elements \( l \in L \) on the affine subbuilding associated with the Levi factor \( M_{\sigma_-, \sigma_+} \), as well as on the fact that, for Moufang groups \( G_k \), the unipotent radical \( U_{\sigma_+}(k) \) acts simply transitively on \( \Opp(\sigma_+) \).

\textbf{Appendix A} presents key results and more detailed proofs from \cite{HW93}, together with a geometric interpretation of the polar decomposition \( K \mathcal{B}_k H_k \) of \( G_k \), as formulated in \cite[Theorem 1.1]{BeOh}. \textbf{Appendix B} recalls the Inverse Function Theorem over non-Archimedean fields and explores its implications for openness and transitivity in the study of Chabauty limits of \( H_k \).

\subsection*{Acknowledgements}  The author was supported by a research grant (VIL53023) from VILLUM FONDEN.

\section{Background material}
We review essential background on the Chabauty topology, which compactifies spaces of closed subgroups, and on Bruhat--Tits buildings associated to reductive groups $G$ over non-Archimedean local fields $k$. Special attention is given to the fixed point sets under involutive $k$-automorphisms of $G$.

\subsection{Short on Chabauty topology}
\label{subsec::short_Chabauty}

For a comprehensive introduction to the Chabauty topology, we refer the reader to~\cite{CoPau, Harpe, GJT} and~\cite[Section 2]{Haettel}, along with the references therein.

In the context of a locally compact topological space $X$, the set $\mathcal{F}(X)$, which consists of all closed subsets of $X$, forms a compact topological space under the Chabauty topology (\cite[Proposition~1.7, p.~58]{CoPau}). When considering a family $\mathcal{T} \subset \mathcal{F}(X)$, a natural question arises: what is the closure $\overline{\mathcal{T}}^{Ch}$ of $\mathcal{T}$ with respect to the Chabauty topology, and do the elements of $\overline{\mathcal{T}}^{Ch} \setminus \mathcal{T}$ retain the same properties as those in $\mathcal{T}$? We refer to the elements of $\overline{\mathcal{T}}^{Ch}$ as the \textbf{Chabauty limits of $\mathcal{T}$}, and $\overline{\mathcal{T}}^{Ch}$ as the  \textbf{Chabauty compactification} of $\mathcal{T}$ (see \cite{Ch}).

For a locally compact group $G$, let $\mathcal{S}(G)$ denote the set of all closed subgroups of $G$. According to \cite[Proposition~1.7, p.~58]{CoPau}, the space $\mathcal{S}(G)$ is closed in $\mathcal{F}(G)$ with respect to the Chabauty topology, and is therefore compact.

\begin{proposition}(\cite[Proposition~1.8, p.~60]{CoPau}, \cite[Proposition I.3.1.3]{CEM})
\label{prop::chabauty_conv}
 Suppose $X$ is a locally compact metric space.
A sequence of closed subsets $\{F_n\}_{n \in \NN} \subset \mathcal{F}(X)$  converges to $F \in \mathcal{F}(X)$ if and only if the following two conditions are satisfied:
\begin{itemize} 
\item[1)] For every $f \in F$ there is a sequence $\{f_n \in F_n\}_{n \in \NN}$ converging to $f$;
\item[2)] For every sequence $\{f_n \in F_n\}_{n \in \NN}$, if there is a strictly increasing subsequence $\{n_k\}_{k \in \NN}$ such that $\{f_{n_k} \in F_{n_k}\}_{k \in \NN}$ converges to $f$, then $f \in F$.
\end{itemize}
\end{proposition}

Proposition~\ref{prop::chabauty_conv} also applies to sequences of closed subgroups $\{H_n\}_{n \in \NN} \subset \mathcal{S}(G)$ converging to $H \in \mathcal{S}(G)$, providing a similar characterization of convergence in $\mathcal{S}(G)$, when $G$ is a locally compact group equipped with a metric inducing its topology.

Let \( H \) be a closed subgroup of \( G \), so that \( H \in \mathcal{S}(G) \). We say that \( L \in \mathcal{S}(G) \) is a \textbf{Chabauty limit of \( H \) in \( G \)} if there exists a sequence \( \{g_n\}_{n \in \mathbb{N}} \) of elements in \( G \) such that the sequence of subgroups \( \{g_n H g_n^{-1}\}_{n \in \mathbb{N}} \) converges to \( L \) with respect to the Chabauty topology on \( \mathcal{S}(G) \).

% A locally compact topological group $G$ is metrizable with a left  invariant metric all of whose spheres are bounded, if and only if $G$ is second countable.

\subsection{The Bruhat--Tits building and the fixed point group of an involution}
\label{section::auto_invol_building}

Let \( G \) be a connected linear reductive group defined over a non-Archimedean local field \( k \), e.g., a finite extension of the \( p \)-adic field \( \mathbb{Q}_p \). It is a well-known fact that all maximal \( k \)-split tori of \( G \) are \( G_k \)-conjugate (see \cite[20.9(ii)]{Borel}). To each maximal \( k \)-split torus of \( G \), one associates a copy of the affine space \( \mathbb{R}^n \), tessellated by copies of a specific \( n \)-dimensional simplex called a chamber. Here, \( n \) denotes the \( k \)-rank of the derived subgroup of \( G \).  According to the theory developed by Bruhat–Tits \cite{BrTi_72,BrTi_84}, these copies of \( \mathbb{R}^n \) are glued together in a coherent way to form a simplicial complex \( \Delta_{G_k} \), on which the group of \( k \)-points \( G_k \) of \( G \) acts by automorphisms that preserve the simplicial structure of \( \Delta_{G_k} \).  This complex \( \Delta_{G_k} \), known as the Bruhat–Tits building of \( G_k \), is a union of a system of simplicial subcomplexes called apartments (each isomorphic to \( \mathbb{R}^n \)) that satisfy three axioms. The construction ensures that the system of apartments in \( \Delta_{G_k} \) is in bijection with the set of maximal \( k \)-split tori of \( G \).  Moreover, by \cite[2.8.4]{BrTi_84}, since the field \( k \) is complete, the system of apartments in the Bruhat–Tits building \( \Delta_{G_k} \) is also complete, in the sense that any other possible apartment of \( \Delta_{G_k} \) belongs to the considered system.  One of the key properties is that \( G_k \) acts strongly transitively on its building \( \Delta_{G_k} \); that is, \( G_k \) acts transitively on the set of all pairs \( (c, A) \), where \( c \) is a chamber and \( A \) is an apartment of \( \Delta_{G_k} \) containing \( c \). The associated Bruhat--Tits building of $G_k$ is not irreducible in general.

\medskip
\textbf{Standing Assumption}: For the rest of this article, we assume that \( G \) is a connected linear reductive group defined over a non-Archimedean local field \( k \) of \textbf{characteristic zero}, and that its associated Bruhat--Tits building \( \Delta_{G_k} \) is \textbf{irreducible}. Consequently, the spherical Tits buildings associated with \( G_k \) exhibits irreducibility as well.

\medskip
The construction of $\Delta_{G_k}$ is quite involved, and made in three steps. First, one needs to construct $\Delta_{G_k}$ for a $k$-split $G$. The next step, called quasi-split descent, and making use of the first step, constructs $\Delta_{G_k}$ for a quasi-split $G$. The third step is to descent from the quasi-split $G$ to a general $G$ and build $\Delta_{G_k}$, this step is called \'{e}tale descent. A different approach is taken in the book \cite{KaPra} by Kaletha--Prasad.

Consider \( \theta \) an involutive \( k \)-automorphism of \( G \), that is, \( \theta^2 = \mathrm{id} \) and \( \theta \) is defined over \( k \). Let \( G^{\theta} \leq G \) be the fixed-point subgroup of \( \theta \):
\[
G^{\theta} := \{ h \in G \mid \theta(h) = h \}
\]
and $H = (G^{\theta})^{o}$ be the connected component of $G^{\theta}$. Then $H$ is defined over $k$, since $\theta$ is defined over $k$ (see \cite[Proposition 1.6]{HW93}). Moreover, by \cite[Theorem 2.1]{PraYu}, $H$ is reductive.

One can observe that \( \theta \) induces an involutive automorphism of the Bruhat–Tits building \( \Delta_{G_k} \), which we will also denote by \( \theta \), see \cite[Section 4.2]{BeOh}, or \cite[Section 2.4]{Pra1}, or \cite[Section 4.1]{WangZou2024}. In principle, such an automorphism \( \theta \) of the building \( \Delta_{G_k} \) may not be unique; see \cite[Chapter 14, pp. 490–491]{KaPra}.  More precisely, since \( \theta \) is a \( k \)-automorphism of \( G \), it maps any maximal \( k \)-split torus of \( G \) to another such torus. Therefore, the set of maximal \( k \)-split tori in \( G \) is preserved under the action of \( \theta \).  Moreover, by \cite[Proposition 2.3]{HW93}, there exists a maximal \( k \)-split torus \( T \subset G \) that is \( \theta \)-stable, i.e., \( \theta(T) = T \). This allows us to define a natural action of \( \theta \) on the spherical root system of \( G \) with respect to the \( \theta \)-stable torus \( T \). One can observe that this induced action stabilizes that root system.  Consequently, \( \theta \) induces a natural action on the Coxeter complex associated with the finite Weyl group of \( (G, T) \), as well as on the affine apartment \( \mathcal{A}_T \subset \Delta_{G_k} \) corresponding to the torus \( T \). Since \( \theta \) stabilizes the spherical root system, it follows that \( \theta(\mathcal{A}_T) = \mathcal{A}_T \), and \( \theta \) maps chambers to chambers, although it may not preserve vertex types.  Finally, because \( \theta \) maps any maximal \( k \)-split torus of \( G \) to another, the induced action of \( \theta \) on \( \Delta_{G_k} \) sends apartments to apartments. Thus, \( \theta \) defines a natural plurisimplicial action on the building \( \Delta_{G_k} \) in the following natural way:
\[
\text{For } g \in G_k \text{ and } x \in \Delta_{G_k}, \text{ we have } \theta(g(x)) = \theta(g)(\theta(x)).
\]

%Alternatively, consider the semidirect product \( G' := \langle \theta \rangle \ltimes G \), where \( \langle \theta \rangle \) is the group of order $2$ generated by \( \theta \). The group \( G' \) can be regarded as a reductive group defined over \( k \) -- it is not reductive in the classical linear algebraic sense, but it behaves reductively in the context of building theory.  By the observation above, the set of maximal \( k \)-split tori of \( G' \) coincides with that of \( G \).  Since the Bruhat--Tits building is constructed from the set of all maximal \( k \)-split tori, it follows that \( \Delta_{G'_k} \) and \( \Delta_{G_k} \) are identical. 

Notice that the same argument as above applies if we consider \( \theta \) to be a \( k \)-automorphism of \( G \) of finite order.

Denote by \( \Delta^{\theta}_{G_k} \) the set of points in \( \Delta_{G_k} \) fixed by the involution \( \theta \):
\[
\Delta^{\theta}_{G_k} := \{ x \in \Delta_{G_k} \mid \theta(x) = x \}.
\]

\begin{remark}
\label{rem::invariance_H}
One can observe that \( \Delta^{\theta}_{G_k} \) is invariant under the action of the \( k \)-points \( H_k \) of \( H \):
\[
H_k(\Delta^{\theta}_{G_k}) = \Delta^{\theta}_{G_k}.
\]
Indeed, by the definition of the action of \( \theta \) on both \( G \) and \( \Delta_{G_k} \), we have
\[
\theta(h(x)) = \theta(h)(\theta(x)) = h(x),
\]
for every \( h \in H_k \) and every \( x \in \Delta^{\theta}_{G_k} \). Hence, \( h(x) \in \Delta^{\theta}_{G_k} \).
\end{remark}

Since \( \Delta_{G_k} \) is a \(\mathrm{CAT}(0)\) space -- as is every affine building -- it is uniquely geodesic. As \( \theta \) is an automorphism of \( \Delta_{G_k} \), for any points \( x, y \in \Delta^{\theta}_{G_k} \), the unique geodesic between \( x \) and \( y \) in \( \Delta_{G_k} \) must lie entirely within \( \Delta^{\theta}_{G_k} \). Therefore, \( \Delta^{\theta}_{G_k} \) is a geodesically convex subset of \( \Delta_{G_k} \). By a similar argument, it is also easy to see that \( \Delta^{\theta}_{G_k} \) is a closed subset of \( \Delta_{G_k} \).

In many cases, one can show that \( \Delta^{\theta}_{G_k} \) is the Bruhat--Tits building of the \( k \)-points \( H_k \). For a clear exposition of this result, the reader is referred to the introduction of \cite{PraYu}.  In principle, certain conditions must be imposed on the non-Archimedean local field \( k \). The most general setting is when \( k \) is a discretely valued field with a Henselian valuation ring and a separably closed (but not necessarily perfect) residue field of characteristic \( p \), such that \( p \) does not divide the order of the \( k \)-automorphism \( \theta \) of \( G \), where \( \theta \) is assumed to be of finite order (see \cite{Pra1} and \cite[Chapter 12]{KaPra}). According to \cite[Example 2.1.4]{KaPra} or \cite[Proposition 2.A.5]{Ach}, any field \( k \) equipped with a valuation \( \nu \) and complete with respect to the topology induced by \( \nu \) is Henselian. In particular, since we are working with non-Archimedean local fields of characteristic zero -- which are complete by definition -- this setting ensures that we are in the setting of Henselian fields. 

Nevertheless, in the remainder of this article, we do not explicitly rely on the fact that \( \Delta^{\theta}_{G_k} \) coincides with the Bruhat--Tits building of \( H_k \).

\section{Some useful lemmas}
\label{sec::some_useful_lemmas}

We establish several technical lemmas concerning projections, fixed points, and convergence in \cat spaces, with a particular focus on the geometry of the Bruhat--Tits buildings. These results form the foundation for the geometric analysis of Chabauty limits developed in the subsequent sections.

In this section, we adopt assumptions that reflect the fundamental properties of maximal \( (\theta, k) \)-split tori in the group \( G \). Specifically, we work under the following hypothesis:

\vspace{1em}
\noindent\textbf{(Hyp(\( \theta,k \))-split):}
Let \( A \leq \Aut(\Delta_{G_k}) \) such that:
\begin{enumerate}
  \item \( A \) is abelian,
  \item For all \( a \in A \), we have \( \theta(a) = a^{-1} \), i.e., \( A \) is \( \theta \)-split,
  \item There exists an apartment \( \mathcal{A}_{\theta,k} \subset \Delta_{G_k} \) such that \( \mathcal{A}_{\theta,k} \subseteq \Min(a) \) and \( a(\mathcal{A}_{\theta,k}) = \mathcal{A}_{\theta,k} \) for all \( a \in A \).
\end{enumerate}

Recall that in the context of \cat spaces \( X \) (see \cite[Chapter II.6]{BH99}), such as the affine building \( \Delta_{G_k} \), the minimal set of an isometry \( g \in \Isom(X) \) is defined by:
\[
\mathrm{Min}(g) := \{ x \in X \mid d(x, g(x)) = |g| \},
\]
where \( |g| \) denotes the \textbf{translation length} of \( g \), given by
\[
|g| := \inf \{ d(x, g(x)) \mid x \in X \}.
\]
Since affine buildings are simplicial complexes and thus behave as discrete geometric objects, the infimum \( |g| \) is always attained for any \( g \in \Aut(\Delta_{G_k}) \). Moreover, elements of \( \Aut(\Delta_{G_k}) \) are either elliptic or hyperbolic.

We say that a sequence $\{a_\ell\}_{\ell \in \NN} \subset \Aut(\Delta_{G_k})$ \textbf{diverges to infinity} if it eventually escapes every compact subset of $\Aut(\Delta_{G_k})$, where the group is endowed with the compact-open topology.

Since the space $\Delta_{G_k} \cup \partial \Delta_{G_k}$ is compact with respect to the cone topology (see \cite[Part II, Theorem 8.13 and Exercise 8.15(2)]{BH99}), the sequence $\{a_\ell\}_{\ell \geq 1}$ admits a subsequence $\{a_{\ell_j}\}_{j \geq 1}$ such that the sequence $\{a_{\ell_j}(x)\}_{j \geq 1}$ converges to an ideal point $\xi \in \partial \Delta_{G_k}$, for some (and any) point $x \in \Delta_{G_k}$.  Note that there exists a unique minimal ideal simplex $\sigma \subset \partial \Delta_{G_k}$ such that $\xi$ lies in the interior of $\sigma$. Depending on the case, $\sigma$ may be an ideal chamber or an ideal vertex of the spherical building $\partial \Delta_{G_k}$.

\medskip
Let $a \in \Aut(\Delta_{G_k})$ be a hyperbolic element such that its minimal set $\Min(a)$ contains an apartment $\mathcal{A}$ of $\Delta_{G_k}$ with $a(\mathcal{A}) = \mathcal{A}$. Then $\mathcal{A}$ contains a translation axis of $a$, whose endpoints $ \xi_{a+}, \xi_{a-} \in \partial \mathcal{A} $
are called the \textbf{attracting} and \textbf{repelling} endpoints of $a$, respectively.  These endpoints lie in the ideal boundary $\partial \mathcal{A}$ of the apartment $\mathcal{A}$.  In particular, there exist two unique, minimal, opposite ideal simplices $\sigma_{a+}, \sigma_{a-} \subset \partial \mathcal{A}$ such that $\xi_{a+}$ lies in the interior of $\sigma_{a+}$ and $\xi_{a-}$ lies in the interior of $\sigma_{a-}$.

\vspace{1em}
Let $A \leq \operatorname{Aut}(\Delta_{G_k})$ be a subgroup satisfying the condition \textbf{(Hyp$(\theta,k)$-split)}. Based on the above observations and notation, and to simplify the exposition, we introduce the following hypothesis concerning a sequence of elements in $A$:

\vspace{1em}
\noindent\textbf{(Same $\sigma_\pm$):}  
Let $\{a_\ell\}_{\ell \in \mathbb{N}} \subset A$ be a sequence of hyperbolic elements with:
\begin{enumerate}
  \item $\displaystyle \lim_{\ell \to \infty} |a_\ell| = \infty$,
  \item There exist ideal simplices $\sigma_\pm \subset \partial \mathcal{A}_{\theta,k}$ such that $\sigma_\pm = \sigma_{a_\ell \pm}$ for every $\ell \in \mathbb{N}$,
  \item There exist points $\xi_+, \xi_- \in \partial \mathcal{A}_{\theta,k}$ such that, for some (and hence any) point $x \in \mathcal{A}_{\theta,k}$, we have
  $\lim\limits_{\ell \to \infty} a_\ell(x) = \xi_+$ and $\lim\limits_{\ell \to \infty} a_\ell^{-1}(x) = \xi_-,$ where the limits are taken with respect to the cone topology on $\Delta_{G_k} \cup \partial \Delta_{G_k}$. In particular, this implies $\lim\limits_{\ell \to \infty} \xi_{a_\ell \pm} = \xi_\pm$, and each point $\xi_\pm$ lies in, possibly on the boundary of, the corresponding ideal simplex $\sigma_\pm$.
\end{enumerate}

\begin{lemma}
\label{lem::projection_apartment}
Let \( \mathcal{A} \) be an apartment of the Bruhat–Tits building \( \Delta_{G_k} \) that is \( \theta \)-stable, i.e., \( \theta(\mathcal{A}) = \mathcal{A} \). Then there exists at least one point in \( \mathcal{A} \) that is fixed by \( \theta \).
\end{lemma}

\begin{proof}
Recall that \( \Delta_{G_k}^{\theta} \) denotes the fixed-point set of \( \theta \), which acts as an automorphism on \( \Delta_{G_k} \). Moreover, \( \Delta_{G_k}^{\theta} \) is a geodesically convex and closed subset of \( \Delta_{G_k} \). Suppose, for contradiction, that \( \mathcal{A} \cap \Delta_{G_k}^{\theta} = \emptyset \).
Consider the projection of \( \mathcal{A} \) onto \( \Delta_{G_k}^{\theta} \), defined as the set of points in \( \Delta_{G_k}^{\theta} \) closest to \( \mathcal{A} \). Since both sets are closed and \( \Delta_{G_k} \) is a \(\mathrm{CAT}(0)\) space, this minimal distance is achieved for at least one pair of points \( x \in \mathcal{A} \) and \( y \in \Delta_{G_k}^{\theta} \), such that
\[
\mathrm{dist}_{\Delta_{G_k}}(\mathcal{A}, \Delta_{G_k}^{\theta}) = \mathrm{dist}_{\Delta_{G_k}}(x, y).
\]
Note that \( \theta(y) = y \), and that \( y \) is the unique projection of \( x \) onto \( \Delta_{G_k}^{\theta} \), due to the uniqueness of projections in \(\mathrm{CAT}(0)\) geometry. Furthermore, since \( \theta \) is an automorphism of \( \Delta_{G_k} \), it preserves distances. Thus,
\[
\mathrm{dist}_{\Delta_{G_k}}(x, y) = \mathrm{dist}_{\Delta_{G_k}}(\theta(x), \theta(y)) = \mathrm{dist}_{\Delta_{G_k}}(\theta(x), y),
\]
which implies that \( y \) is also the projection of \( \theta(x) \) onto \( \Delta_{G_k}^{\theta} \).

Since \( \mathcal{A} \) is geodesically convex and closed, the projection of \( y \) onto \( \mathcal{A} \) is unique and must be \( x \). But \( \theta(x) \) satisfies the same minimality condition, so it must also be the projection of \( y \) onto \( \mathcal{A} \), implying \( \theta(x) = x \). Therefore, \( x \in \Delta_{G_k}^{\theta} \), contradicting our assumption that \( \mathcal{A} \cap \Delta_{G_k}^{\theta} = \emptyset \).  Hence, \( \mathcal{A} \cap \Delta_{G_k}^{\theta} \neq \emptyset \), and the lemma is proven.
\end{proof}

It may happen that, for an apartment \( \mathcal{A} \) of \( \Delta_{G_k} \) that is \( \theta \)-stable, the nontrivial intersection \( \mathcal{A} \cap \Delta_{G_k}^{\theta} \) is either compact or noncompact, but in both cases it is a geodesically convex and closed subset.

\begin{lemma}
\label{lem::projection_point_building}
For any point \( x \in \Delta_{G_k} \), the projection \( \mathrm{proj}_{\Delta_{G_k}^{\theta}}(x) \) of \( x \) onto \( \Delta_{G_k}^{\theta} \) is the midpoint \( \mathrm{mid}(x, \theta(x)) \in \Delta_{G_k} \) between \( x \) and \( \theta(x) \). More precisely,
\[
\mathrm{proj}_{\Delta_{G_k}^{\theta}}(x) = \mathrm{mid}(x, \theta(x)) = \mathrm{proj}_{\Delta_{G_k}^{\theta}}(\theta(x)).
\]
\end{lemma}
\begin{proof}
Let \( x \in \Delta_{G_k} \). Since \( \Delta_{G_k} \) is a \(\mathrm{CAT}(0)\) space and \( \Delta_{G_k}^{\theta} \) is a geodesically convex subset, the projection \( \mathrm{proj}_{\Delta_{G_k}^{\theta}}(y) \) of any point \( y \in \Delta_{G_k} \) onto \( \Delta_{G_k}^{\theta} \) is unique. In particular, this holds for \( \mathrm{proj}_{\Delta_{G_k}^{\theta}}(x) \).

Because \( \theta \) fixes \( \mathrm{proj}_{\Delta_{G_k}^{\theta}}(x) \in \Delta_{G_k}^{\theta} \), and \( \theta \) is an automorphism of \( \Delta_{G_k} \), we have:
\[
\mathrm{dist}_{\Delta_{G_k}}(x, \mathrm{proj}_{\Delta_{G_k}^{\theta}}(x)) = \mathrm{dist}_{\Delta_{G_k}}(\theta(x), \mathrm{proj}_{\Delta_{G_k}^{\theta}}(x)).
\]
Moreover, it is easy to see that the projection of \( \theta(x) \) onto \( \Delta_{G_k}^{\theta} \) is exactly \( \mathrm{proj}_{\Delta_{G_k}^{\theta}}(x) \).

Now, since geodesics in a \(\mathrm{CAT}(0)\) space are unique, and the geodesic segment \( [x, \theta(x)] \) is \( \theta \)-stable, the midpoint \( \mathrm{mid}(x, \theta(x)) \in \Delta_{G_k} \) between \( x \) and \( \theta(x) \) is also unique and must be fixed by \( \theta \). This implies that \( \mathrm{mid}(x, \theta(x)) \in \Delta_{G_k}^{\theta} \).

Suppose, for contradiction, that \( \mathrm{proj}_{\Delta_{G_k}^{\theta}}(x) \neq \mathrm{mid}(x, \theta(x)) \). Consider the geodesic triangle in \( \Delta_{G_k} \) with vertices \( x \), \( \theta(x) \), and \( \mathrm{proj}_{\Delta_{G_k}^{\theta}}(x) \), and its comparison triangle in the Euclidean plane with vertices \( \overline{x} \), \( \overline{\theta(x)} \), and \( \overline{\mathrm{proj}_{\Delta_{G_k}^{\theta}}(x)} \). In this comparison triangle, the segments \( [\overline{x}, \overline{\mathrm{proj}_{\Delta_{G_k}^{\theta}}(x)}] \) and \( [\overline{\theta(x)}, \overline{\mathrm{proj}_{\Delta_{G_k}^{\theta}}(x)}] \) are equal in length, and by the definition of comparison triangles, we have:
\[
\mathrm{mid}(\overline{x}, \overline{\theta(x)}) = \overline{\mathrm{mid}(x, \theta(x))}.
\]
Thus, the segment \( [\overline{\mathrm{proj}_{\Delta_{G_k}^{\theta}}(x)}, \mathrm{mid}(\overline{x}, \overline{\theta(x)})] \) is perpendicular to \( [\overline{x}, \overline{\theta(x)}] \) in the Euclidean plane, and therefore:
\[
\mathrm{dist}_{\mathrm{Eucl}}(\overline{x}, \overline{\mathrm{mid}(x, \theta(x))}) < \mathrm{dist}_{\mathrm{Eucl}}(\overline{x}, \overline{\mathrm{proj}_{\Delta_{G_k}^{\theta}}(x)}) = \mathrm{dist}_{\Delta_{G_k}}(x, \mathrm{proj}_{\Delta_{G_k}^{\theta}}(x)).
\]
By the \(\mathrm{CAT}(0)\) condition, we then have:
\[
\mathrm{dist}_{\Delta_{G_k}}(x, \mathrm{mid}(x, \theta(x))) \leq \mathrm{dist}_{\mathrm{Eucl}}(\overline{x}, \overline{\mathrm{mid}(x, \theta(x))}) < \mathrm{dist}_{\Delta_{G_k}}(x, \mathrm{proj}_{\Delta_{G_k}^{\theta}}(x)),
\]
which contradicts the minimality of \( \mathrm{dist}_{\Delta_{G_k}}(x, \mathrm{proj}_{\Delta_{G_k}^{\theta}}(x)) \). Therefore, \( \mathrm{proj}_{\Delta_{G_k}^{\theta}}(x) = \mathrm{mid}(x, \theta(x)) \), and the lemma is proven.
\end{proof}

\begin{lemma}
\label{lem::projection_point_apartment}
Let $A \leq \operatorname{Aut}(\Delta_{G_k})$ be a subgroup satisfying the condition \textbf{(Hyp$(\theta,k)$-split)}. Suppose that $\mathcal{A}:= \mathcal{A}_{\theta,k}$  is \( \theta \)-stable, i.e., \( \theta(\mathcal{A}) = \mathcal{A} \). Let \( \{a_\ell\}_{\ell \geq 1} \subset A \) be a sequence of hyperbolic elements satisfying the condition \textbf{(Same $\sigma_\pm$)}. Then:
\begin{enumerate}
    \item For any point \( x \in \mathcal{A} \), the projection \( \mathrm{proj}_{\Delta_{G_k}^{\theta}}(x) \) of \( x \) onto \( \Delta_{G_k}^{\theta} \) is the midpoint \( \mathrm{mid}(x, \theta(x)) \in \mathcal{A} \) between \( x \) and \( \theta(x) \).
    
    \item For every point \( x \in \mathcal{A} \), the sequences of projections
    \[
    \left\{ \mathrm{proj}_{\Delta_{G_k}^{\theta}}(a_\ell(x)) \right\}_{\ell \geq 1} \subset \Delta_{G_k}^{\theta} \cap \mathcal{A}, \quad \text{and} \quad \left\{ \mathrm{proj}_{\Delta_{G_k}^{\theta}}(a_\ell^{-1}(x)) \right\}_{\ell \geq 1} \subset \Delta_{G_k}^{\theta} \cap \mathcal{A}
    \]
    remain at bounded \( \Delta_{G_k} \)-distance from \( \mathrm{mid}(x, \theta(x)) \). Hence, they form bounded subsets of \( \Delta_{G_k}^{\theta} \cap \mathcal{A} \), where the bound depends on \( x \).
\end{enumerate}
\end{lemma}
\begin{proof}
We retain the notation from the conditions \textbf{(Hyp$(\theta,k)$-split)} and \textbf{(Same $\sigma_\pm$)}.
\begin{enumerate}
\item
This follows directly from Lemma~\ref{lem::projection_point_building} as a special case.

\item
Let \( x \in \mathcal{A} \).  
It suffices to prove the second part of the lemma for the sequence of projections \( \{\mathrm{proj}_{\Delta_{G_k}^{\theta}}(a_\ell(x))\}_{\ell \geq 1} \subset \Delta_{G_k}^{\theta} \). Observe that
\[
\theta(a_\ell(x)) = \theta(a_\ell)(\theta(x)) = a_\ell^{-1}(\theta(x)),
\]
and the sequences \( \{a_\ell(x)\}_{\ell \geq 1} \subset \mathcal{A} \) and \( \{a_\ell^{-1}(\theta(x))\}_{\ell \geq 1} \subset \mathcal{A} \) converge respectively to the boundary points \( \xi_+ \in \sigma_+ \subset \partial \mathcal{A} \) and \( \xi_- \in \sigma_- \subset \partial \mathcal{A} \).

Let \( \mathrm{mid}(x, \theta(x)) \in \mathcal{A} \) be the midpoint of the geodesic segment \( [x, \theta(x)] \subset \mathcal{A} \). As before, \( \mathrm{mid}(x, \theta(x)) \) is fixed by \( \theta \), hence \( \mathrm{mid}(x, \theta(x)) \in \mathcal{A} \cap \Delta_{G_k}^{\theta} \). Moreover, the sequences \( \{a_\ell(\mathrm{mid}(x, \theta(x)))\}_{\ell \geq 1} \) and \( \{a_\ell^{-1}(\mathrm{mid}(x, \theta(x)))\}_{\ell \geq 1} \) converge respectively to \( \xi_+ \) and \( \xi_- \). Applying $\theta$, one obtains that $\theta(\xi_\pm)= \xi_\mp$.

Consider the bi-infinite geodesic in \( \mathcal{A} \) passing through \( \mathrm{mid}(x, \theta(x)) \) with endpoints \( \xi_\pm \). Since \( \theta(\xi_\pm) = \xi_\mp \), the geodesic \( [\xi_-, \mathrm{mid}(x, \theta(x)), \xi_+] \subset \mathcal{A} \cup \partial \mathcal{A} \) is \( \theta \)-stable. Furthermore, the sequences of geodesic segments
\[
\{[a_\ell(\mathrm{mid}(x, \theta(x))), \mathrm{mid}(x, \theta(x))]\}_{\ell \geq 1}, \quad \{[a_\ell^{-1}(\mathrm{mid}(x, \theta(x))), \mathrm{mid}(x, \theta(x))]\}_{\ell \geq 1}
\]
converge in the cone topology of \( \mathcal{A} \cup \partial \mathcal{A} \) to the geodesic rays \( [\xi_+, \mathrm{mid}(x, \theta(x))] \) and \( [\xi_-, \mathrm{mid}(x, \theta(x))] \), respectively.

In particular, as the geodesic segment \( [a_\ell^{-1}(\mathrm{mid}(x, \theta(x))), a_\ell(\mathrm{mid}(x, \theta(x)))] \) becomes asymptotically parallel to the bi-infinite geodesic \( [\xi_-, \mathrm{mid}(x, \theta(x)), \xi_+] \) as \( \ell \to \infty \). Moreover, 
\begin{equation*}
\begin{split}
\mathrm{dist}_{\mathcal{A}}(a_\ell(\mathrm{mid}(x, \theta(x))), \mathrm{mid}(x, \theta(x))) 
&= \mathrm{dist}_{\mathcal{A}}(\theta(a_\ell(\mathrm{mid}(x, \theta(x)))), \theta(\mathrm{mid}(x, \theta(x)))) \\
&= \mathrm{dist}_{\mathcal{A}}(a_\ell^{-1}(\mathrm{mid}(x, \theta(x))), \mathrm{mid}(x, \theta(x))),
\end{split}
\end{equation*}
and a simple Euclidean geometric argument using isosceles triangles shows that the midpoint of the segment
\[
[a_\ell^{-1}(\mathrm{mid}(x, \theta(x))), a_\ell(\mathrm{mid}(x, \theta(x)))]
\]
-- which, by Lemma~\ref{lem::projection_point_building}, is exactly \( \mathrm{proj}_{\Delta_{G_k}^{\theta}}(a_\ell(\mathrm{mid}(x, \theta(x)))) \) -- remains at bounded distance from \( \mathrm{mid}(x, \theta(x)) \), more precisely,  \( \mathrm{proj}_{\Delta_{G_k}^{\theta}}(a_\ell(\mathrm{mid}(x, \theta(x)))) =  \mathrm{mid}(x, \theta(x)) \). 

Since \( a_\ell \in A \) and \( x \in \mathcal{A} \), we have \( a_\ell(x) \in \mathcal{A} \), and by Lemma~\ref{lem::projection_point_building}, \( \mathrm{proj}_{\Delta_{G_k}^{\theta}}(a_\ell(x)) \in \Delta_{G_k}^{\theta} \cap \mathcal{A} \). Moreover, because each \( a_\ell \) is an isometry of \( \Delta_{G_k} \), and projections onto geodesically convex, closed subsets of \(\mathrm{CAT}(0)\) spaces do not increase distances (see \cite[Chapter II.2, Prop. 2.4(4)]{BH99}), we have:
\begin{equation*}
\begin{split}
\mathrm{dist}_{\Delta_{G_k}}(\mathrm{proj}_{\Delta_{G_k}^{\theta}}(a_\ell(x)), \mathrm{mid}(x, \theta(x))) 
&\leq \mathrm{dist}_{\Delta_{G_k}}(\mathrm{proj}_{\Delta_{G_k}^{\theta}}(a_\ell(x)), \mathrm{proj}_{\Delta_{G_k}^{\theta}}(a_\ell(\mathrm{mid}(x, \theta(x))))) \\
&\quad + \mathrm{dist}_{\Delta_{G_k}}(\mathrm{proj}_{\Delta_{G_k}^{\theta}}(a_\ell(\mathrm{mid}(x, \theta(x)))), \mathrm{mid}(x, \theta(x))) \\
&\leq \mathrm{dist}_{\Delta_{G_k}}(a_\ell(x), a_\ell(\mathrm{mid}(x, \theta(x)))) + \text{const} \\
&= \mathrm{dist}_{\Delta_{G_k}}(x, \mathrm{mid}(x, \theta(x))) + \text{const} < \infty,
\end{split}
\end{equation*}
as desired. This completes the proof of the lemma.
\end{enumerate}
\end{proof}

The following lemma extends Lemma~\ref{lem::projection_point_apartment} by considering arbitrary points \( x \in \Delta_{G_k} \).

\begin{lemma}
\label{lem::projection_sequence}
Let \( A \leq \Aut(\Delta_{G_k}) \) be a subgroup satisfying the condition \textbf{(Hyp$(\theta,k)$-split)}. Suppose that the apartment \( \mathcal{A} := \mathcal{A}_{\theta,k} \) is \( \theta \)-stable, i.e., \( \theta(\mathcal{A}) = \mathcal{A} \). Let \( \{a_\ell\}_{\ell \geq 1} \subset A \) be a sequence of hyperbolic elements satisfying the condition \textbf{(Same $\sigma_\pm$)}.

Then, for each point \( x \in \mathcal{A} \), the sequences of projections
\[
\left\{ \mathrm{proj}_{\Delta_{G_k}^{\theta}}(a_\ell(x)) \right\}_{\ell \geq 1} \subset \Delta_{G_k}^{\theta} \quad \text{and} \quad \left\{ \mathrm{proj}_{\Delta_{G_k}^{\theta}}(a_\ell^{-1}(x)) \right\}_{\ell \geq 1} \subset \Delta_{G_k}^{\theta}
\]
form bounded subsets of \( \Delta_{G_k}^{\theta} \), where the bound depends on \( x \).
\end{lemma}

\begin{proof}
We retain the notation from the conditions \textbf{(Hyp$(\theta,k)$-split)} and \textbf{(Same $\sigma_\pm$)}.

Let \( x \in \Delta_{G_k} \), and let \( y := \mathrm{proj}_{\mathcal{A}}(x) \in \mathcal{A} \) be the projection of \( x \) onto the apartment \( \mathcal{A} \). Since \( \mathcal{A} \) is a geodesically convex and closed subset of \( \Delta_{G_k} \), and we are working in a \(\mathrm{CAT}(0)\) setting, this projection $y$ is well-defined and unique.

By the assumptions \textbf{(Hyp$(\theta,k)$-split)} and \textbf{(Same $\sigma_\pm$)}, the sequences \( \{a_\ell(x)\}_{\ell \geq 1} \) and \( \{a_\ell(y)\}_{\ell \geq 1} \) both converge to \( \xi_+ \in \sigma_+ \), and similarly, \( \{a_\ell^{-1}(x)\}_{\ell \geq 1} \) and \( \{a_\ell^{-1}(y)\}_{\ell \geq 1} \) converge to \( \xi_- \in \sigma_- \), with respect to the cone topology on \( \Delta_{G_k} \cup \partial \Delta_{G_k} \).

It suffices to prove the lemma for the sequence \( \{\mathrm{proj}_{\Delta_{G_k}^{\theta}}(a_\ell(x))\}_{\ell \geq 1} \).

Since each \( a_\ell \in \operatorname{Aut}(\Delta_{G_k}) \) is an isometry, we have
\[
\mathrm{dist}_{\Delta_{G_k}}(x, y) = \mathrm{dist}_{\Delta_{G_k}}(a_\ell(x), a_\ell(y)) \quad \text{for all } \ell \geq 1.
\]
Moreover, projections onto geodesically convex, closed subsets of \(\mathrm{CAT}(0)\) spaces do not increase distances (see \cite[Chapter II.2, Proposition 2.4(4)]{BH99}), so
\[
\mathrm{dist}_{\Delta_{G_k}}(\mathrm{proj}_{\Delta_{G_k}^{\theta}}(a_\ell(x)), \mathrm{proj}_{\Delta_{G_k}^{\theta}}(a_\ell(y))) \leq \mathrm{dist}_{\Delta_{G_k}}(x, y).
\]

Applying the triangle inequality and Lemma~\ref{lem::projection_point_apartment}.2 to the point \( y \in \mathcal{A} \), we obtain:
\begin{equation*}
\begin{split}
\mathrm{dist}_{\Delta_{G_k}}(\mathrm{proj}_{\Delta_{G_k}^{\theta}}(a_\ell(x)), \mathrm{mid}(y, \theta(y))) 
&\leq \mathrm{dist}_{\Delta_{G_k}}(\mathrm{proj}_{\Delta_{G_k}^{\theta}}(a_\ell(x)), \mathrm{proj}_{\Delta_{G_k}^{\theta}}(a_\ell(y))) \\
&\quad + \mathrm{dist}_{\Delta_{G_k}}(\mathrm{proj}_{\Delta_{G_k}^{\theta}}(a_\ell(y)), \mathrm{mid}(y, \theta(y))) \\
&\leq \mathrm{dist}_{\Delta_{G_k}}(x, y) + \text{const} < \infty.
\end{split}
\end{equation*}
This establishes the desired boundedness.
\end{proof}

\begin{lemma}
\label{lem::trans_Min_hyperbolic}
Let \( \{g_n\}_{n \geq 1} \subset G_k \) be a sequence of elements converging to an element \( g \in G_k \). Then, there exists \( N > 0 \) such that \( |g_n| = |g| \) for all \( n \geq N \).

Moreover, the sequence of closed, geodesically convex subsets \( \{\operatorname{Min}(g_n)\}_{n \geq 1} \subset \Delta_{G_k} \) converges to \( \operatorname{Min}(g) \) in the following sense:
\begin{enumerate}
    \item For every point \( x \in \operatorname{Min}(g) \), there exists a sequence \( \{x_{n} \in \operatorname{Min}(g_{n})\}_{n \geq 1} \) such that \( x_{n} \to x \) in the compact-open topology on \( \Delta_{G_k} \).
      \item If a sequence \( \{x_n \in \operatorname{Min}(g_n)\}_{n \geq 1} \) admits a convergent subsequence \( \{x_{n_j}\}_{j \geq 1} \) with strictly increasing indices, and \( x_{n_j} \to x \in \Delta_{G_k} \) in the compact-open topology, then \( x \in \operatorname{Min}(g) \).
\end{enumerate}
\end{lemma}

\begin{proof}
Fix a basepoint \( x_0 \in \operatorname{Min}(g) \subset \Delta_{G_k} \). By the definition of the compact-open topology on \( G_k \) and \( \Delta_{G_k} \), for every open ball \( B(x_0, r) \subset \Delta_{G_k} \) centered at \( x_0 \) with radius \( r > 0 \), there exists \( N > 0 \) such that
\[
g_n(B(x_0, r)) = g(B(x_0, r)) \quad \text{pointwise for all } n \geq N.
\]
From this, it follows that \( |g_n| = |g| \) for all sufficiently large \( n \), proving the first part of the lemma.  Indeed, if \( g \) is elliptic, then \( g(x_0) = x_0 \) since \( x_0 \in \operatorname{Min}(g) \), and hence \( g_n(x_0) = g(x_0) = x_0 \), implying that \( g_n \) is also elliptic.

If \( g \) is hyperbolic, we can choose \( r > 0 \) large enough so that \( g^2(x_0) \in B(x_0, r) \), i.e., the geodesic segment \( [x_0, g(x_0)] \) is translated by \( g \) within \( B(x_0, r) \), and the same holds for \( g_n \) for all sufficiently large \( n \). Since the orientation of the segment \( [x_0, g(x_0)] \) is preserved under translation by $g$ to \( [g(x_0), g^2(x_0)] \), it follows that \( g_n \) is hyperbolic with the same translation length \( |g_n| = |g| \).

\vspace{1em}
We now prove the second part of the lemma. Let \( x \in \operatorname{Min}(g) \). Choose \( r > 0 \) such that \( x \in B(x_0, r) \). Then, for all sufficiently large \( n \), we have \( g_n(x) = g(x) \). Therefore,
\[
|g| = \mathrm{dist}_{\Delta_{G_k}}(g(x), x) = \mathrm{dist}_{\Delta_{G_k}}(g_n(x), x) = |g_n|,
\]
which implies \( x \in \operatorname{Min}(g_n) \) for large \( n \). Thus, we may take \( x_n := x \) for all sufficiently large \( n \), yielding the desired sequence.

Now consider a subsequence \( \{x_{n_j} \in \operatorname{Min}(g_{n_j})\}_{j \geq 1} \subset \Delta_{G_k} \) converging to a point \( x \in \Delta_{G_k} \). Since the action of \( G_k \) on \( \Delta_{G_k} \) is continuous, the sequence \( \{g_{n_j}(x_{n_j})\}_{j \geq 1} \) converges to \( g(x) \in \Delta_{G_k} \). Moreover, because we are working in an affine building, the translation lengths \( \{|g_n|\}_{n \geq 1} \) stabilize to \( |g| \) for large \( n \). Hence,
\[
\mathrm{dist}_{\Delta_{G_k}}(g(x), x) = \lim_{j \to \infty} \mathrm{dist}_{\Delta_{G_k}}(g_{n_j}(x_{n_j}), x_{n_j}) = \lim_{j \to \infty} |g_{n_j}| = |g|,
\]
which implies \( x \in \operatorname{Min}(g) \) by the definition of the minimal set. This completes the proof.
\end{proof}

%\red{Add things about restricted simple roots from the other paper.} 
%\red{Show that $\sigma$ and $\theta(\sigma)$ are in $\partial_\infty \Delta_{G_k} \setminus \partial_\infty \Delta_{G_k}^{\theta}$.}
%\red{TO SHOW that $\xi_-,\xi_+$ are not ideal points of $\partial_\infty \Delta_{G_k}^{\theta}$. This is because $\theta$ is swapping $\xi_+$ and $\xi_-$, and $\theta$ fixes pointwise the geodesically convex and closed subset $\Delta_{G_k}^{\theta}$ as well as its visual boundary $\partial_\infty \Delta_{G_k}^{\theta}$.}
%\red{Show that $\Min(h)$ is $\theta$-stable if $h \in H_k$, but it does not mean that $\Min(h)$ is pointwise fixed by $\theta$.}

\section{Chabauty limits are subgroups of parabolics}

We prove that Chabauty limits of the fixed point group \( H_k \), associated with a \( k \)-involution \( \theta \) of \( G \), are contained within parabolic subgroups of \( G_k \). The analysis distinguishes between elliptic and hyperbolic cases and identifies the relevant parabolic subgroup through the limiting behavior of sequences in the Bruhat--Tits building of $G_k$.

In this section, we adopt assumptions that reflect the structural properties of maximal \( (\theta, k) \)-split tori in the group \( G \). Specifically, we work under the following hypotheses and notation:

\medskip
Let \( A \leq G_k \) be a subgroup satisfying the condition \textbf{(Hyp(\( \theta,k \))-split)}. Suppose that the apartment \( \mathcal{A} := \mathcal{A}_{\theta,k} \) is \( \theta \)-stable, i.e., \( \theta(\mathcal{A}) = \mathcal{A} \). Let \( \{a_\ell\}_{\ell \geq 1} \subset A \) be a sequence of hyperbolic elements satisfying the condition \textbf{(Same \( \sigma_\pm \))}.

\medskip
We now consider the sequence of closed subgroups \( \{ a_\ell H_k a_\ell^{-1} \}_{\ell \in \mathbb{N}} \). Since the space \( \mathcal{S}(G_k) \) of closed subgroups of \( G_k \) is compact with respect to the Chabauty topology, we may, after passing to a subsequence, assume that
\[
\{ a_\ell H_k a_\ell^{-1} \}_{\ell \in \mathbb{N}} \quad \text{converges to a subgroup} \quad L \in \mathcal{S}(G_k),
\]
which we denote by
\[
L := \limch_{\ell \to \infty} a_\ell H_k a_\ell^{-1} \leq G_k.
\]

The main result of this section is the following theorem.

\begin{theorem}
\label{thm::chabauty_in_parabolic}
We retain the notation from the conditions \textbf{(Hyp$(\theta,k)$-split)} and \textbf{(Same $\sigma_\pm$)}, and above hypotheses. Suppose that the apartment \( \mathcal{A} := \mathcal{A}_{\theta,k} \) is \( \theta \)-stable, i.e., \( \theta(\mathcal{A}) = \mathcal{A} \).  Then $L \leq P_{\sigma'_+}(k)$, where $P_{\sigma'_+}$ is the parabolic subgroup of $G$ associated with the unique ideal simplex $\sigma'_+$ containing $\xi_+$ in its interior and $P_{\xi_+}(k)$ are its $k$-points, i.e., $P_{\sigma'_+}(k) =\{ g \in G_k  \; \vert \; g(\sigma'_+)= \sigma'_+ \text{ setwise}\}$. %\red{Maybe $P_{\xi_+}$ instead of $P_{\sigma_+}$???}
\end{theorem}
\begin{proof}

To prove the theorem, it suffices to show that every element \( l \in L \) lies in the parabolic subgroup \( P_{\sigma'_+}(k) \).

By Proposition~\ref{prop::chabauty_conv}, for each \( l \in L \), there exists a sequence \( \{a_{\ell} h_{\ell} a_{\ell}^{-1}\}_{\ell \in \mathbb{N}} \), with \( h_{\ell} \in H_k \), such that $l = \lim\limits_{\ell \to \infty} a_{\ell} h_{\ell} a_{\ell}^{-1}$
with respect to the compact-open topology on \( G_k \). Applying the involution \( \theta \), we obtain $\theta(l) = \lim\limits_{\ell \to \infty} a_{\ell}^{-1} h_{\ell} a_{\ell}$.

Since \( G_k \), and hence \( H_k \), contains both two type of elements, elliptic and hyperbolic, we distinguish two cases.

\medskip
\noindent\textbf{Elliptic Case.} Suppose \( l \in L \leq G_k \) is an elliptic element. Then there exists a point \( x \in \Delta_{G_k} \) such that \( l(x) = x \), i.e., \( l \) fixes \( x \). By Lemma~\ref{lem::trans_Min_hyperbolic}.1 applied to \( l \), the sequence \( \{h_{\ell}\}_{\ell \geq 1} \subset H_k \) consists of elliptic elements of \( G_k \), and hence of \( H_k \) as well.  In particular, for sufficiently large \( \ell \), we have $a_{\ell} h_{\ell} a_{\ell}^{-1}(x) = x$, which implies $h_{\ell} a_{\ell}^{-1}(x) = a_{\ell}^{-1}(x)$ for all sufficiently large \( \ell \).

Since the point \( x \) is fixed, the condition \textbf{(Same \( \sigma_\pm \))} implies that the sequences \( \{a_{\ell_j}(x)\}_{j \geq 1} \) and \( \{a_{\ell_j}^{-1}(x)\}_{j \geq 1} \) converge respectively to the interior ideal point \( \xi_+ \) of the ideal simplex \( \sigma'_+ \), and to the interior ideal point \( \xi_- \) of the ideal simplex \( \theta(\sigma'_+) = \sigma'_- \), in the boundary at infinity \( \partial \Delta_{G_k} \).

Note that for sufficiently large \( \ell \geq 1 \), the element \( h_{\ell} \in H_k \) is elliptic. Consequently, it stabilizes the geodesically convex and closed subset \( \Delta_{G_k}^{\theta} \subset \Delta_{G_k} \), and fixes pointwise the element \( a_{\ell}^{-1}(x) \in \Delta_{G_k} \).

Consider the projection
\[
\proj_{\Delta_{G_k}^{\theta}}(a_{\ell}^{-1}(x)) \in \Delta_{G_k}^{\theta}
\]
of the point \( a_{\ell}^{-1}(x) \) onto \( \Delta_{G_k}^{\theta} \). Since both \( \Delta_{G_k} \) and \( \Delta_{G_k}^{\theta} \) are CAT(0) spaces, this projection is unique. Moreover, because \( h_{\ell} \) stabilizes \( \Delta_{G_k}^{\theta} \) and fixes \( a_{\ell}^{-1}(x) \), it also fixes the projection \( \proj_{\Delta_{G_k}^{\theta}}(a_{\ell}^{-1}(x)) \) pointwise.

Furthermore, by Lemma~\ref{lem::projection_sequence}, the sequence of projections
\[
\left\{ \proj_{\Delta_{G_k}^{\theta}}(a_{\ell}^{-1}(x)) \right\}_{\ell \geq 1} \subset \Delta_{G_k}^{\theta}
\]
forms a bounded subset of \( \Delta_{G_k}^{\theta} \). Therefore, the sequence $\left\{ a_{\ell} \left( \proj_{\Delta_{G_k}^{\theta}}(a_{\ell}^{-1}(x)) \right) \right\}_{\ell \geq 1}$ also converges to the ideal point \( \xi_+ \).

Now, since the action of \( G_k \) on \( \Delta_{G_k} \cup \partial \Delta_{G_k} \) is continuous with respect to the compact-open topology of \( G_k \) and the cone topology on \( \Delta_{G_k} \cup \partial \Delta_{G_k} \), we obtain:
\[
\lim_{\ell \to \infty} a_{\ell} h_{\ell} a_{\ell}^{-1} \left( a_{\ell} \left( \proj_{\Delta_{G_k}^{\theta}}(a_{\ell}^{-1}(x)) \right) \right)
= l(\xi_+) = \lim_{\ell \to \infty} a_{\ell} \left( \proj_{\Delta_{G_k}^{\theta}}(a_{\ell}^{-1}(x)) \right) = \xi_+.
\]
From this, we deduce that \( l \in P_{\xi_+} \), and since $\xi_+$ lies in the interior of $\sigma'_+$ we also obtain that $l \in  P_{\sigma'_+}$.

\medskip
\noindent\textbf{Hyperbolic Case.} Suppose \( l \in L \leq G_k \) is a hyperbolic element. Then \( l \) admits a translation axis and has positive translation length, i.e., \( |l| > 0 \).  By Lemma~\ref{lem::trans_Min_hyperbolic}.1 applied to \( l \), the sequence \( \{h_{\ell}\}_{\ell \geq 1} \subset H_k \) consists of hyperbolic elements of \( G_k \), and hence of \( H_k \), each with the same translation length as \( l \).  Since \( \theta \) is an automorphism of both \( G_k \) and \( \Delta_{G_k} \), we have $\theta(l) = \lim\limits_{\ell \to \infty} \theta(a_{\ell} h_{\ell} a_{\ell}^{-1}) = \lim\limits_{\ell \to \infty} a_{\ell}^{-1} h_{\ell} a_{\ell}$.  It follows that \( \theta(l) \in G_k \) is also a hyperbolic element with the same translation length as \( l \), i.e., \( |l| = |\theta(l)| \).  Therefore, all the properties established for \( l \) also hold for \( \theta(l) = \lim\limits_{\ell \to \infty} a_{\ell}^{-1} h_{\ell} a_{\ell} \).

Now consider a point \( y \in \Min(l) \subset \Delta_{G_k} \), meaning that \( y \) lies on a translation axis of the hyperbolic element \( l \). By Lemma~\ref{lem::trans_Min_hyperbolic}, there exists \( N > 0 \) such that for every \( \ell \geq N \), we have:
\[
a_{\ell} h_{\ell} a_{\ell}^{-1}(y) = l(y), \quad \text{and} \quad y \in \Min(a_{\ell} h_{\ell} a_{\ell}^{-1}).
\]
Moreover, it follows that \( a_{\ell}^{-1}(y) \in \Min(h_{\ell}) \) for all \( \ell \geq N \).  Since \( \theta \) is an automorphism of \( \Delta_{G_k} \), we also have \( \theta(y) \in \Min(\theta(l)) \). Consequently, \( \theta(y) \in \Min(a_{\ell}^{-1} h_{\ell} a_{\ell}) \), and thus $a_{\ell}(\theta(y)) \in \Min(h_{\ell})$
for every \( \ell \geq N \). In particular, since \( \Min(h_{\ell}) \) is a geodesically convex subset of \( \Delta_{G_k} \), it follows that the geodesic segment  $[a_{\ell}^{-1}(y), \; a_{\ell}(\theta(y))]$ is entirely contained in \( \Min(h_{\ell}) \) for all \( \ell \geq N \).

Since the chosen point \( y \), and consequently \( \theta(y) \), are fixed points in \( \Delta_{G_k} \),  the condition \textbf{(Same \( \sigma_\pm \))} implies that 
\[
\lim_{\ell \to \infty} a_{\ell}^{-1}(y) = \xi_-, \quad \text{and} \quad \lim_{\ell \to \infty} a_{\ell}(\theta(y)) = \xi_+.
\]
%Because \( \xi_- \) and \( \xi_+ \) are not ideal points of the boundary \( \partial \Delta_{G_k}^{\theta} \), we may assume that $a_{\ell}(\theta(y)), \; a_{\ell}^{-1}(y) \in \Min(h_{\ell}) \setminus \Delta_{G_k}^{\theta}$ for every \( \ell \geq N \).

%Consider the geodesic ray \( [y, \xi_+) \). Then there exists a flat \( F \subset \Delta_{G_k} \), not necessarily of maximal dimension, contained in some apartment \( \mathcal{A}' \subset \Delta_{G_k} \), such that the entire geodesic ray \( [y, \xi_+) \) is contained in \( F \subset \mathcal{A}' \). Indeed, by the properties of affine buildings and their visual boundaries, there exists an apartment \( \mathcal{A}' \subset \Delta_{G_k} \) such that \( y \in \mathcal{A}' \) and \( \sigma'_+ \subset \partial_\infty \mathcal{A}' \). From this, one can select a flat \( F \subset \mathcal{A}' \) that contains the geodesic ray \( [y, \xi_+) \) and satisfies the desired properties.

We claim that the geodesic ray \( [y, \xi_+) \) is entirely contained in \( \Min(l) \). Indeed, since the geodesic segment \( [a_{\ell}^{-1}(y), a_{\ell}(\theta(y))] \) is entirely contained in \( \Min(h_{\ell}) \) for every \( \ell \geq N \), it follows that the segment
\[
[y, a_{\ell}^{2}(\theta(y))]
\]
is entirely contained in \( \Min(a_{\ell} h_{\ell} a_{\ell}^{-1}) = a_{\ell} \Min(h_{\ell}) \) for all \( \ell \geq N \). Now, by Lemma~\ref{lem::trans_Min_hyperbolic}, we have $\lim_{\ell \to \infty} \Min(a_{\ell} h_{\ell} a_{\ell}^{-1}) = \Min(l)$.
Since \( y \in \Min(l) \), and \( \lim_{\ell \to \infty} a_{\ell}^{2}(\theta(y)) = \xi_+ \), and because \( \Delta_{G_k} \) is locally finite, there must exist a point \( w \in (y, \xi_+) \) such that \( w \in \Min(l) \).  By induction and repeating the same argument, we can find points arbitrarily far along the geodesic ray \( [y, \xi_+) \) that lie in \( \Min(l) \). Hence, every point on the ray belongs to \( \Min(l) \), and the claim follows.

From the claim above, we know that the geodesic rays \( [y, \xi_+) \) and \( [l(y), l(\xi_+)) \) are both contained in \( \Min(l) \). This implies that they are asymptotic, and therefore determine the same endpoint in the visual boundary \( \partial \Delta_{G_k} \). In particular, we have $
\xi_+ = l(\xi_+)$.  Since \( \xi_+ \) lies in the interior of the ideal simplex \( \sigma'_+ \), it follows that, $l(\sigma'_+) = \sigma'_+$, and hence \( l \in P_{\sigma'_+}(k) \), as desired.
\end{proof}

\section{Elements of $L$ in the Levi factor}

In this section, we work under general assumptions where \( A \leq G_k \) is a subgroup satisfying the hypothesis \textbf{(Hyp(\( \theta, k \))-split)}. We assume that the apartment \( \mathcal{A} := \mathcal{A}_{\theta,k} \) is \( \theta \)-stable, i.e., \( \theta(\mathcal{A}) = \mathcal{A} \). Let \( \{a_\ell\}_{\ell \geq 1} \subset A \) be a sequence of hyperbolic elements satisfying the condition \textbf{(Same \( \sigma_\pm \))}, and define the Chabauty limit group
\[
L := \limch_{\ell \to \infty} a_\ell H_k a_\ell^{-1} \leq G_k.
\]

Retaining the notation from \textbf{(Hyp(\( \theta, k \))-split)} and \textbf{(Same \( \sigma_\pm \))}, let \( \sigma'_\pm \) denote the unique ideal simplices in the visual boundary \( \partial \Delta_{G_k} \) that contain the endpoints \( \xi_\pm \) in their interiors. Note that \( \sigma'_+ \) is opposite to \( \sigma'_- \).

The goal of this section is to analyze the structure of the Chabauty limit group \( L \), particularly its intersection with the Levi factor \( M_{\sigma'_-, \sigma'_+} \) of the parabolic subgroup \( P_{\sigma'_+} \). We show how elements of \( L \) arise from sequences in \( H_k \), and we characterize their elliptic and hyperbolic components within the Levi decomposition of $P_{\sigma'_+}$.

\subsection{On elements of $H_k \cap M_{\sigma'_-,\sigma'_+}$}

\begin{proposition}
\label{prop::conv_geod_rays}
Let \( l \in L \) be such that $l = \lim\limits_{\ell \to \infty} a_{\ell} h_{\ell} a_{\ell}^{-1}$,
with \( h_{\ell} \in H_k \). Fix a point \( y \in \Min(l) \). Then the sequence of geodesic segments
\[
\left\{ [a_{\ell}^{-1}(y), a_{\ell}(\theta(y))] \right\}_{\ell \geq 0} \subset \Delta_{G_k}
\]
admits a convergent subsequence, reindexed by \( \ell_j \), which converges to the bi-infinite geodesic line \( (\xi_-, \xi_+) \subset \Delta_{G_k} \).  In particular, there exists a point \( z \in (\xi_-, \xi_+) \) and a sequence $\left\{ z_{\ell_j} \in [a_{\ell_j}^{-1}(y), a_{\ell_j}(\theta(y))] \right\}_{j \geq 0}$
such that $z = \lim\limits_{j \to \infty} z_{\ell_j}$.
Moreover, one can extract a subsequence \( \{h_{\ell_j}\}_{j \geq 1} \subset H_k \) such that \( z \in \Min(h_{\ell_j}) \) and $h_{\ell_j}(z) = h_{\ell_i}(z)$
for all sufficiently large \( j, i \).
\end{proposition}

\begin{proof}
Take the fixed point \( y \in \Min(l) \) as the basepoint for the compact-open topology on \( \Delta_{G_k} \). Then, for every open ball \( B(y, r) \subset \Delta_{G_k} \) centered at \( y \) with radius \( r > 0 \), there exists \( N > 0 \) such that
\[
a_{\ell} h_{\ell} a_{\ell}^{-1}(B(y, r)) = l(B(y, r))
\]
pointwise for all \( \ell \geq N \). In particular, we can choose \( r \) such that \( l(y) \in B(y, r) \).  Since \( |a_{\ell} h_{\ell} a_{\ell}^{-1}| = |h_{\ell}| = |l| \) for all sufficiently large \( \ell \) (see Lemma~\ref{lem::trans_Min_hyperbolic}), it follows that
\[
|l| = \dist_{\Delta_{G_k}}(l(y), y) = \dist_{\Delta_{G_k}}(a_{\ell} h_{\ell} a_{\ell}^{-1}(y), y) = \dist_{\Delta_{G_k}}(h_{\ell} a_{\ell}^{-1}(y), a_{\ell}^{-1}(y)) = |h_{\ell}|.
\]
This implies that \( a_{\ell}^{-1}(y) \in \Min(h_{\ell}) \) for all sufficiently large \( \ell \).

Now apply the \( k \)-involution \( \theta \), which acts as an automorphism of \( \Delta_{G_k} \). Under the hypothesis \textbf{(Hyp(\( \theta, k \))-split)}, we obtain: $\theta(l) = \lim\limits_{\ell \to \infty} \theta(a_{\ell} h_{\ell} a_{\ell}^{-1}) = \lim\limits_{\ell \to \infty} a_{\ell}^{-1} h_{\ell} a_{\ell}$.  Moreover, we have \( \Min(\theta(l)) = \theta(\Min(l)) \), since for every \( x \in \Delta_{G_k} \),
\[
\dist_{\Delta_{G_k}}(l(x), x) = \dist_{\Delta_{G_k}}(\theta(l(x)), \theta(x)) = \dist_{\Delta_{G_k}}(\theta(l)(\theta(x)), \theta(x)),
\]
which also implies that \( |l| = |\theta(l)| \).  Then, as \( \theta(y) \in \Min(\theta(l)) \), it follows that \( a_{\ell}(\theta(y)) \in \Min(h_{\ell}) \) for all sufficiently large \( \ell \).

Combining the previous two facts, there exists \( N > 0 \) such that the geodesic segment
\[
[a_{\ell}^{-1}(y), a_{\ell}(\theta(y))] \subset \Min(h_{\ell})
\]
for every \( \ell > N \).  In addition, by Lemma~\ref{lem::projection_point_building}, we have
\[
\proj_{\Delta_{G_k}^{\theta}}(a_{\ell}^{-1}(y)) = \operatorname{mid}(a_{\ell}^{-1}(y), a_{\ell}(\theta(y))) = \proj_{\Delta_{G_k}^{\theta}}(a_{\ell}(\theta(y))) \in \Min(h_{\ell}),
\]
for every \( \ell > N \).   Moreover, by Lemma~\ref{lem::projection_sequence}, the sequence of projections
\[
\{z_{\ell} := \proj_{\Delta_{G_k}^{\theta}}(a_{\ell}^{-1}(y))\}_{\ell \geq 1} \subset \Delta_{G_k}^{\theta}
\]
is bounded in \( \Delta_{G_k}^{\theta} \), and therefore admits a convergent subsequence to some point \( z \in \Delta_{G_k}^{\theta} \). We denote this convergent subsequence by the index \( \ell_j \).  Observe that since \( y \) is fixed and $\{a_{\ell_j}\}_{j\in \NN}$ satisfies hypothesis \textbf{(Same \( \sigma_\pm \))}, we have \( a_{\ell_j}^{-1}(y) \to \xi_- \) and \( a_{\ell_j}(\theta(y)) \to \xi_+ \) as \( \ell_j \to \infty \). Up to re-extracting a subsequence, the sequence of geodesic segments
\[
\{[a_{\ell_j}^{-1}(y), a_{\ell_j}(\theta(y))]\}_{j \geq 0}
\]
converges to the bi-infinite geodesic line \( (\xi_-, \xi_+) \subset \Delta_{G_k} \) that passes through the point \( z \).

It remains to show that, up to re-extracting a subsequence, we have \( z \in \Min(h_{\ell_j}) \) for all sufficiently large \( \ell_j \). Since \( \Delta_{G_k} \) is a simplicial complex, there exists a unique simplex \( \sigma \subset \Delta_{G_k} \) such that \( z \) lies in its interior (if \( z \) is a vertex, then \( \sigma \) is that vertex).  For the sequence \( z_{\ell_j} \) to converge to \( z \), each \( z_{\ell_j} \) must lie either in the interior of a larger simplex that has \( \sigma \) as a face, or in the interior of \( \sigma \) itself.  Now, although the action of \( G_k \) preserves the simplicial structure of \( \Delta_{G_k} \), it is not necessarily type-preserving. Thus, an elliptic element \( g \in G_k \) does not necessarily fix pointwise a simplex containing a point from \( \Min(g) \) in its interior, but it should fix pointwise its barycenter. Similarly, if \( g \in G_k \) is hyperbolic, a simplex containing a point from \( \Min(g) \) in its interior may not be entirely contained in \( \Min(g) \), since \( g \) may not preserve vertex types.  However, by considering the barycentric subdivision of the simplicial complex \( \Delta_{G_k} \), we obtain a refinement where, for both elliptic and hyperbolic elements \( g \in G_k \), any barycentric subdivision simplex that contains a point from \( \Min(g) \) in its interior is entirely contained in \( \Min(g) \), including its closure.

With the above facts in mind, and since the limit point \( z \) can lie only in the closure of finitely many barycentric subdivision simplices, we conclude -- after possibly extracting a further subsequence from \( \{\ell_j\}_{j \geq 0} \) -- that \( z \in \Min(h_{\ell_j}) \) for all sufficiently large \( \ell_j \).

To prove the final part of the proposition, recall that the first part of the proof established \( |h_{\ell}| = |l| \) for all sufficiently large \( \ell \). If \( |l| = 0 \), i.e., \( l \) is elliptic, then we immediately obtain \( z = h_{\ell_j}(z) \) for all large \( \ell_j \).

Now consider the case \( |l| > 0 \), i.e., \( l \) is hyperbolic. Since \( z \in \Min(h_{\ell_j}) \) for all large \( \ell_j \), the set
\[
\{h_{\ell_j}(z) \in \Min(h_{\ell_j})\}_{j \text{ large enough}}
\]
is contained in the closed ball \( \overline{B(z, 2|l|)} \) centered at \( z \) with radius \( |l| \). Because we are working in a locally finite simplicial complex, we can extract a further subsequence to obtain the desired conclusion.
\end{proof}

Retaining the notation from the conditions \textbf{(Hyp(\( \theta, k \))-split)} and \textbf{(Same \( \sigma_\pm \))}, recall that by Theorem~\ref{thm::chabauty_in_parabolic}, we have \( L \leq P_{\sigma'_+}(k) \). Note that \( \sigma'_+ \) is a subsimplex of \( \sigma_+ \), and let \( \sigma'_- \) be the subsimplex of \( \sigma_- \) that is opposite to \( \sigma'_+ \). Then \( \xi_- \) lies in the interior of \( \sigma'_- \).

Define $M_{\sigma'_-, \sigma'_+} := P_{\sigma'_+} \cap P_{\sigma'_-}$, which is called a Levi factor of \( P_{\sigma'_+} \).

\begin{corollary}
\label{cor::find_elements_H}
Let \( l \in L \) be such that $l = \lim\limits_{\ell \to \infty} a_{\ell} h_{\ell} a_{\ell}^{-1}$,
with \( h_{\ell} \in H_k \).  Then there is an element $h \in H_k\cap M_{\sigma'_-,\sigma'_+}$ such that (up to extracting a subsequence $\{\ell_j\}_{j \in\NN}$) $\{h_{\ell_j}\}_{j\geq 1}$ converges to $h \in H_k$. If moreover $M_{\sigma'_-,\sigma'_+}(k) \leq Z_{G_k}(A)$, then $h=\lim\limits_{j \to \infty} a_{\ell_j} h a_{\ell_j}^{-1} $ is an element of $L$, and $l h^{-1}\in L$ is an elliptic element in $G_k$. 
\end{corollary}

\begin{proof}
From Proposition~\ref{prop::conv_geod_rays}, applied to $l = \lim\limits_{\ell \to \infty} a_{\ell} h_{\ell} a_{\ell}^{-1}$,
we can select a point \( z \in \Delta_{G_k} \) and extract a subsequence \( \{h_{\ell_j}\}_{j \geq 1} \subset H_k \)  satisfying
\[
z \in \Min(h_{\ell_j}) \quad \text{and} \quad h_{\ell_j}(z) = h_{\ell_i}(z),
\]
for all sufficiently large \( j, i \).

Fixing an index \( \ell_i \), we obtain \( h_{\ell_j}(z) = h_{\ell_i}(z) \) for all large \( j \), which implies $h_{\ell_i}^{-1} h_{\ell_j}(z) = z$.
This shows that $h_{\ell_j} = h_{\ell_i} h'_{\ell_j}$,
for all large \( j \), where \( \{h'_{\ell_j}\}_{j \geq 1} \subset H_k \) is a sequence of elliptic elements satisfying \( h'_{\ell_j}(z) = z \) for every \( j \).

Therefore, \( z \in \Min(h'_{\ell_j}) \) for all \( j \), and moreover, we can extract a convergent subsequence (denoted with the same indexing) such that
\[
\lim_{j \to \infty} h'_{\ell_j} = h' \in H_k,
\]
where \( h' \) is elliptic and \( z \in \Min(h') \). This implies that the sequence \( \{h_{\ell_j}\}_{j \geq 1} \subset H_k \) also converges in \( H_k \) to the element \( h := h_{\ell_i} h' \in H_k \) with $z \in \Min(h)$, i.e.,
\[
\lim_{j \to \infty} h_{\ell_j} = h \in H_k.
\]
 
By Lemma~\ref{lem::trans_Min_hyperbolic}, the sequence \( \{\Min(h_{\ell_j})\}_{j \geq 1} \subset \Delta_{G_k} \) also converges to \( \Min(h) \subset \Delta_{G_k} \), and we have \( |h_{\ell_j}| = |h| \) for all sufficiently large \( j \).  Moreover, as shown in the proof of Lemma~\ref{lem::trans_Min_hyperbolic}, the construction of the point \( z \) arises from a fixed point \( y \in \Min(l) \), with the property that the sequence of geodesic segments
\[
\left\{ [a_{\ell_j}^{-1}(y), a_{\ell_j}(\theta(y))] \right\}_{j \geq 1} \subset \Min(h_{\ell_j})
\]
converges to a specific bi-infinite geodesic line \( (\xi_-, \xi_+) \subset \Delta_{G_k} \), with
\[
z = \lim_{j \to \infty} \proj_{\Delta_{G_k}^{\theta}}(a_{\ell_j}^{-1}(y)).
\]
From these facts, it follows that \( h \) must fix both \( \xi_- \) and \( \xi_+ \), and therefore \( h \in H_k \cap M_{\xi_-, \xi_+} \), where $M_{\xi_-, \xi_+} := P_{\xi_+} \cap P_{\xi_-}$.  Since $\xi_\pm$ in the interior of $\sigma'_\pm$, we have $P_{\xi_\pm}= P_{\sigma'_\pm}$ and so $M_{\xi_-,\xi_+}=M_{\sigma'_-,\sigma'_+}$.

If moreover $M_{\sigma'_-,\sigma'_+}(k) \leq Z_{G_k}(A)$, then $h \in H_k \cap M_{\xi_-, \xi_+}$ commutes with all the elements of $A$, in particular with $\{a_{\ell_j}\}_{j \in \NN}$ as well. Then $\lim\limits_{j\to \infty} h_{\ell_j} h^{-1}=\id \in H_k$ and 
$$\lim\limits_{j\to \infty} a_{\ell_j} h_{\ell_j} a_{\ell_j}^{-1} h^{-1}=\lim\limits_{k\to \infty} a_{\ell_j} h_{\ell_j} h^{-1} a_{\ell_j}^{-1}=l h^{-1} \in G_k. $$
Since $\{h_{\ell_j}h^{-1}\}_{j \geq 1} \subset H_k$, by Proposition \ref{prop::chabauty_conv}, $lh^{-1}$ is an elliptic element of  $L$ with  $h\in L$. 
 \end{proof}

\begin{remark}
If \( M_{\sigma'_-, \sigma'_+}(k) \leq Z_{G_k}(A) \), then for every nontrivial element \( l \in L \), we can find an element \( h \in H_k \cap M_{\sigma'_-, \sigma'_+} \) such that \( h \in L \) and \( lh^{-1} \in L \) is elliptic in \( G_k \).  It is straightforward to verify that every element \( g \in H_k \cap M_{\sigma'_-, \sigma'_+} \) lies in \( L \). Moreover, by Corollary~\ref{cor::find_elements_H}, any hyperbolic element of \( L \leq P_{\sigma'_+} \) arises from a hyperbolic element in \( H_k \cap M_{\sigma'_-, \sigma'_+} \). Regarding hyperbolic elements of \( L \), we do not obtain any beyond those.
\end{remark}

\subsection{On elliptic elements of $L \cap M_{\sigma'_-,\sigma'_+}$}

Before discussing the elliptic elements of \( L \cap M_{\sigma'_-,\sigma'_+} \), we first recall from \cite[Section 3, Further decomposition of \( M_I \)]{CiLe_p} several key facts about the sub-buildings within \( \Delta_{G_k} \) that correspond to Levi factors of parabolic subgroups of \( G_k \).

Let \( \Sigma \) be an apartment of \( \Delta_{G_k} \), and let \( \eta_{+} \) and \( \eta_{-} \) be two opposite simplices in \( \partial \Sigma \). Let \( (W_I, I) \) denote the subgroup of the finite Weyl group \( (W_{\mathrm{fin}}, S) \) associated with \( G_k \) that stabilizes the simplex \( \eta_{+} \). The \textbf{residue} \( \operatorname{res}(\eta_{+}) \) is defined as the set of all ideal chambers in \( \Ch(\partial \Delta_{G_k}) \) that contain \( \eta_{+} \). Let \( \Delta(\eta_{+}, \eta_{-}) \) denote the union of all apartments in \( \Delta_{G_k} \) whose ideal boundaries contain both \( \eta_{+} \) and \( \eta_{-} \). Let \( n \) denote the dimension of \( \Delta_{G_k} \).

\begin{proposition}[Proposition 3.15 in \cite{CiLe_p}]
\label{prop::res_building}
We retain the notation introduced above. Then the union of apartments \( \Delta(\eta_{+}, \eta_{-}) \) forms a closed convex subset of \( \Delta_{G_k} \). It is an extended, locally finite, thick affine building, and satisfies
\[
\Delta(\eta_{+}, \eta_{-}) \cong \mathbb{R}^{n - |I|} \times \Delta_I,
\]
where \( \Delta_I \) is a locally finite affine building of dimension \( |I| \). Moreover, the residue \( \operatorname{res}(\eta_{+}) \) is a compact subset of \( \Ch(\partial \Delta_{G_k}) \), it is a spherical building, and
\[
\operatorname{res}(\sigma_{+}) \cong \Ch(\partial \Delta_I).
\]
\end{proposition}

When \( \eta_{+} \) and \( \eta_{-} \) are two opposite chambers in \( \partial \Sigma \), we have \( \Delta(\eta_{+}, \eta_{-}) \cong \mathbb{R}^{n} \cong \Sigma \), and \( \Delta_I \) reduces to a point. Note that the Euclidean factor \( \mathbb{R}^{n - |I|} \) in the splitting of \( \Delta(\eta_{+}, \eta_{-}) \) contains the simplices \( \eta_{+} \) and \( \eta_{-} \) in its boundary at infinity.

Recall the notation: \( M_I := M_{\eta_-, \eta_+} := P_{\eta_+} \cap P_{\eta_-} \), which is referred to as a Levi factor of \( P_{\eta_+} \), with \( M_I(k) \) denoting its \( k \)-points.

\begin{proposition}[Proposition 3.17 in \cite{CiLe_p}]
\label{prop::res_building_str_tran}
We retain the notation introduced above. Then \( M_I(k) \) acts on \( \Delta(\eta_{+}, \eta_{-}) \) and preserves the splitting $\Delta(\eta_{+}, \eta_{-}) \cong \mathbb{R}^{n - |I|} \times \Delta_I$.  Moreover, the induced map
\[
\alpha : M_I(k) \to \Isom(\Delta(\eta_{+}, \eta_{-})) = \Isom(\mathbb{R}^{n - |I|}) \times \Isom(\Delta_I)
\]
is a group homomorphism, and the normal subgroups
\[
H_I := \alpha^{-1}\left( \alpha|_{\Isom(\Delta_I)}(M_I(k)) \right), \quad \text{and} \quad H^I := \alpha^{-1}\left( \alpha|_{\Isom(\mathbb{R}^{n - |I|})}(M_I(k)) \right)
\]
of \( M_I(k) \) act by automorphisms on \( \Delta_I \), and by translations on \( \mathbb{R}^{n - |I|} \), respectively.
\end{proposition}

\begin{remark}
\label{rem::alpha_map}
In general the map $\alpha : M_{I}(k) \to \Isom(\RR^{n-\vert I \vert}) \times \Isom(\Delta_{I})$ from Proposition \ref{prop::res_building_str_tran} is not injective since the subgroup $\Ker(\alpha) \leq M_{I}(k)$ that acts trivially on $\Delta(\eta_{+}, \eta_{-}) \cong \RR^{n-\vert I \vert} \times \Delta_{I}$ might be non-trivial.
In addition, the image subgroup $\alpha(M_I(k)$ might not split as a direct product in  $\Isom(\RR^{n-\vert I \vert}) \times \Isom(\Delta_{I})$. 
Both issues occur in the special case $ M_{I}(\QQ_p) \leq G_{\QQ_p}=\SL(n, \QQ_p)$.
\end{remark}

We are now ready to prove the following result, when $\sigma_\pm= \sigma'_\pm$, i.e., $\xi_\pm$ are in the interior of $\sigma_\pm$.
 
\begin{proposition}
\label{prop::levi_factors}
Assume that $\sigma'_\pm = \sigma_\pm$ and $M_{\sigma_-,\sigma_+}(k) \leq Z_{G_k}(A)$, and let \( l \in L \) be such that $l = \lim\limits_{\ell \to \infty} a_{\ell} h_{\ell} a_{\ell}^{-1}$, with \( h_{\ell} \in H_k \). If \( l \in P_{\sigma_-} \) (and hence \( l \in M_{\sigma_-,\sigma_+} \)), then
\[
l \in \Ker(\alpha)^0\cdot \left( H_k \cap M_{\sigma_-,\sigma_+} \right),
\]
where \( \alpha \) is the map from Proposition~\ref{prop::res_building_str_tran}, and $\Ker(\alpha)^0 := \left\{ g \in \Ker(\alpha) \;\middle|\; g \text{ is elliptic} \right\}$.
\end{proposition}

\begin{proof}
Since \( l \in L \leq P_{\sigma_+} \) and, the by hypothesis \( l \in P_{\sigma_-} \), it follows that $l \in P_{\sigma_-} \cap P_{\sigma_+} =: M_{\sigma_-,\sigma_+}$. Now, by Corollary~\ref{cor::find_elements_H}, and possibly after extracting a subsequence of indices \( \{ \ell_j \}_{j \geq 1} \), we know there exists an element \( h \in H_k \cap M_{\sigma_-,\sigma_+} \leq L \) such that
\[
h = \lim_{j \to \infty} h_{\ell_j}, \quad lh^{-1} \in L, \quad \text{and} \quad lh^{-1} \text{ is elliptic}.
\]
In particular, by hypothesis, \( lh^{-1} \in M_{\sigma_-,\sigma_+} \), and the sequence \( \{ h_{\ell_j} h^{-1} \}_{j \geq 1} \) converges to the identity in \( H_k \). Therefore, we are reduced to the following case: \( l \in L \cap M_{\sigma_-,\sigma_+} \) is elliptic, with
\[
l = \lim_{j \to \infty} a_{\ell_j} h_{\ell_j} a_{\ell_j}^{-1}, \quad \text{and} \quad \{ h_{\ell_j} \}_{j \geq 1} \to e \text{ in } H_k.
\]
This implies that, for sufficiently large \( \ell_j \), the elements \( h_{\ell_j} \) act by type-preserving automorphisms on \( \Delta_{G_k} \cup \partial \Delta_{G_k} \). In particular, this means that \( l \) also acts type-preservingly on \( \Delta_{G_k} \cup \partial \Delta_{G_k} \).

Next, we aim to compare the actions of \( l \) and the identity element \( e \), both of which lie in \( M_{\sigma_-,\sigma_+} \), on the locally finite affine building \( \Delta_I \), associated with \( \eta_{\pm} = \sigma_{\pm} \), as described in Propositions~\ref{prop::res_building} and~\ref{prop::res_building_str_tran}.

Since \( l \in M_{\sigma_-,\sigma_+}(k) \) is elliptic and type-preserving on \( \Delta_{G_k} \cup \partial \Delta_{G_k} \), there exists a flat \( F \subset \Delta_{G_k} \) that is fixed pointwise by \( l \) and has the simplices \( \sigma_+ \) and \( \sigma_- \) in its boundary at infinity \( \partial F \). Note that this flat \( F \) is contained in \( \Delta(\sigma_+, \sigma_-) \). Moreover, \( F \) is parallel to the apartment \( \mathcal{A} = \mathcal{A}_{\theta,k} \) from the conditions \textbf{(Hyp(\( \theta, k \))-split)} and \textbf{(Same \( \sigma_\pm \))}.

In particular, since we assume \( \sigma'_\pm = \sigma_\pm \), and the attracting and repelling endpoints of \( \{a_{\ell}\}_{\ell \in \mathbb{N}} \) lie in the interior of \( \sigma_\pm \), the flats \( \{a_\ell^{\pm1}(F)\}_{\ell \geq 1} \) are all parallel to \( \mathcal{A} \) (because \( a_{\ell}^{\pm1}(\sigma_\pm) = \sigma_\pm \) setwise), and each is at the same bounded distance from \( \mathcal{A} \) as \( F \) is.  Now, because \( \Delta_{G_k} \) is a locally finite simplicial complex, we may extract a subsequence if necessary and assume that the sequences \( \{a_{\ell_j}(F)\}_{j \geq 1} \) and \( \{a_{\ell_j}^{-1}(F)\}_{j \geq 1} \) are constant and equal to flats \( F' \) and \( F'' \), respectively, in \( \Delta(\sigma_+, \sigma_-) \subset \Delta_{G_k} \), each having \( \sigma_\pm \) in its boundary at infinity.

Let \( x \in F \), and let \( B(x, r) \) denote the open ball of radius \( r > 0 \) in \( \Delta_{G_k} \) centered at \( x \). We may choose \( r \) to be arbitrarily large. Since $l = \lim\limits_{j \to \infty} a_{\ell_j} h_{\ell_j} a_{\ell_j}^{-1}$, there exists \( N > 0 \) such that for every \( \ell_j \geq N \), we have
\[
l(B(x, r)) = a_{\ell_j} h_{\ell_j} a_{\ell_j}^{-1}(B(x, r)) \text{ pointwise,}
\] with \( l(B(x, r)) = B(x, r) \) setwise and \( l(B(x, r) \cap F) = B(x, r) \cap F \) pointwise.  In particular, this implies
\[
a_{\ell_j} h_{\ell_j} a_{\ell_j}^{-1}(B(x, r)) = B(x, r) \text{ setwise,} 
\]
and
\[
a_{\ell_j} h_{\ell_j} a_{\ell_j}^{-1}(B(x, r) \cap F) = B(x, r) \cap F \text{ pointwise,}
\]
for every \( \ell_j \geq N \). Consequently, we obtain
\[
h_{\ell_j} a_{\ell_j}^{-1}(B(x, r)) = a_{\ell_j}^{-1}(B(x, r)) \text{ setwise,}
\]
and
\[
h_{\ell_j} a_{\ell_j}^{-1}(B(x, r) \cap F) = a_{\ell_j}^{-1}(B(x, r) \cap F) = a_{\ell_j}^{-1}(B(x, r)) \cap a_{\ell_j}^{-1}(F) = a_{\ell_j}^{-1}(B(x, r)) \cap F''
\]
pointwise.

We now recall that the sequence \( \{h_{\ell_j}\}_{j \geq 1} \) converges to the identity element \( e \) in \( H_k \), and hence in \( G_k \). Let \( y \in \mathcal{A} \), and consider the open ball \( B(y, R) \) of radius \( R > 0 \) in \( \Delta_{G_k} \), centered at \( y \). We choose \( R \) large enough so that the flats \( F' \) and \( F'' \), which lie at bounded distance from \( \mathcal{A} \), intersect \( B(y, R) \). This is possible because through any point \( y \in \mathcal{A} \), one can draw a flat in \( \mathcal{A} \) that is parallel to \( F \), and hence to \( F' \) and \( F'' \), and which contains \( \sigma_\pm \) in its boundary at infinity.

Then, there exists \( M > 0 \) such that for every \( \ell_j \geq M \), we have
\[
h_{\ell_j}(B(y, R)) = B(y, R) \text{ pointwise.}
\] In particular, restricting the trivial action of $h_{\ell_j}$ on $ B(y,R) \cap \Delta(\sigma_+, \sigma_-)$, we can abuse notation and apply Proposition~\ref{prop::res_building_str_tran} to obtain
\[
\alpha(h_{\ell_j})|_{\Delta_I} \left( \proj_{\Delta_I}(B(y, R)) \right) = \proj_{\Delta_I}(B(y, R)) \text{ pointwise,}
\]
for every \( \ell_j \geq M \). Here, \( \proj_{\Delta_I}(B(y, R)) \) denotes the projection of \( B(y, R) \cap \Delta(\sigma_+, \sigma_-) \) onto the factor \( \Delta_I \) in the decomposition of \( \Delta(\sigma_+, \sigma_-) \). Thus, as \( \ell_j \) increases, the elements \( \alpha(h_{\ell_j}) \) converge to the identity automorphism of the affine building \( \Delta_I \).

We now recall that the sequence \( \{h_{\ell_j}\}_{j \geq 1} \) converges to the identity element \( e \) in \( H_k \), and hence in \( G_k \). Consider the action of \( h_{\ell_j} \) on \( \Delta_I \), viewed from the ball \( a_{\ell_j}^{-1}(B(x, r)) \), which intersects \( F'' \) non-trivially since \( a_{\ell_j}^{-1}(x) \in F'' = a_{\ell_j}^{-1}(F) \). One can choose \( r \) sufficiently large so that the intersection \( a_{\ell_j}^{-1}(B(x, r)) \cap F'' \) has the same dimension as the flat \( F'' \).

Because \( h_{\ell_j} \) fixes pointwise both \( F'' \cap a_{\ell_j}^{-1}(B(x, r)) \) and \( F'' \cap B(y, R) \), as well as the entire ball \( B(y, R) \), we can consider transversal sections of \( \Delta_I \) through points in these intersections. Specifically, taking the transversal section of \( \Delta_I \) in \( a_{\ell_j}^{-1}(B(x, r)) \) through a point in \( F'' \cap a_{\ell_j}^{-1}(B(x, r)) \), and similarly in \( B(y, R) \) through a point in \( F'' \cap B(y, R) \), we observe that the corresponding actions of \( \alpha(h_{\ell_j}) \) must coincide. Since the action on \( \proj_{\Delta_I}(B(y, R)) \) is trivial, the action on the transversal section in \( a_{\ell_j}^{-1}(B(x, r)) \) must also be trivial.  Indeed, as \( h_{\ell_j} \) is an automorphism of \( \Delta_{G_k} \), it preserves parallelism of flats.

Therefore, by choosing \( R \) sufficiently large, we conclude that the action of \( \alpha(l) \) on the transversal section of \( \Delta_I \) in \( B(x, r) \), through a point in \( F \cap B(x, r) \), is also trivial. Since this holds for every radius \( r > 0 \), it follows that \( \alpha(l)|_{\Delta_I} \) is trivial.

Moreover, because \( l \) is elliptic, the action \( \alpha(l)|_{\mathbb{R}^{n - |I|}} \) is also trivial. Hence, \( \alpha(l) \) is trivial, and we conclude that
\[
l \in \Ker(\alpha)^0 := \left\{ g \in \Ker(\alpha) \;\middle|\; g \text{ is elliptic} \right\},
\]
as desired.
\end{proof}

\begin{lemma}
\label{lem::stab_H_sigma}
Let \( \eta_+, \eta_- \) be two opposite simplices in the ideal boundary of \( \Delta_{G_k} \) such that \( \theta(\eta_+) = \eta_- \). Then
\[
(H_k)_{\eta_-} = (H_k)_{\eta_-, \eta_+} = (H_k)_{\eta_+} = H_k \cap M_{\eta_-, \eta_+},
\]
where \( M_{\eta_-, \eta_+} := P_{\eta_+} \cap P_{\eta_-} \).
\end{lemma}

\begin{proof}
Since \( \theta(\eta_\pm) = \eta_\mp \), it follows from the definitions that \( \theta(P_{\eta_\pm}) = P_{\eta_\mp} \).  Given that \( \theta(h) = h \) for every \( h \in H_k \), and \( \theta(\eta_\pm) = \eta_\mp \), we compute:
\[
\theta(h(\eta_\pm)) = \theta(h)(\theta(\eta_\pm)) = h(\eta_\mp).
\]
Now, if \( h \in (H_k)_{\eta_\pm} \), then \( h(\eta_\pm) = \eta_\pm \), and hence
\[
\theta(h(\eta_\pm)) = \theta(\eta_\pm) = \eta_\mp = h(\eta_\mp),
\]
which implies \( h \in (H_k)_{\eta_\mp} \). Therefore, \( (H_k)_{\eta_\pm} \subseteq (H_k)_{\eta_\mp} \), and by symmetry, equality holds. The lemma follows.
\end{proof}

The conclusion of Lemma~\ref{lem::stab_H_sigma} does not hold if \( \eta_+ \) and \( \eta_- \) are opposite ideal simplices for which \( \theta(\eta_\pm) \neq \eta_\mp \). In such cases, we do \textbf{not} have
\[
(H_k)_{\eta_-} = (H_k)_{\eta_-, \eta_+} = (H_k)_{\eta_+},
\]
and elements in \( (H_k)_{\eta_-} \) may not commute with \( a_{\ell_j} \) for any index \( \ell_j \).

% The Lemma \ref{lem::stab_H_sigma} is not true if we consider two opposite ideal simplices $\eta_-, \eta_+$ for which $\theta(\eta_\pm)$ does not equal $\eta_\mp$, i.e., we will not have $(H_k)_{\eta_-}= (H_k)_{\eta_-,\eta_+}= (H_k)_{\eta_+}$, and \red{elements in $(H_k)_{\eta_-}$ might not commute with $a_{\ell_j}$, for any index $\ell_j$.}

\section{Construction of unipotent elements}
\label{sec::const_unipot}

In this section, we work under general assumptions where \( A \leq G_k \) is a subgroup satisfying the hypothesis \textbf{(Hyp(\( \theta, k \))-split)}. We assume that the apartment \( \mathcal{A} := \mathcal{A}_{\theta,k} \) is \( \theta \)-stable, i.e., \( \theta(\mathcal{A}) = \mathcal{A} \). Let \( \{a_\ell\}_{\ell \geq 1} \subset A \) be a sequence of hyperbolic elements satisfying the condition \textbf{(Same \( \sigma_\pm \))}, and define the Chabauty limit group
\[
L := \limch_{\ell \to \infty} a_\ell H_k a_\ell^{-1} \leq G_k.
\]

Retaining the notation from \textbf{(Hyp(\( \theta, k \))-split)} and \textbf{(Same \( \sigma_\pm \))}, let \( \sigma'_\pm \) denote the unique ideal simplices in the visual boundary \( \partial \Delta_{G_k} \) that contain the endpoints \( \xi_\pm \) in their interiors. Note that \( \sigma'_+ \) is opposite to \( \sigma'_- \).

In addition, we assume the following refinement:

\vspace{1em}
\noindent\textbf{(Hyp \( \sigma'_\pm = \sigma_\pm \)):} We have \( \sigma'_\pm = \sigma_\pm \); that is, the endpoints \( \xi_\pm \) lie in the interior of the simplices \( \sigma_\pm \).

Under these conditions, we construct nontrivial unipotent elements in the Chabauty limit group \( L \) by leveraging the geometry of the Bruhat--Tits building and the Moufang property. These constructions demonstrate that \( L \) acts transitively on the set of ideal simplices in \( \partial \Delta_{G_k} \) that are opposite to a fixed simplex.

\subsection{$L$ acts transitively on $\Opp(\sigma_+)$}
\label{subsec::L_trans_sigma_opp}

According to the definitions given in \cite[Section 4, p. 13, Section 13, p. 55]{HW93}, a parabolic subgroup $P$ of $G$ is called \textbf{$\theta$-split} if $P$ and $\theta(P)$ are opposite parabolic subgroups, and a minimal parabolic \( k \)-subgroup \( P \) of \( G \) is called \textbf{quasi-\( \theta \)-split} if \( P \) is contained in a minimal \( \theta \)-split parabolic \( k \)-subgroup of \( G \). By \cite[Proposition 9.2]{HW93}, \( P \) is quasi-\( \theta \)-split if and only if \( HP \) is open in \( G \). Furthermore, by \cite[Proposition 13.4]{HW93}, for a minimal parabolic \( k \)-subgroup \( P \) of \( G \), the following equivalence holds:
\[
P \text{ is quasi-\( \theta \)-split} \quad \Longleftrightarrow \quad H_k P_k \text{ is open in } G_k.
\]

\begin{proposition}[Version of Proposition 3.3 from \cite{PlaRa}]
\label{prop::open_neigh_compact}
Let \( G', G \) be two connected linear algebraic groups defined over \( k \). Suppose there exists a dominant \( k \)-morphism \( f: G' \to G \), and a point \( x \in G'_k \) such that the differential \( d_x f : T_x G' \to T_{f(x)} G \) is surjective. Then the induced map \( f_k : G'_k \to G_k \) is open at \( x \). Moreover, there exists a Zariski-open subset \( U \subset G' \) such that \( f_k \) is open at every point of \( U_k \).
\end{proposition}

We apply Proposition~\ref{prop::open_neigh_compact} in the setting where \( G \) is a connected linear reductive algebraic group defined over \( k \), that is equipped with a \( k \)-involution \( \theta \).

\begin{proposition}
\label{prop::fixed_H_open}
Let \( G \) be a connected reductive linear algebraic group defined over \( k \), equipped with a \( k \)-involution \( \theta \). Let \( H^0 \) denote the connected component of the identity in the fixed point subgroup of \( G \) under \( \theta \). Let \( P \) be a parabolic \( k \)-subgroup of \( G \) that is \( \theta \)-split, or a minimal parabolic \( k \)-subgroup that is quasi-\( \theta \)-split. Then there exist compact-open subgroups \( K_1 \subseteq H_k \) and \( K_2 \subseteq P_k \) such that \( K_1 K_2 \) is open in \( G_k \).
\end{proposition}

\begin{proof} Let \( P \) be a parabolic \( k \)-subgroup of \( G \) that is \( \theta \)-split.  By Theorem~\ref{thm::Vust_open_split} (or Proposition 9.2 from \cite{HW93}) and its proof, the map
\[
\mu : H^0 \times P \to G, \quad (h, p) \mapsto hp
\]
is dominant, meaning that \( H^0 P \) is Zariski-dense in \( G \), and the differential \( d_e \mu \) is surjective at the identity element \( e \in H^0 \times P \). Since the identity element of an algebraic group over a field \( k \) is a \( k \)-point, we can apply Proposition~\ref{prop::open_neigh_compact} to the connected linear algebraic groups \( G' := H^0 \times P \) and \( G \), and the \( k \)-morphism \( \mu \).

It follows that the induced map
\[
\mu_k : H^0_k \times P_k \to G_k
\]
is open at the identity \( e \in H^0_k \times P_k \). 

If \( P \) is a minimal parabolic \( k \)-subgroup that is quasi-\( \theta \)-split, then by \cite[Proposition 13.4]{HW93}, the product \( H_k P_k \) is open in \( G_k \). Thus, in both cases -- whether \( P \) is \( \theta \)-split or quasi-\( \theta \)-split -- we deduce the following.
Since \( H^0_k \) and \( P_k \) are equipped with the compact-open topology inherited from \( G_k \), and \( H^0_k \times P_k \) carries the product topology, it follows that there exist compact-open subgroups \( K_1 \leq H_k \) and \( K_2 \leq P_k \) such that \( K_1 K_2 \) is open in \( G_k \).

Finally, note that \( G' := H^0 \times P \) is connected because, over an algebraically closed field, the tensor product of integral domains is again an integral domain.
\end{proof}

% The product of V and W is defined as the algebraic set V × W = V( f1,..., fN, g1,..., gM) in A^n+m. 
%over an algebraically closed field, the tensor product of integral domains is an integral domain. Let K be an algebraically closed field and A, B two K-algebras which are integral domains. Then A⊗KB is an integral domain. Let x,x′∈A⊗KB such that xx′=0. Write x=∑ai⊗bi and x′=∑a′j⊗b′j. By taking minimal representations (as sums of monomial tensors) for x and y one can assume that (ai), (a′j), (bi) and (b′j) are linearly independent over K. Now considering A′ the K-subalgebra of A generated by (ai) and (a′j), and analogously B′ the K-subalgebra of B generated by (bi) and (b′j) we reduce the problem to the affine case that is proved in Milne's book. 

%http://tomlr.free.fr/Math%E9matiques/Fichiers%20Claude/Auteurs/aaaDivers/Milne%20Js%20Algebraic%20Geometry%20(2003,%20206S).pdf page 60, Proposition 3.16

%Since simple points always exist in a variety (see Chapter 2, Section 2.4.3 from \cite{PlaRa}), and since groups are homogeneous, one gets that an algebraic group is always smooth, meaning that every of its points is simple (see page 98 from \cite{PlaRa}). Thus a linear algebraic group carries a structure of an analytic variety. 
%
%
%Write also about the algebraic $k$-varieties $X$ and their $k$-points $X_k$ and their associated analytic manifolds $X_{an}$. 

Fix, for what follows, a minimal parabolic \( k \)-subgroup \( P \) of \( G \) that is quasi-\( \theta \)-split. Since \( P \) is minimal, there exists an ideal chamber \( c_- \) in the boundary \( \partial \Delta_{G_k} \) of the Bruhat–Tits building \( \Delta_{G_k} \) of \( G_k \) such that \( P_k = \Stab_{G_k}(c_-) \), i.e., $P_k$ is the Borel subgroup of $G_k$ corresponding to the ideal chamber $c_-$.

By definition, the canonical projection map \( p : G_k \to G_k / P_k \) is continuous and open. By Proposition~\ref{prop::fixed_H_open}, there exist compact-open subgroups \( K_1 \subseteq H_k \) and \( K_2 \subseteq P_k \) such that \( K_1 K_2 \) is open in \( G_k \). It follows that \( p(K_1 K_2) \) is open in \( G_k / P_k \), which is equivalent to saying that \( K_1 P_k / P_k \cong K_1 c_- \) is open in the $G_k$-orbit \( G_k c_- \) of $c_-$.

Therefore, there exists a standard open neighborhood \( V_{c_-} \subset  \Ch(\partial \Delta_{G_k}) \) of \( c_- \) in the cone topology of \( \Ch(\partial \Delta_{G_k}) \), such that the compact-open subgroup \( K_1 \subseteq H_k \) acts transitively on \( V_{c_-} \). That is, for each \( c \in V_{c_-} \), there exists an element \( k \in K_1 \) such that \( k(c_-) = c \).

\medskip
In addition to the notation and assumptions introduced at the beginning of Section~\ref{sec::const_unipot} -- namely, \textbf{(Hyp(\( \theta, k \))-split)}, \textbf{(Same \( \sigma_\pm \))}, \textbf{(Hyp \( \sigma'_\pm = \sigma_\pm \))}, and \( \mathcal{A} := \mathcal{A}_{\theta,k} \) is \( \theta \)-stable -- we also assume the following:

\vspace{1em}
\noindent\textbf{(Hyp \( \sigma_- \subseteq c_- \in \Ch(\partial \mathcal{A}) \)):} The chamber \( c_- \) is an ideal chamber of \( \mathcal{A} \) that contains the simplex \( \sigma_- \).

\vspace{1em}
With this setup, let \( c_+ \) be the ideal chamber in \( \mathcal{A} \) opposite to \( c_- \). Then \( \sigma_+ \) is a subsimplex of \( c_+ \). Choose a point \( x \in \mathcal{A} \), and let \( F \) be the unique minimal flat in \( \mathcal{A} \) that passes through \( x \) and has \( \sigma_\pm \) in its ideal boundary.  Consider also the bi-infinite geodesic line in \( \mathcal{A} \) passing through \( x \), with endpoints given by the barycenters \( \eta_+ \) and \( \eta_- \) of \( \sigma_+ \) and \( \sigma_- \), respectively. Denote this line by \( (\eta_-, x, \eta_+) \). In particular, we have \( (\eta_-, x, \eta_+) \subset F \).

\vspace{1em}
Recall also the definition:

\begin{definition}[Stars in apartments and \( \sigma \)-cones]
\label{def::st_cone_N}
Let \( \sigma \) be an ideal simplex in \( \partial \Delta_{G_k} \), and let \( \mathcal{B} \) be an apartment of \( \Delta_{G_k} \) such that \( \sigma \subset \partial \mathcal{B} \). Let \( x \in \mathcal{B} \) be a point. We define
\[
\st(\sigma, \partial \mathcal{B}) := \{ c \in \Ch(\partial \mathcal{B}) \mid \sigma \text{ lies in the closure of } c \},
\]
and refer to it as the \textbf{\( \partial \mathcal{B} \)-star of \( \sigma \)}.  The associated \textbf{\( \sigma \)-cone in \( \mathcal{B} \) with base point \( x \in \mathcal{B} \)} is defined by
\[
Q(x, \sigma, \mathcal{B}) := \bigcup_{c \in \st(\sigma, \partial \mathcal{B})} Q(x, c),
\]
where \( Q(x, c) \subset \mathcal{B} \) is the closed Weyl sector emanating from \( x \) and pointing toward the ideal chamber \( c \).  A \textbf{\( \sigma \)-cone in \( \mathcal{B} \)} is any set of the form \( Q(x, \sigma, \mathcal{B}) \) for some \( x \in \mathcal{B} \).
\end{definition}

\textbf{Construction:}
\begin{enumerate}
\item[]
Let \( \st(\sigma_+, \partial \mathcal{A}) \) be the \( \partial \mathcal{A} \)-star of \( \sigma_+ \), and let \( Q(x, \sigma_+, \mathcal{A}) \) be the \( \sigma_+ \)-cone in \( \mathcal{A} \) based at \( x \). Note that the geodesic ray \( [x, \eta_+) \) lies in the interior of \( Q(x, \sigma_+, \mathcal{A}) \).

Let \( \sigma \in \Opp(\sigma_+) \) be an ideal simplex opposite to \( \sigma_+ \).  By \cite[Lemma 3.4]{Cio24}, consider the unique \( \st(\sigma, \partial \mathcal{A}_{\sigma, \sigma_+}) \) in the \( \partial \Delta_{G_k} \)-star of \( \sigma \), such that \( \st(\sigma, \partial \mathcal{A}_{\sigma, \sigma_+}) \) and \( \st(\sigma_+, \partial \mathcal{A}) \) determine a unique apartment \( \mathcal{A}_{\sigma, \sigma_+} \subset \Delta_{G_k} \) containing both \( \sigma \) and \( \sigma_+ \) at infinity.  Then \( \mathcal{A}_{\sigma, \sigma_+} \cap \mathcal{A} \) is a geodesically convex subset and contains a \( \sigma_+ \)-subcone of \( Q(x, \sigma_+, \mathcal{A}) \).

We can now extend the bi-infinite geodesic line \( (\eta_-, x, \eta_+) \subset \mathcal{A} \) to a bi-infinite geodesic line \( (\eta, \eta_+) \subset \mathcal{A}_{\sigma, \sigma_+} \), such that
\[
(\eta_-, x, \eta_+) \cap (\eta, \eta_+) = [x_\sigma, \eta_+),
\]
where \( x_\sigma \in \mathcal{A}_{\sigma, \sigma_+} \cap \mathcal{A} \) is a point. The geodesic ray \( (\eta, x_\sigma) \subset \mathcal{A}_{\sigma, \sigma_+} \) has empty intersection with \( \mathcal{A} \).

Note that \( \eta \) is the barycenter of the simplex \( \sigma \), $\eta$ is opposite to \( \eta_+ \), and the point \( x_\sigma \) is uniquely determined by \( \eta \), the fixed point \( x \), and the endpoint \( \eta_+ \).

Moreover, there exists a unique ideal chamber \( c_\sigma \in \st(\sigma, \partial \mathcal{A}_{\sigma, \sigma_+}) \) that is opposite to \( c_+ \). Then, by \cite[Proposition 4.8]{Cio24}, applied  to the sequence $\{a_\ell\}_{\ell \geq 1}$ and apartment $\mathcal{A}$, we have
\[
\lim_{\ell \to \infty} a_\ell^{-1}(c_\sigma) = c_-.
\]
\end{enumerate}

\begin{proposition}
\label{prop::exist_l_sigma}
Retain the notation and assumptions from \textbf{(Hyp(\( \theta, k \))-split)}, \textbf{(Same \( \sigma_\pm \))}, \textbf{(Hyp \( \sigma'_\pm = \sigma_\pm \))},  \( \mathcal{A} := \mathcal{A}_{\theta,k} \) is \( \theta \)-stable, and \textbf{(Hyp \( \sigma_- \subseteq c_- \in \Ch(\partial \mathcal{A}) \))}.

Let \( \sigma \in \Opp(\sigma_+) \) with \( \sigma \neq \sigma_- \). Then there exists an elliptic and type-preserving element \( l \in L  \in P_{\sigma_+}\) such that \( l(\sigma_-) = \sigma \).
\end{proposition}

\begin{proof}
 From the \textbf{Construction} above, consider the corresponding point \( x_\sigma \) and the chamber \( c_\sigma \in \Opp(c_+) \) associated with \( \sigma \). Let \( F_{\sigma} \subset \mathcal{A}_{\sigma,\sigma_+} \) denote the unique minimal flat that passes through \( x_\sigma \) and whose ideal boundary contains both \( \sigma \) and \( \sigma_+ \).

Moreover, there exists a constant \( N > 0 \) such that for every \( \ell > N \), we have \( a_{\ell}^{-1}(c_\sigma) \in V_{c_-} \). Since \( K_1 \) is compact, all its elements are elliptic. For \( K_1 \) sufficiently small, every element of \( K_1 \) acts as a type-preserving automorphism on \( \Delta_{G_k} \cup \partial \Delta_{G_k} \).  Therefore, for our chosen \( c_\sigma \in \Opp(c_+) \), we can select a sequence of elliptic elements \( \{h_{\ell}\}_{\ell} \subset K_1 \leq H_k \) such that
\[
h_{\ell}(c_-) = a_{\ell}^{-1}(c_\sigma), \quad \text{for all } \ell > N.
\]
By the hypothesis that $\mathcal{A} \subseteq \Min(a_{\ell})$, we have $a_{\ell}(c_\pm)= c_\pm$ pointwise, for every $\ell\in \NN$. It follows that
\[
a_{\ell} h_{\ell} a_{\ell}^{-1}(c_-) = c_\sigma.
\]
In particular, we obtain
\[
a_{\ell} h_{\ell} a_{\ell}^{-1}(\sigma_-) = \sigma, \quad \text{for all } \ell > N.
\]
\textit{Note:} Since the action of $h_\ell$ is type-preserving, the type of \( \sigma_- \) remains unchanged under conjugation by \( a_{\ell} \), ensuring that the image is indeed \( \sigma \).

We claim that the sequence \( \{a_{\ell} h_{\ell} a_{\ell}^{-1}\}_{\ell \geq 1} \) admits a convergent subsequence indexed by a strictly increasing sequence \( \{\ell_j\}_{j \geq 1} \), such that
\[
\lim_{j \to \infty} a_{\ell_j} h_{\ell_j} a_{\ell_j}^{-1} = l \in L,
\]
where the limit element \( l \) satisfies the desired properties: it is elliptic, type-preserving, and maps \( \sigma_- \) to \( \sigma \), i.e., \( l(\sigma_-) = \sigma \).

To establish the claim, it suffices to show that the sequence \( \{a_{\ell} h_{\ell} a_{\ell}^{-1}\}_{\ell \geq 1} \) admits a convergent subsequence in \( G_k \). The rest of the claim then follows from Proposition~\ref{prop::chabauty_conv}.

Recall that \( K_1 \) is a compact-open subgroup of \( H_k \), where $H_k$ is equipped with the compact-open topology inherited from \( G_k \). In particular, we can choose \( K_1 \) sufficiently small so that it fixes the chosen point \( x_\sigma \in \mathcal{A} \) from the \textbf{Construction}.

Observe that both flats \( F \) and \( F_\sigma \) contain the point \( x_\sigma \) and share the same simplex at infinity \( \sigma_+ \). Therefore, the intersection \( F \cap F_\sigma \) is a geodesically convex subset of \( F \), and hence also of \( \mathcal{A} \), since \( F \subset \mathcal{A} \). Moreover, this intersection contains \( x_\sigma \) in its boundary.

Moreover, since \( F \subseteq \mathcal{A} \subseteq \Min(a_\ell) \) with $\sigma_\pm \subset \partial F$, the flat \( F \) contains a translation axis of \( a_\ell \), for every \( \ell \geq 1 \). In particular, the point \( a_\ell^{-1}(x_\sigma) \) remains in \( F \) for all \( \ell \geq 1 \).

Furthermore, consider the cone \( Q(x_\sigma, \sigma_-) \subset F \), based at \( x_\sigma \) and directed toward the ideal simplex \( \sigma_- \subset \partial_\infty F \). By the choice of the sequence \( \{a_\ell\}_{\ell \geq 1} \), the geodesic segment \( (x_\sigma, a_\ell^{-1}(x_\sigma)] \) lies in the interior of the cone \( Q(x_\sigma, \sigma_-) \).

In order for \( h_\ell \) to send \( \sigma_- \subset c_- \) to \( a_\ell^{-1}(\sigma) \subset a_\ell^{-1}(c_\sigma) \), it must -- being elliptic, simplicial, and fixing \( x_\sigma \) -- pointwise fix the geodesically convex intersection of the flat \( a_\ell^{-1}(F_\sigma) \) with the cone \( Q(x_\sigma, \sigma_-) \). This ensures that the cone \( Q(x_\sigma, \sigma_-) \) is mapped by $h_{\ell}$ into a cone lying within the flat \( a_\ell^{-1}(F_\sigma) \), thereby sending \( \sigma_- \) to \( a_\ell^{-1}(\sigma) \).

In particular, the point \( a_\ell^{-1}(x_\sigma) \) lies on the boundary of this convex intersection and is fixed by \( h_\ell \). Therefore, applying the hyperbolic element \( a_\ell \) back yields
\[
a_\ell h_\ell a_\ell^{-1}(x_\sigma) = x_\sigma.
\]
This shows that the sequence \( \{a_\ell h_\ell a_\ell^{-1}\}_{\ell \geq 1} \) is bounded by its action on \( x_\sigma \), and hence we can extract a convergent subsequence indexed by a strictly increasing sequence \( \{\ell_j\}_{j \geq 1} \). The claim as well as the proposition follow.
\end{proof}

\subsection{On the unipotent part of $L$}
 
Let $\eta$ be an ideal simplex in the boundary $\partial \Delta_{G_k}$ of the building $\Delta_{G_k}$. The \textbf{star of $\eta$ in $\partial \Delta_{G_k}$}, denoted by $\mathrm{St}(\eta, \partial \Delta_{G_k})$, is defined as the set of all ideal chambers in $\partial \Delta_{G_k}$ that contain $\eta$ as a subsimplex.

Let $\alpha$ be a root of $\partial \Delta_{G_k}$, i.e., a half-apartment in $\partial \Delta_{G_k}$. The associated \textbf{root group} $U_{\alpha}(G_k)$ of $G_k$ is defined as the set of elements $g \in G_k$ satisfying the following conditions:
\begin{itemize}
    \item[(a)] $g$ fixes the half-apartment $\alpha$ pointwise
    \item[(b)] $g$ fixes $\mathrm{St}(P, \partial \Delta_{G_k})$ pointwise for every panel $P$ in the interior of $\alpha$, i.e., within the half-apartment $\alpha$ excluding its boundary wall.
\end{itemize}

Let $\mathcal{A}(\alpha, \partial \Delta_{G_k})$ denote the set of all apartments in $\partial \Delta_{G_k}$ that contain the root $\alpha$ as a half-apartment.  We say that the group $G_k$ is \textbf{Moufang} if, for every root $\alpha$ of $\partial \Delta_{G_k}$, the corresponding root group $U_{\alpha}(G_k)$ acts transitively on the set $\mathcal{A}(\alpha, \partial \Delta_{G_k})$.  For further details on root groups and the Moufang property, see \cite[Chapter 7]{AB} or \cite{Cio_M}.

\begin{proposition}
\label{prop::find_unipotent}
Retain the notation and assumptions from \textbf{(Hyp(\( \theta, k \))-split)}, \textbf{(Same \( \sigma_\pm \))}, \textbf{(Hyp \( \sigma'_\pm = \sigma_\pm \))},  \( \mathcal{A} := \mathcal{A}_{\theta,k} \) is \( \theta \)-stable, and \textbf{(Hyp \( \sigma_- \subseteq c_- \in \Ch(\partial \mathcal{A}) \))}. Suppose that $G_k$ is Moufang and   $M_{\sigma_-,\sigma_+}(k) \leq Z_{G_k}(A)$.

Let \( u \in U_{\sigma_+}(k) \) be a non-trivial unipotent element in the unipotent radical \( U_{\sigma_+}(k) \) of the parabolic subgroup \( P_{\sigma_+}(k) \). Then there exists a subsequence \( \{h_{\ell_j}\}_{j \geq 1} \subset H_k \), indexed by a strictly increasing sequence \( \{\ell_j\}_{j \geq 1} \subset \{\ell\}_{\ell \geq 1} \), such that
\[
\lim_{j \to \infty} a_{\ell_j} h_{\ell_j} a_{\ell_j}^{-1} = u m \in L\leq P_{\sigma_+}(k),
\]
for some \( m \in \Ker(\alpha)^0 \leq M_{\sigma_-,\sigma_+}(k) \), where \( \alpha \) is the map from Proposition~\ref{prop::res_building_str_tran}, and $
\Ker(\alpha)^0 := \{ g \in \Ker(\alpha) \mid g \text{ is elliptic} \}$.
\end{proposition}

\begin{proof}
The assumption that $G_k$ is Moufang implies that the spherical building $\partial \Delta_{G_k}$ is Moufang in the sense of \cite[Definition 7.27]{AB}. Recall that, as stated in Section~\ref{section::auto_invol_building}, the buildings $\Delta_{G_k}$ and $\partial \Delta_{G_k}$ are assumed to be irreducible.

By \cite[Theorem 6.18 and its proof]{Ro}, \cite[Lemma 7.25(3)]{AB}, and \cite[Corollary 2.3]{Cio_M}, the unipotent radical \( U_{\sigma_+}(k) \) of the parabolic subgroup \( P_{\sigma_+}(k) \) acts simply transitively on the set
\[
\Opp(\sigma_+) \subset \Delta_{G_k}
\]
of all ideal simplices opposite to \( \sigma_+ \). Consequently, there exists a unique simplex \( \sigma \in \Opp(\sigma_+) \) such that \( u(\sigma_-) = \sigma \). Note that \( \sigma \neq \sigma_- \).
Let \( c_{\sigma} \) be the unique ideal chamber in \( \partial \Delta_{G_k} \) such that \( \sigma \subset c_{\sigma} \) and \( c_{\sigma} \) is opposite to \( c_+ \). Then, by Proposition~\ref{prop::exist_l_sigma}, we obtain an elliptic and type-preserving element
\[
l = \lim_{j \to \infty} a_{\ell_j} h_{\ell_j} a_{\ell_j}^{-1} \in L
\]
such that \( l(\sigma_-) = \sigma \). We aim to show that
\[
l \in u\cdot \Ker(\alpha)^0 \cdot (H_k \cap M_{\sigma_-,\sigma_+}),
\]
and we can simplify this further using Corollary~\ref{cor::find_elements_H} applied to $l$. Thus, from now on, we consider
\[
l = \lim_{j \to \infty} a_{\ell_j} h_{\ell_j} a_{\ell_j}^{-1} \in L \text{ type-preserving}, \quad \text{with } \lim_{j \to \infty} h_{\ell_j} = e, \quad \text{and } l(\sigma_-) = \sigma,
\]
and we wish to prove that \( l \in u \Ker(\alpha)^0 \).

Indeed, referring to the \textbf{Construction} and its notation, consider the geodesic ray \( [x_{\sigma}, \eta_+) \subset \mathcal{A} \), which is pointwise fixed by \( l \), and the flat \( F \subset \mathcal{A} \) containing the bi-infinite geodesic line \( (\eta_-, x_\sigma, \eta_+) \). The visual boundary of \( F \) contains \( \sigma_- \) and \( \sigma_+ \), and \( F \) is the unique minimal flat in \( \mathcal{A} \) that contains \( x_\sigma \) and has \( \sigma_\pm \) in its ideal boundary.

Moreover, since \( l(\sigma_+) = \sigma_+ \) pointwise and \( l \) is elliptic and type-preserving, it follows that \( l \) fixes pointwise the cone \( Q(x_\alpha, \sigma_+, F) \) in the flat \( F \subset \mathcal{A} \), with basepoint \( x_\alpha \) and pointing toward \( \sigma_+ \).

Consider the sequence of points \( \{a_{\ell_j}(x_\alpha)\} \subset Q(x_\alpha, \sigma_+, F) \). These points lie in the interior of the cone \( Q(x_\alpha, \sigma_+, F) \), since each \( a_\ell \) has a translation axis contained in the flat \( F \subseteq \mathcal{A} \subseteq \Min(a_\ell) \), for every \( \ell \geq 1 \). Therefore,
\[
l(a_{\ell_j}(x_\alpha)) = a_{\ell_j}(x_\alpha), \quad \text{for all } j \geq 1,
\]
and consequently,
\[
a_{\ell_j}^{-1} l a_{\ell_j}(x_\alpha) = x_\alpha, \quad \text{for all } j \geq 1.
\]
This implies that the sequence \( \{a_{\ell_j}^{-1} l a_{\ell_j}\}_{j \geq 1} \) is bounded in \( G_k \), and hence admits a convergent subsequence. For simplicity, we continue to index this subsequence by \( \{\ell_j\}_j \).

\medskip
One can show that the limit $m := \lim\limits_{j \to \infty} a_{\ell_j}^{-1} l a_{\ell_j}$
is an element of \( M_{\sigma_-, \sigma_+}(k) \), and is moreover elliptic and type-preserving, since \( l \) is. Indeed, for every \( j \geq 1 \),
\[
a_{\ell_j}^{-1} l a_{\ell_j}(\sigma_+) = a_{\ell_j}^{-1} l(\sigma_+) = a_{\ell_j}^{-1}(\sigma_+) = \sigma_+,
\]
so \( m(\sigma_+) = \sigma_+ \). For \( \sigma_- \), observe that
\[
\lim_{j \to \infty} a_{\ell_j}^{-1} l a_{\ell_j}(a_{\ell_j}^{-1}(x_\alpha)) = \lim_{j \to \infty} a_{\ell_j}^{-1}(x_\alpha) = \xi_- = m(\xi_-),
\]
where \( \xi_- \) lies in the interior of \( \sigma_- \). It follows that \( m \in M_{\sigma_-, \sigma_+}(k) \).

\medskip
In order to conclude that \( l \in U_{\sigma_+}(k) \), we need to show that \( m \) is the identity element in \( G_k \). However, we will only be able to show that \( m \in \Ker(\alpha)^0 \leq M_{\sigma_-,\sigma_+}(k) \), meaning that \( m \) acts trivially on the sub-building \( \Delta(\sigma_+, \sigma_-) \subset \Delta_{G_k} \). This does not imply that \( m \) acts trivially on the entire building \( \Delta_{G_k} \), nor that \( m \) is the identity element in \( G_k \).

Since \( M_{\sigma_-, \sigma_+}(k) \leq Z_{G_k}(A) \), it follows that \( m \) commutes with \( a_\ell \) for every \( \ell \). Therefore, we obtain
\[
l m^{-1} \in U_{\sigma_+}(k),
\]
with the property that
\[
l m^{-1}(\sigma_-) = \sigma = u(\sigma_-).
\]
Here, we have used the characterization of unipotent radicals as contraction subgroups to deduce that $l m^{-1} \in U_{\sigma_+}(k)$  (see, for example, \cite[Appendix B]{Cio24}).

Because \( U_{\sigma_+}(k) \) acts simply transitively on \( \Opp(\sigma_+) \), it follows that \( l m^{-1} = u \), and hence \( l = u m \), as claimed in the statement of the proposition.

To show that \( m \in \Ker(\alpha)^0 \), it remains to compare the actions of \( l \), \( m \), and the identity element \( e \) on the locally finite thick affine building \( \Delta(\sigma_+, \sigma_-) \), or equivalently on the building \( \Delta_I \), as described in Propositions~\ref{prop::res_building} and~\ref{prop::res_building_str_tran}.

First, observe that by the definition \( m = \lim_{j \to \infty} a_{\ell_j}^{-1} l a_{\ell_j} \), and due to the compatible compact-open topologies on \( G_k \) and \( \Delta_{G_k} \), for any point \( x \in \Delta_{G_k} \) and radius \( r > 0 \), there exists an index \( N_{x,r} > 0 \) such that
\[
m(B(x, r)) = a_{\ell_j}^{-1} l a_{\ell_j}(B(x, r)) \text{ pointwise, for every } \ell_j \geq N_{x,r}.
\] Since \( m \) commutes with \( a_{\ell_j} \), we have:
\[
a_{\ell_j} m(B(x, r)) = m(B(a_{\ell_j}(x), r)) = l(B(a_{\ell_j}(x), r)) \text{ pointwise, for every } \ell_j \geq N_{x,r}.
\]
This shows that the action of \( l \) on the balls \( B(a_{\ell_j}(x), r) \) coincides with that of \( m \), for sufficiently large \( \ell_j \).

For simplicity, we now take \( x = x_\sigma \). Since \( l(a_{\ell_j}(x_\sigma)) = a_{\ell_j}(x_\sigma) \), it follows that \( m(a_{\ell_j}(x_\sigma)) = a_{\ell_j}(x_\sigma) \) for all \( \ell_j \), which in particular implies \( m(x_\sigma) = x_\sigma \), as \( m \) commutes with \( a_{\ell_j} \).  Moreover, \( m \) fixes the flat \( F \) pointwise, since \( m(x_\sigma) = x_\sigma \in F \) and \( m(\sigma_\pm) = \sigma_\pm \), and $m$ is type-preserving and simplicial.
 
 \medskip
%Given that, let us study the action of $m$ on $\Delta(\sigma_+, \sigma_-) \cong F \times \Delta_{I}$, and thus on $\Delta_I$. First notice that by the choice of $\{a_\ell\}_{\ell \in \NN}$, we have $a_\ell(F)=F$ setwise, for every $\ell\geq 1$. So if we restrict a ball $B(a_{\ell_j}(x_\sigma),r)$ to a transversal section corresponding to $\Delta_I$ through the point $a_{\ell_j}(x_\sigma)$, we will get the same action of $m$, independently from $a_{\ell_j}$, as $m$ fixes pointwise $F$.  Therefore, we can deduce that for $\ell_j$ large enough, the action of $l$ on the balls $B(a_{\ell_j}(x_\sigma),r)$ restricted to a transversal section corresponding to $\Delta_I$ through the point $a_{\ell_j}(x_\sigma)$ will be all the same as the one of $m$ on the balls $B(a_{\ell_j}(x_\sigma),r)$, and thus on $B(x_\sigma, r)$, restricted to a transversal section corresponding to $\Delta_I$ through the point $x_\sigma$, and so constant. Thus for every $\ell_j \geq N_{x_\sigma, r}$  
%$$l(\proj_{\Delta_I}(B(a_{\ell_j}(x_\sigma),r))) = m(\proj_{\Delta_I}(B(x_\sigma, r))) \text{ pointwise},$$ 
%where $\proj_{\Delta_I}(B(y,R)$ is $B(y,R) \cap \Delta(\sigma_+,\sigma_-)$ projected to the factor $\Delta_I$ of $\Delta(\sigma_+,\sigma_-)$.
% 

Given this setup, let us study the action of \( m \) on the sub-building \( \Delta(\sigma_+, \sigma_-) \cong F \times \Delta_I \), and in particular on the factor \( \Delta_I \).  First, note that by the choice of the sequence \( \{a_\ell\}_{\ell \in \mathbb{N}} \), we have \( a_\ell(F) = F \) setwise for every \( \ell \geq 1 \). Therefore, if we restrict a ball \( B(a_{\ell_j}(x_\sigma), r) \) to a transversal section corresponding to \( \Delta_I \) through the point \( a_{\ell_j}(x_\sigma) \), the action of \( m \) on this section will be independent of \( a_{\ell_j} \), since \( m \) fixes \( F \) pointwise.

It follows that, for sufficiently large \( \ell_j \), the action of \( l \) on the balls \( B(a_{\ell_j}(x_\sigma), r) \), when restricted to the transversal section corresponding to \( \Delta_I \) through the point \( a_{\ell_j}(x_\sigma) \), coincides with the action of \( m \) on the balls \( B(x_\sigma, r) \), restricted to the transversal section through \( x_\sigma \). Hence, the action of $l$ is constant across these sections.

Thus, for every \( \ell_j \geq N_{x_\sigma, r} \), we have:
\[
l\left( \proj_{\Delta_I}(B(a_{\ell_j}(x_\sigma), r)) \right) = m\left( \proj_{\Delta_I}(B(x_\sigma, r)) \right) \quad \text{pointwise},
\]
where \( \proj_{\Delta_I}(B(y, R)) := B(y, R) \cap \Delta(\sigma_+, \sigma_-) \) projected to the factor \( \Delta_I \) of \( \Delta(\sigma_+, \sigma_-) \) through the point $y$.
 
\medskip
We now aim to compare the action of \( l \) on \( \proj_{\Delta_I}(B(a_{\ell_j}(x_\sigma), r)) \) with that of the identity element \( e \) on the same set. To do this, we use the hypothesis that
\[
l = \lim_{j \to \infty} a_{\ell_j} h_{\ell_j} a_{\ell_j}^{-1}, \quad \text{with } \lim_{j \to \infty} h_{\ell_j} = e.
\]
Our approach follows the strategy used in the proof of Proposition~\ref{prop::levi_factors}.  Fix an index \( \ell_{j_0} \geq N_{x_\sigma, r} \). Since \( l = \lim_{j \to \infty} a_{\ell_j} h_{\ell_j} a_{\ell_j}^{-1} \), there exists \( N_{j_0} > 0 \) such that for all \( \ell_j \geq N_{j_0} \), we have:
\[
l(B(a_{\ell_{j_0}}(x_\sigma), r)) = a_{\ell_j} h_{\ell_j} a_{\ell_j}^{-1}(B(a_{\ell_{j_0}}(x_\sigma), r)) \text{ pointwise},
\]
with $l(B(a_{\ell_{j_0}}(x_\sigma), r)) = B(a_{\ell_{j_0}}(x_\sigma), r)$ setwise and
\[
l(B(a_{\ell_{j_0}}(x_\sigma), r) \cap F) = B(a_{\ell_{j_0}}(x_\sigma), r) \cap F) \text{ pointwise}.
\]
 
In particular, we get 
$$a_{\ell_j} h_{\ell_j} a_{\ell_j}^{-1}(B(a_{\ell_{j_0}}(x_\sigma), r)) = B(a_{\ell_{j_0}}(x_\sigma), r) \text{ setwise,}$$
with $a_{\ell_j} h_{\ell_j} a_{\ell_j}^{-1}(B(a_{\ell_{j_0}}(x_\sigma), r) \cap F) = B(a_{\ell_{j_0}}(x_\sigma), r)\cap F$ pointwise, for every  $\ell_j \geq N_{j_0}$. This implies that for every  $\ell_j \geq N_{j_0}$
\begin{equation*}
\begin{split}
h_{\ell_j} a_{\ell_j}^{-1}(B(a_{\ell_{j_0}}(x_\sigma), r)) &= a_{\ell_j}^{-1}(B(a_{\ell_{j_0}}(x_\sigma), r)) \text{ setwise, with}\\
h_{\ell_j} a_{\ell_j}^{-1}(B(a_{\ell_{j_0}}(x_\sigma), r) \cap F) &= a_{\ell_j}^{-1}(B(a_{\ell_{j_0}}(x_\sigma), r)\cap F) \\
&=a_{\ell_j}^{-1}(B(a_{\ell_{j_0}}(x_\sigma), r)) \cap a_{\ell_j}^{-1}(F)\\
& =a_{\ell_j}^{-1}(B(a_{\ell_{j_0}}(x_\sigma), r)) \cap F \text{ pointwise.}\\
\end{split}
\end{equation*} 

We recall that the sequence $\{h_{\ell_j}\}_{j \geq 1}$ converges to $e$ in $H_k$, and hence in $G_k$. Therefore, for any $R > 0$, there exists $M_R > 0$ such that $h_{\ell_j}(B(x_\sigma, R)) = B(x_\sigma, R)$ pointwise for every $\ell_j \geq M_R$. In particular, we can apply Proposition \ref{prop::res_building_str_tran} and obtain
$$
\alpha(h_{\ell_j})\vert_{\Delta_I}(\proj_{\Delta_I}(B(x_\sigma, R))) = \proj_{\Delta_I}(B(x_\sigma, R)) \quad \text{pointwise}
$$
for every $\ell_j \geq M_R$. Thus, as $\ell_j$ increases, the elements $\alpha(h_{\ell_j})$ converge to the identity automorphism of the affine building $\Delta_I$.

Now we consider the action of $h_{\ell_j}$ on $\Delta_I$, viewed from the ball $a_{\ell_j}^{-1}(B(a_{\ell_{j_0}}(x_\sigma), r))$. Since $h_{\ell_j}$ fixes pointwise $F \cap a_{\ell_j}^{-1}(B(a_{\ell_{j_0}}(x_\sigma), r))$, $F \cap B(x_\sigma, R)$, and $B(x_\sigma, R)$ for every $\ell_j \geq \max(M_R, N_{j_0}, N_{x_\sigma, r})$, we take a transversal section corresponding to $\Delta_I$ through a point in $F \cap a_{\ell_j}^{-1}(B(a_{\ell_{j_0}}(x_\sigma), r))$, and another through a point in $F \cap B(x_\sigma, R)$. We observe that the corresponding actions of $\alpha(h_{\ell_j})$ on $\proj_{\Delta_I}(B(x_\sigma, R))$ and $\proj_{\Delta_I}(B(a_{\ell_j}(x_\sigma), r))$ must be the same, and hence trivial, since the one on $B(x_\sigma, R)$ is. Therefore, by choosing $R$ large enough, we conclude that the action of $\alpha(l)$ on a transversal section corresponding to $\Delta_I$ through a point in $F \cap B(a_{\ell_{j_0}}(x_\sigma), r)$ is also trivial. This implies that
$$
l(\proj_{\Delta_I}(B(a_{\ell_j}(x_\sigma), r))) = m(\proj_{\Delta_I}(B(x_\sigma, r))) \text{ pointwise}
$$
is in fact the trivial action. Since this holds for every radius $r > 0$, we deduce that $\alpha(m)\vert_{\Delta_I}$ is trivial. As $m$ fixes $F$ pointwise, it follows that $\alpha(m)\vert_F$ is trivial. Hence, $\alpha(m)$ is trivial, and so $m \in \Ker(\alpha)^0 := \{ g \in \Ker(\alpha) \mid g \text{ is elliptic} \}$, as desired.
 \end{proof}
 
 %Moreover, recalling that a bases for the compact-open topology on $G_k$, and so on $H_k$, is given by the $G_k$-pointwise stabilizers of closed balls around points of $\Delta_{G_k}$ (those are compact since $\Delta_{G_k}$ is locally finite), we can consider the compact-open subgroup $K'_1:= K_1 \cap \Stab_{G_k}(x) \cap H_k \leq H_k$.  

\section{The Main Theorem}

We state and prove the main theorem, which characterizes all nontrivial Chabauty limits of \( H_k \) in \( G_k \). These limits are shown to be $G_k$-conjugate to specific subgroups of 
\[
U_{\sigma_+}(k) \rtimes \Ker(\alpha)^0 \cdot (H_k \cap M_{\sigma_-, \sigma_+}),
\]
where \( \alpha \) is the projection map introduced in Proposition~\ref{prop::res_building_str_tran} and associated with the Levi decomposition $ U_{\sigma_+} \rtimes M_{\sigma_-, \sigma_+}$ of a specific parabolic subgroup $P_{\sigma_+}$, and where $\Ker(\alpha)^0 := \{ g \in \Ker(\alpha) \mid g \text{ is elliptic} \}$.

To establish this result, we first recall several key notions and structural facts. By \cite[Proposition 2.3]{HW93}, there exists a maximal \( k \)-split torus \( T \subset G \) that is \( \theta \)-stable; that is, \( \theta(T) = T \) setwise. For any such \( \theta \)-stable torus \( T \), we have the decomposition
\[
T = T_+ T_-,
\]
where
\begin{align*}
T_+ &:= \{ t \in T \mid \theta(t) = t \}, \\
T_- &:= \{ t \in T \mid \theta(t) = t^{-1} \}.
\end{align*}
Note that this decomposition is not necessarily a direct product, as the intersection \( T_+ \cap T_- \) may be nontrivial. For instance, it may contain elements \( t \in T \) satisfying \( t = t^{-1} \), which is equivalent to \( t^2 = 1 \).

A $k$-torus $A$ in $G$ is called \textbf{$(\theta,k)$-split} if it is $k$-split and satisfies $\theta(x) = x^{-1}$ for every $x \in A$. According to \cite[Proposition 4.3]{HW93}, if $\theta$ acts nontrivially on the isotropic factor of $G$ over $k$, then $G$ admits nontrivial $(\theta,k)$-split tori.

Moreover, we claim the existence of a $\theta$-stable maximal $k$-split torus $T$ in $G$ such that its $T_-$ component has the maximal possible dimension among all $(\theta,k)$-split tori in $G$. Indeed, by \cite[Proposition 4.5(iii)]{HW93} or \cite[Proposition 3.2(a)]{BeOh}, any maximal $k$-split torus of $G$ that contains a maximal $(\theta,k)$-split torus is $\theta$-stable. Now, assuming that $\theta$ acts nontrivially on the isotropic factor of $G$ over $k$, it follows from \cite[Proposition 4.3]{HW93}, see also \cite[Lemma 3.1]{BeOh}, that $[G,G]$, and hence $G$, admits $(\theta,k)$-split tori. Since every $k$-split torus --particularly a $(\theta,k)$-split torus -- is contained in some maximal $k$-split torus of $G$, the claim follows.

By \cite[Section 6.9, Proposition 6.10, and Corollary 6.16]{HW93} (see also \cite[Proposition 3.2(c)]{BeOh}), it is known that there are only finitely many $H_k$-conjugacy classes of $\theta$-stable maximal $k$-split tori in $G$. This finiteness result holds because $k$ is a local field.  Furthermore, by \cite[Proposition 10.3]{HW93}, any two maximal $(\theta,k)$-split tori in $G$ are conjugate under $G_k$. Consequently, there are only finitely many $H_k$-conjugacy classes of maximal $(\theta,k)$-split tori in $G$.  Choose a set $\{A_i \; \vert \; 1\leq i \leq m\}$ of representatives of $H_k$-conjugacy class of maximal  $(\theta,k)$-split tori of $G$ and set $\mathcal{B} = \bigcup\limits_{i=1}^{m}A_i$, and $\mathcal{B}_k = \bigcup\limits_{i=1}^{m}A_i(k)$ where $A_i(k)$ are the $k$-points of the group $A_i$.

We have the following theorem due to Benoist--Oh:
\begin{theorem}[Theorem 1.1 of \cite{BeOh}]
\label{thm::kbh_decom}
There exists a compact subset \( K \subset G_k \) such that
\[
G_k = K \mathcal{B}_k H_k.
\]
\end{theorem}
It is important to note that all the subgroups \( A_i \) for \( 1 \leq i \leq m \) are required -- up to \( H_k \)-conjugacy -- in the above decomposition. Moreover, in general, one cannot choose \( K \) to be a compact subgroup of \( G_k \).

\medskip
In Appendix \ref{appen::A}, and in accordance with the approaches outlined in \cite{HW93} and \cite{BeOh}, we present the majority of the steps involved in establishing the $K \mathcal{B} H$ decomposition.

\medskip
The following theorem is one of the main results of the article. 

\begin{theorem}
\label{thm::main_thm}
Let \( G \) be a connected linear reductive group defined over a non-Archimedean local field \( k \) of characteristic zero. Assume that its associated Bruhat--Tits building \( \Delta_{G_k} \) is irreducible and that \( G_k \) is Moufang. Let \( \theta \) be an involutive \( k \)-automorphism of \( G \), let \( G^{\theta} := \{ h \in G \mid \theta(h) = h \} \leq G \) be the fixed-point subgroup of \( \theta \), and $H = (G^{\theta})^{o}$ be the connected component of $G^{\theta}$.

Then for any Chabauty limit \( L' \) of \( H_k \) in \( G_k \) that is not \( G_k \)-conjugate to \( H_k \), there exist a $(\theta, k)$-split torus \( A \) from \( \mathcal{B} \) and a sequence \( \{a_\ell\}_{\ell \in \mathbb{N}} \subset A_k \) of hyperbolic elements satisfying the assumptions \textbf{(Same \( \sigma_\pm \))} and \textbf{(Hyp \( \sigma'_\pm = \sigma_\pm \))}, such that \( L' \) is \( G_k \)-conjugate to the limit group
\[
L := \limch_{\ell \to \infty} a_\ell H_k a_\ell^{-1} \leq G_k.
\]

In particular:
\begin{enumerate}
\item[1)]
\label{thm::main_1}
\( L \leq U_{\sigma_+}(k) \rtimes \left( \Ker(\alpha)^0 \cdot \left( H_k \cap M_{\sigma_-,\sigma_+} \right) \right) \leq P_{\sigma_+} = U_{\sigma_+} \rtimes M_{\sigma_-,\sigma_+} \), where \( \alpha \) is the map from Proposition~\ref{prop::res_building_str_tran}, and \( \Ker(\alpha)^0 := \{ g \in \Ker(\alpha) \mid g \text{ is elliptic} \} \),
\item[2)]
\( H_k \cap M_{\sigma_-,\sigma_+} \leq L \),
\item[3)] the group \( L \) acts transitively on the set \( \Opp(\sigma_+) \), which consists of all ideal simplices in \( \partial \Delta_{G_k} \) that are opposite to \( \sigma_+ \),
\item[4)]
for every nontrivial unipotent element \( u \in U_{\sigma_+}(k) \), there exists \( m \in \Ker(\alpha)^0 \) such that \( um \in L \).
\end{enumerate}
\end{theorem}
\begin{proof}
By the definition of a Chabauty limit of \( H_k \) in \( G_k \), as given in Section~\ref{subsec::short_Chabauty}, there exists a sequence \( \{g_\ell\}_{\ell \in \mathbb{N}} \) of elements in \( G_k \) such that the sequence of closed subgroups \( \{g_\ell H_k g_\ell^{-1}\}_{\ell \in \mathbb{N}} \) converges to \( L' \) with respect to the Chabauty topology on \( \mathcal{S}(G_k) \). 

Applying the \( K\mathcal{B}H \) decomposition from Theorem~\ref{thm::kbh_decom} to each element \( g_\ell \), we obtain \( g_\ell = k_\ell a_\ell h_\ell \), where \( k_\ell \in K \), \( a_\ell \in \mathcal{B}_k \), and \( h_\ell \in H_k \), for every \( \ell \in \mathbb{N} \). Since \( K \) is a compact subset of \( G_k \), we may, after extracting a subsequence, assume that the sequence \( \{k_\ell\}_{\ell \in \mathbb{N}} \) converges to some \( k \in K \). Moreover, since \( \mathcal{B}_k \) is a finite union of maximal \((\theta, k)\)-split tori, we may, after extracting another subsequence, assume that \( \{a_\ell\}_{\ell \in \mathbb{N}} \) lies in a fixed maximal $(\theta,k)$-split torus \( A \subset \mathcal{B}_k \). The sequence \( \{h_\ell\}_{\ell \in \mathbb{N}} \) does not affect the limit, as \( h_\ell H_k h_\ell^{-1} = H_k \). Finally, after possibly extracting a further subsequence, we may assume that \( \{a_\ell H_k a_\ell^{-1}\}_{\ell \in \mathbb{N}} \) converges to some \( L'' \in \mathcal{S}(G_k) \).

Thus, we have:
\[
\begin{split}
L' &= \limch_{\ell \to \infty} k_\ell a_\ell h_\ell H_k h_\ell^{-1} a_\ell^{-1} k_\ell^{-1} \\
&= \limch_{\ell \to \infty} k_\ell a_\ell H_k a_\ell^{-1} k_\ell^{-1} = k L'' k^{-1}.
\end{split}
\]
This means that there exist a \((\theta, k)\)-split torus \( A \subset \mathcal{B} \) and a sequence \( \{a_\ell\}_{\ell \in \mathbb{N}} \subset A \) (not necessarily consisting of hyperbolic elements) such that \( L' \) is \( G_k \)-conjugate to the Chabauty limit \( L'' = \limch_{\ell \to \infty} a_\ell H_k a_\ell^{-1} \).

Choose a maximal \( k \)-split torus \( T \subset G \) that contains \( A \). Such a torus exists because any \( k \)-split torus -- such as \( A \) -- can be extended to a maximal one. Then, by \cite[Proposition 4.5(iii)]{HW93} or \cite[Proposition 3.2(a)]{BeOh}, any maximal \( k \)-split torus of \( G \) that contains a maximal \((\theta, k)\)-split torus is \( \theta \)-stable. Therefore, \( T \) is \( \theta \)-stable.

Moreover, by the construction of the Bruhat–Tits building \( \Delta_{G_k} \) of \( G_k \), there is a one-to-one correspondence between maximal \( k \)-split tori of \( G \) and apartments of \( \Delta_{G_k} \). Thus, to \( T \) corresponds a unique apartment \( \mathcal{A}_T \subset \Delta_{G_k} \), and furthermore, \( t(\mathcal{A}_T) = \mathcal{A}_T \) setwise and \( \mathcal{A}_T \subset \Min(t) \) for every \( t \in T(k) \) (see \cite[Axiom 4.1.4]{KaPra}). Since \( T \) is \( \theta \)-stable, its corresponding apartment \( \mathcal{A}_T \) is also \( \theta \)-stable.

We deduce that \textbf{(Hyp\((\theta,k)\)-split)} is satisfied with the maximal \((\theta,k)\)-split torus \( A \) and the apartment \( \mathcal{A}_T \). Moreover, the apartment \( \mathcal{A} \) is \( \theta \)-stable.

\medskip
Let us now turn our attention to the condition \textbf{(Same \( \sigma_\pm \))}. So far, the sequence \( \{a_\ell\}_{\ell \in \mathbb{N}} \subset A(k) \) may consist of both elliptic and hyperbolic elements. If we can extract a subsequence \( \{a_{\ell_j}\}_{j \in \mathbb{N}} \) consisting of elliptic elements, then, by the property that \( \mathcal{A}_T \subset \Min(t) \) for every \( t \in T(k) \), the sequence \( \{a_{\ell_j}\}_{j \in \mathbb{N}} \) admits a convergent subsequence to some elliptic element \( a \in G_k \). Consequently, \( L'' = a H_k a^{-1} \), which means that \( L'' \) is \( G_k \)-conjugate to \( H_k \).

If instead we can only extract a subsequence \( \{a_{\ell_j}\}_{j \in \mathbb{N}} \) consisting of hyperbolic elements, we are faced with two possibilities: either \( \lim\limits_{\ell \to \infty} |a_\ell| < \infty \), or \( \lim\limits_{\ell \to \infty} |a_\ell| = \infty \). In the first case, by extracting a further subsequence if necessary and again using the property that \( \mathcal{A}_T \subset \Min(t) \) for every \( t \in T(k) \), we obtain a convergent subsequence to some hyperbolic element \( a \in G_k \). Consequently, we again have \( L'' = a H_k a^{-1} \), meaning that \( L'' \) is \( G_k \)-conjugate to \( H_k \).

\medskip
We are thus left with the non-trivial case in which $\lim\limits_{\ell \to \infty} |a_\ell| = \infty$.
In this scenario, by further extracting subsequences and utilizing the facts that $
\Delta_{G_k} \cup \partial \Delta_{G_k}$
comprises only finitely many types of simplices, and that $\mathcal{A}_T$
contains translation axes for each element \( a_\ell \), the conditions of 
\textbf{(Same \( \sigma_\pm \))} are satisfied for a subsequence $\{a_{\ell_j}\}_{j \in \mathbb{N}}$
of the hyperbolic elements $\{a_\ell\}_{\ell \in \mathbb{N}} \subset A(k)$.  For simplicity, we shall continue to index the subsequence $\{a_{\ell_j}\}_{j \in \mathbb{N}}$ by $\{a_\ell\}_{\ell \in \mathbb{N}}$.

\medskip
From this point onward, we work under the assumption of the condition \textbf{(Hyp \((\theta,k)\)-split)} for the maximal \((\theta,k)\)-split torus \( A \) 
and the \( \theta \)-stable apartment \( \mathcal{A}_T \),  as well as the condition \textbf{(Same \( \sigma_\pm \))} 
for the sequence of hyperbolic elements $\{a_\ell\}_{\ell \in \mathbb{N}} \subset A(k)$.

\medskip
We are now in a position to deduce the condition \textbf{(Hyp \( \sigma'_\pm = \sigma_\pm \))} for the sequence $\{a_\ell\}_{\ell \in \mathbb{N}} \subset A(k)$.  Since both \( T(k) \) and \( A(k) \) are abelian groups whose subgroups of hyperbolic elements  decompose as direct products of finitely many infinite cyclic groups, we can select, in each of \( T(k) \) and \( A(k) \), a finite set of hyperbolic elements such that every hyperbolic element in \( T(k) \), respectively in \( A(k) \), can be expressed as a product of powers of elements from this basis. To construct a generating set for \( A(k) \), one utilizes the existence of a set of simple roots for the root system \( \Phi(G, T) \) of \( G \) with respect to the torus \( T \). This set comprises both simple imaginary roots and simple non-imaginary roots (see \cite[Lemma 1.2]{DeConPro} or \cite[Chapter 36]{TaYu}). The latter subset -- the simple non-imaginary roots -- is used to generate \( A(k) \). 

We apply this decomposition to each element of the sequence $\{a_\ell\}_{\ell \in \mathbb{N}}$, which, as previously assume, satisfies the condition \textbf{(Same \( \sigma_\pm \))}. From there, we recall that $\lim\limits_{\ell \to \infty} \xi_{a_\ell \pm} = \xi_\pm$,
where \( \xi_{a_\ell +} \) and \( \xi_{a_\ell -} \) denote the attracting and repelling endpoints of \( a_\ell \), respectively.  Each point \( \xi_\pm \) therefore lies within, or possibly on the boundary of, the corresponding ideal simplex \( \sigma_\pm \). If \( \xi_\pm \) lie in the interior of \( \sigma_\pm \), the  condition \textbf{(Hyp \( \sigma'_\pm = \sigma_\pm \))} follows immediately.  If they lie on the boundary, then some of the powers in the decomposition of \( a_\ell \) must remain bounded. 
Since \( A(k) \) is abelian, we may isolate the bounded components on the left-hand side of the decomposition of \( a_\ell \). 
Furthermore, by passing to a subsequence if necessary, we may assume that these components remain constant for all \( a_\ell \). 
The remaining terms \( a'_\ell \) in the decomposition thus form a sequence of hyperbolic elements  \( \{a'_\ell\}_{\ell \in \mathbb{N}} \subset A(k) \) 
that satisfies the condition \textbf{(Hyp \( \sigma'_\pm = \sigma_\pm \))}. 
Consequently, the limit $L''=\limch_{\ell \to \infty} a_\ell H_k a_\ell^{-1}$ is \( G_k \)-conjugate to $L:=\limch_{\ell \to \infty} a'_\ell H_k (a'_\ell)^{-1}$. 

Since we have identified a \((\theta, k)\)-split torus \( A \) contained in \( \mathcal{B} \), 
along with a sequence \( \{a'_\ell\}_{\ell \in \mathbb{N}} \subset A \) of hyperbolic elements 
satisfying the conditions \textbf{(Same \( \sigma_\pm \))} and \textbf{(Hyp \( \sigma'_\pm = \sigma_\pm \))}, 
and given that \( L' \) is \( G_k \)-conjugate to $L:=\lim^{\mathrm{ch}}_{\ell \to \infty} a'_\ell H_k (a'_\ell)^{-1}$, the first part of the theorem follows.

Part 1) of the remainder of the theorem is deduced by first applying Theorem \ref{thm::chabauty_in_parabolic} 
to the sequence \( \{a'_\ell\}_{\ell \in \mathbb{N}} \subset A \), where \( A \) is a \((\theta,k)\)-split torus 
and \( \mathcal{A}_T \) is a \( \theta \)-stable apartment. 
By the condition \textbf{(Hyp \( \sigma'_\pm = \sigma_\pm \))}, we obtain \( L \leq P_{\sigma_+} \).
Moreover, by \cite[Proposition 9.2, Section 4 and Proposition 4.7]{HW93}, there exists a minimal \( \theta \)-split parabolic \( k \)-subgroup \( P' \) of \( G \) that contains a \( \theta \)-stable maximal \( k \)-split torus of \( G \), which in turn contains the maximal \( (\theta, k) \)-split torus \( A \). Furthermore, the centralizer \( Z_G(A) \) is equal to \( P' \cap \theta(P') \).

Since \( P' \subseteq P_{\sigma_+} \), it follows from \cite[Corollaries 4.16 and 4.19]{BoTi} that the intersection \( A' \) of \( A \) with the radical \( \mathcal{R}(P_{\sigma_+}) \) of \( P_{\sigma_+} \) is a maximal \( k \)-split torus of \( \mathcal{R}(P_{\sigma_+}) \). Moreover, by \cite[Theorem 4.15]{BoTi}, we have \( Z_G(A') = P_{\sigma_+} \cap P_{\sigma_-} \).  In addition, the sequence \( \{a'_\ell\}_{\ell \in \mathbb{N}} \) consists of hyperbolic elements in \( A'(k) \).

Then, by \cite[Theorem 4.15]{BoTi}, we obtain $M_{\sigma_-,\sigma_+}(k) = Z_{G_k}(A')$. Therefore, we may apply Corollary \ref{cor::find_elements_H} and Propositions \ref{prop::levi_factors} and \ref{prop::find_unipotent}, 
using the \((\theta, k)\)-split torus \( A' \) and the sequence \( \{a'_\ell\}_{\ell \in \mathbb{N}} \), 
to deduce Part 1) of the theorem.

Part 2) follows directly from the identity $M_{\sigma_-,\sigma_+}(k) = Z_{G_k}(A')$ and Corollary \ref{cor::find_elements_H}. 

Part 3) corresponds to Proposition \ref{prop::exist_l_sigma}, and Part 4) is established by Proposition \ref{prop::find_unipotent}. 

This completes the proof of the theorem.
\end{proof} 

\section{Appendix A: The $K\mathcal{B}H$ decomposition and symmetric varieties}
\label{appen::A}

Let \( G \) be a connected reductive linear algebraic group defined over a non-Archimedean local field \( k \). Due to the applicability of the Inverse Function Theorem in the context of analytic manifolds over non-Archimedean fields -- which holds only in characteristic zero (see Appendix~\ref{appen:B}, and more precisely \cite{PlaRa}, middle of page 110) -- we assume moreover that the characteristic of \( k \) is zero.

Let \( \theta \) be an involutive \( k \)-automorphism of \( G \), that is, \( \theta^2 = \mathrm{id} \) and \( \theta \) is defined over \( k \).
Let \( G^{\theta} \leq G \) be the fixed-point subgroup of \( \theta \):
\[
G^{\theta} := \{ h \in G \mid \theta(h) = h \}
\]
and $H = (G^{\theta})^{o}$ be the connected component of $G^{\theta}$. Then $H$ is defined over $k$, since $\theta$ is defined over $k$ (see \cite[Proposition 1.6]{HW93}). Moreover, by \cite[Theorem 2.1]{PraYu}, $H$ is reductive. By \( G_k \), \( H_k \) we denote the groups of \( k \)-points.

We present the \( K\mathcal{B}_kH_k \) decomposition of \( G_k \) as a finite union of double cosets involving compact subsets, maximal \( (\theta,k) \)-split tori, and the fixed point group \( H_k \) of the \( k \)-involution \( \theta \) on \( G \). This decomposition plays a central role in understanding the structure of Chabauty limits and the geometry of symmetric varieties.

We also introduce the following subsets of \( G \):
\[
Q' := \{ g \in G \mid \theta(g) = g^{-1} \}, \quad Q := \{ g\theta(g^{-1}) \mid g \in G \},
\]
noting that \( Q \subset Q' \).

There is a natural \( (G,\theta) \)-twisted action on \( Q' \) given by:
\[
(g,x) \in G \times Q' \mapsto g x \theta(g^{-1}) \in Q'.
\]
Under this action, \( Q' \) is the \( (G,\theta) \)-orbit of the identity element \( e \in G \), and the stabilizer  $\Stab_{(G,\theta)}(e)$ of \( e \) is precisely \( H \). This identification allows us to view \( Q \) as the symmetric variety \( G/H \) associated with the pair \( (G,H) \) and the involution \( \theta \).

To formalize this, we define the map:
\begin{equation}
\label{equ::tau_map}
\tau : G \to G, \quad x \mapsto \tau(x) := x \theta(x^{-1}),
\end{equation}
whose image \( \tau(G) \) coincides with the variety \( Q \).

According to \cite[Section 9]{Richar}, the variety \( Q' \) contains only finitely many \( (G,\theta) \)-twisted orbits, each of which is a closed \( k \)-subvariety of \( G \). In particular, this finiteness also applies to the symmetric variety \( Q \).

\subsection{$G_k$ as a disjoint union of double cosets $H_kvP_k$}

Let us fix a minimal parabolic \( k \)-subgroup \( P \subset G \). It is known (see \cite[Lemma 2.4]{HW93} and \cite[Proposition 3.3(b)]{BeOh}) that \( P \) contains a \( \theta \)-stable maximal \( k \)-split torus \( A \), which we fix for the remainder of the discussion. For completeness, we recall the proof of this fact below.

\begin{lemma}[Lemma 2.4 of \cite{HW93}]
\label{lem::if_correct}
Every minimal parabolic \( k \)-subgroup \( P \) of \( G \) contains a \( \theta \)-stable maximal \( k \)-split torus of $P$, unique up to conjugation by an element of \( (H \cap U)_k \), where \( U \) is the unipotent radical of \( P \). Consequently, \( A \) is a \( \theta \)-stable maximal \( k \)-split torus of \( G \).
\end{lemma}

\begin{proof}
Consider the subgroup \( L = P \cap \theta(P) \). By \cite[Lemma 2.2]{HW93}, \( L \) admits a \( \theta \)-stable Levi \( k \)-subgroup \( M \). Since \( M \) is the intersection of two opposite parabolic subgroups of \( G \), we may apply \cite[Corollary 4.18]{BoTi} to conclude that \( M \) contains the centralizer \( Z_G(T) \) of a maximal \( k \)-split torus \( T \) of \( G \).

Because \( P \) is minimal, \cite[Corollary 4.16(iii), (iv)]{BoTi} implies that the centralizer \( Z_G(T) \) is a Levi subgroup of \( P \). Therefore, the \( \theta \)-stable Levi subgroup \( M \) must coincide with \( Z_G(T) \).

Since \( \theta \) is defined over \( k \), it preserves the property of being a maximal \( k \)-split torus. Thus, \( \theta(T) \) is also a maximal \( k \)-split torus of \( G \). The \( \theta \)-stability of \( M \) implies that
\[
M = Z_G(\theta(T)) = Z_G(T),
\]
which in turn shows that \( \theta(T) = T \). Hence, \( T \) is a \( \theta \)-stable maximal \( k \)-split torus of \( G \).

The remainder of the proof follows as in \cite{HW93}.
\end{proof}

Furthermore, by \cite[Lemma 4.8, Propositions 6.8 and 9.2]{HW93}, there exists a minimal \( \theta \)-split parabolic \( k \)-subgroup of \( G \). The definition of a \( \theta \)-split parabolic subgroup was recalled at the beginning of Section~\ref{subsec::L_trans_sigma_opp}. According to \cite[Proposition 4.7]{HW93}, any minimal \( \theta \)-split parabolic \( k \)-subgroup of \( G \) contains a \( \theta \)-stable maximal \( k \)-split torus \( T \subset G \), such that the \( T_- \) component of \( T \) is a maximal \( (\theta,k) \)-split torus of \( G \). Since any two maximal \( k \)-split tori of \( G \) are conjugate under the action of \( G_k \), the same holds for any two maximal \( (\theta,k) \)-split tori of \( G \).

Although not yet formulated in a fully rigorous mathematical framework, we heuristically interpret the \( k \)-points \( Q'_k \) of the variety \( Q' \) as an analogue of an affine building where ``apartments'' correspond to the normalizers \( N_{G_k}(T) \) for specific maximal \( k \)-split tori \( T \subset G \).

In particular, consider the fixed \( \theta \)-stable maximal \( k \)-split torus \( A \subset P \) of $G$. Every element of the variety \( Q'_k \) can be transported into the ``apartment'' associated with \( N_{G_k}(A) \) via the \( (G,\theta) \)-twisted action of the unipotent subgroup \( U_k \subset P_k \). Since \( A \) is \( \theta \)-stable, its normalizer \( N_G(A) \) is still \( \theta \)-stable. This mechanism is precisely described in \cite[Proposition 6.6]{HW93}.

\begin{proposition}[Proposition 6.6, \cite{HW93}]
\label{prop:unipotent_moves_Q}
Let \( g \in G_k \) be such that \( \theta(g) = g^{-1} \), i.e., \( g \in Q'_k \). Then there exists \( u \in U_k \) such that $
ug\, \theta(u)^{-1} \in N_{G_k}(A)$.
\end{proposition}
\begin{proof}
We recall the proof from \cite[Proposition 6.6, Section 6]{HW93}, reformulated for clarity.

\textbf{Step 1: The \(\theta\)-twisted Bruhat decomposition.}  
We use the \(\theta\)-twisted version of the Bruhat decomposition that is given by $G_k = P_k N_{G_k}(A) \theta(P_k) = U_k N_{G_k}(A) \theta(U_k)$,
which holds because \( A \) is \( \theta \)-stable, and hence its normalizer \( N_G(A) \) is also \( \theta \)-stable. Moreover, the \( N_{G_k}(A) \)-component in this decomposition is unique; see \cite[Section 6.5]{HW93}.  Using this decomposition, we can write
\[
g = u_1 n u_2,
\]
with \( u_1 \in U_k \), \( u_2 \in \theta(U_k) \), and \( n \in N_{G_k}(A) \). Applying \( \theta \) to both sides and using the condition \( \theta(g) = g^{-1} \), we obtain:
\[
\theta(u_1)\theta(n)\theta(u_2) = \theta(g) = g^{-1} = u_2^{-1} n^{-1} u_1^{-1}.
\]
By the uniqueness of the \( N_{G_k}(A) \)-component in the twisted Bruhat decomposition, it follows that
\[
\theta(n) = n^{-1}.
\]
\textbf{Step 2: The involutions \(\theta\) and \(\theta_n\) on  \( N_{G_k}(A) \).}  
Recall that \( \theta \) is a \( k \)-involution and that \( N_G(A) \) is \( \theta \)-stable, so \( \theta \) restricts to a \( k \)-involution on \( N_G(A) \).

Given \( n \in N_{G_k}(A) \) with \( \theta(n) = n^{-1} \), we define a new involution \( \theta_n \) on \( G \) by
\[
\theta_n(x) := n \theta(x) n^{-1}, \quad \text{for } x \in G.
\]
By \cite[Lemma 6.1]{HW93}, \( \theta_n \) is a \( k \)-involution that also preserves \( N_G(A) \), i.e., \( N_G(A) \) is \( \theta_n \)-stable.

Both $k$-involutions $\theta$ and $\theta_n$ can be regarded as inducing involutional automorphisms on the Bruhat--Tits building $\Delta_{G_k}$ of $G_k$ (see Section \ref{section::auto_invol_building}). 

%And since they stabilize $N_G(A)$, they also induce involutional automorphisms on the apartment of $\Delta_{G_k}$ corresponding to the $\theta$-stable maximal $k$-split torus $A$ of $G$.

\textbf{Step 3: The unipotent radical \( U_k \) as a product of intersections with opposite unipotent radicals.}  
Let \( P' \) be another minimal parabolic \( k \)-subgroup of \( G \), conjugate to \( P \) via an element of \( N_{G_k}(A) \). Let \( U' \) denote the unipotent radical of \( P' \), and let \( (P')^{-} \) be the parabolic subgroup opposite to \( P' \) with respect to the ``apartment'' \( N_{G_k}(A) \), with corresponding unipotent radical \( (U')^{-} \), which is opposite to \( U' \).

It is shown in \cite[proof of Lemma 4.1]{Springer} that the product map
\[
(U_k \cap U') \times (U_k \cap (U')^{-}) \longrightarrow U_k
\]
is a bijection.

This result is applied in the proof of \cite[Proposition 6.6]{HW93} to the unipotent radical \( U_k \) and the opposite unipotent radicals \( \theta_n(U) \) and \( \theta_n(U^{-}) \). Consequently, any element \( u \in U_k \) can be written as
\[
u = v_1 v_2,
\]
where \( v_1 \in U_k \cap \theta_n(U^{-}) \) and \( v_2 \in U_k \cap \theta_n(U) \).

From this, we deduce two key observations:
\begin{itemize}
    \item First, consider the element \( u_1 \) from the decomposition \( g = u_1 n u_2 \). We can write \( u_1 = w_1 w_2 \), with \( w_1 \in U_k \cap \theta_n(U) \) and \( w_2 \in U_k \cap \theta_n(U^{-}) \). Applying the \( (G, \theta) \)-twisted action of \( w_1^{-1} \) to \( g \), we obtain:
    \[
    w_1^{-1} g \theta(w_1) = w_1^{-1} w_1 w_2 n u_2 \theta(w_1) = w_2 n u_2 \theta(w_1),
    \]
    where \( w_2 \in U_k \cap \theta_n(U^{-}) \) and \( u_2 \theta(w_1) \in \theta(U_k) \). Thus, we may assume that \( u_1 \in U_k \cap \theta_n(U^{-}) = U_k \cap n \theta(U^{-}) n^{-1} \). This refinement will be crucial in the next step.
    
    Note that the element \( w_1^{-1} g \theta(w_1) \) remains in \( Q'_k \subset G_k \).
\item
We now consider the second key observation. Recall that \( \theta(u_2)^{-1} \) can be written as
\[
\theta(u_2)^{-1} = v_1 v_2 \in U_k,
\]
where \( v_1 \in U_k \cap \theta_n(U^{-}) \) and \( v_2 \in U_k \cap \theta_n(U) \). Substituting this into the identity \( \theta(g) = g^{-1} \), and using the decomposition \( g = u_1 n u_2 \), we obtain: 
$$
\theta(u_1) n^{-1} v_2^{-1} v_1^{-1} = \theta(v_1) \theta(v_2) n^{-1} u_1^{-1}.$$

To exploit the opposition between \( \theta(U) \) and \( \theta(U^{-}) \), we multiply both sides on the right by \( n \), yielding: $
(\theta(u_1) n^{-1} v_2^{-1} n)(n^{-1} v_1^{-1} n) = \theta(v_1) \theta(v_2) n^{-1} u_1^{-1} n$.

Observe that:
\begin{itemize}
    \item \( \theta(u_1) n^{-1} v_2^{-1} n \in \theta(U) \),
    \item \( \theta(v_1) \theta(v_2) \in \theta(U) \),
    \item \( n^{-1} v_1^{-1} n, n^{-1} u_1^{-1} n \in \theta(U^{-}) \).
\end{itemize}

This implies that both sides of the equation lie in the intersection \( \theta(U) \cap \theta(U^{-}) \). Given our earlier assumption that \( u_1 \in U_k \cap \theta_n(U^{-}) = U_k \cap n \theta(U^{-}) n^{-1} \), the equality above must reduce to the identity. Therefore, we conclude:
\[
n^{-1} v_1^{-1} n = n^{-1} u_1^{-1} n \quad \Rightarrow \quad v_1 = u_1, \quad  \text{ and } \theta_n(v_2) = v_2^{-1}.
\]
This yields the expression: $g = v_1 n \theta(v_2^{-1} v_1^{-1})$.

And applying the \( (G, \theta) \)-twisted action of \( v_1^{-1} \) to \( g \), we obtain:
$$v_1^{-1} g \theta(v_1) = n \theta(v_2^{-1}) = v_2 n.$$
\end{itemize}

\textbf{Step 4: Expressing \( v_2 \) as \( y \theta_n(y^{-1}) \) and identifying the desired element \( u \in U_k \).}  
We now use the fact that the subgroup \( U \cap \theta_n(U) \) is defined over \( k \) and is \( \theta_n \)-stable. Therefore, we can apply the decomposition result from \cite[Lemma 0.6(ii)]{HW93} to the element \( v_2 \in U_k \cap \theta_n(U) \), which satisfies the relation \( \theta_n(v_2) = v_2^{-1} \). This implies that there exists an element \( y \in U_k \cap \theta_n(U) \) such that
\[
v_2 = y \theta_n(y^{-1}).
\]
We now apply the \( \theta \)-twisted action of \( y^{-1} \) to the expression \( v_1^{-1} g \theta(v_1) = v_2 n \), yielding:
\[
y^{-1} v_1^{-1} g \theta(v_1) \theta(y) = y^{-1} v_2 n \theta(y) = y^{-1} v_2 \theta_n(y) n = n.
\]
From this, we deduce that the element $u := (v_1 y)^{-1} \in U_k$ satisfies the desired property stated in the proposition.
\end{proof}

As noted in \cite[Section 6]{HW93}, there may not exist \( \theta \)-stable minimal parabolic \( k \)-subgroups \( P \subset G \). To address this, one considers the pair \( (P, \theta(P)) \), which yields the \( \theta \)-stable subgroup \( P \cap \theta(P) \).

Since we are interested in the variety $Q := \{ g \theta(g^{-1}) \mid g \in G$, and given that elements \( x \in Q \) can be moved via the \( (U_k, \theta) \)-twisted action to the \( \theta \)-stable  \( N_{G_k}(A) \) (by Proposition~\ref{prop:unipotent_moves_Q}), it is natural to investigate which elements of \( N_{G_k}(A) \) lie in the variety \( Q_k \). To this end, we consider the map \( \tau \) from (\ref{equ::tau_map}) and define the preimage of \( Q \cap N_{G_k}(A) \) under \( \tau \) as
\[
\tau^{-1}(N_{G_k}(A)) := \{ g \in G_k \mid \tau(g) = g \theta(g^{-1}) \in N_{G_k}(A) \}.
\]
This set \( \tau^{-1}(N_{G_k}(A)) \subset G_k \) admits a natural action by the group \( Z_{G_k}(A) \times H_k \), defined by
\begin{equation}
\label{equ::action_Z_H}
(y, x, z) \in Z_{G_k}(A) \times H_k \times \tau^{-1}(N_{G_k}(A)) \quad \mapsto \quad (y, x) \cdot z := y z x.
\end{equation}
To verify that \( y z x \in \tau^{-1}(N_{G_k}(A)) \), observe that
\[
\tau(y z x) = y z x \theta(x^{-1} z^{-1} y^{-1}) = y z \theta(z^{-1} y^{-1}) \in N_{G_k}(A),
\]
by the \( \theta \)-stability of \( N_G(A) \) and because $Z_{G}(A) \leq  N_G(A)$. Let \( V \) denote a set of representatives for the \( Z_{G_k}(A) \times H_k \)-orbits in \( \tau^{-1}(N_{G_k}(A)) \).

Let us recall again that $P$ is a minimal parabolic $k$-subgroup of $G$ and that $A \subset P$ is a $\theta$-stable maximal $k$-split torus of $G$. Denote by $U$ the unipotent radical of $P$. 

\begin{proposition}[Proposition 6.8 in \cite{HW93}]
\label{prop::HvP_orbits}
Under the above notation, the group \( G_k \) is the disjoint union of the double cosets \( P_k v H_k \), with \( v \in V \).
\end{proposition}

\begin{proof}
Let \( g \in G_k \). By Proposition~\ref{prop:unipotent_moves_Q}, there exists \( u \in U_k \) such that $
u g \theta(g^{-1}) \theta(u^{-1}) \in N_{G_k}(A)$.

This implies that \( u g \in \tau^{-1}(N_{G_k}(A)) \), and hence \( g \in U_k \tau^{-1}(N_{G_k}(A) \). Since \( P_k = U_k \rtimes Z_{G_k}(A) \), and \( Z_{G_k}(A) \times H_k \) acts on \( \tau^{-1}(N_{G_k}(A)) \), we obtain the (not necessarily disjoint) union
\[
G_k = \bigcup_{v \in V} P_k v H_k.
\]
We now show that this union is in fact disjoint. Suppose there exist \( g_1, g_2 \in \tau^{-1}(N_{G_k}(A)) \) such that \( g_2 \in P_k g_1 H_k \). To prove disjointness, we must show that \( g_2 \in Z_{G_k}(A) g_1 H_k \).

Write \( g_2 = y g_1 x \), with \( x \in H_k \) and \( y \in P_k \). Since \( P_k = U_k \rtimes Z_{G_k}(A) \), we can write \( y = u z \), with \( u \in U_k \), \( z \in Z_{G_k}(A) \). Then
\[
g_2 \theta(g_2^{-1}) = u z g_1 x \theta((u z g_1 x)^{-1}) = u (z g_1) \theta((z g_1)^{-1}) \theta(u^{-1}) \in U_k N_{G_k}(A) \theta(U_k).
\]
This expression lies in the twisted Bruhat decomposition \( G_k = U_k N_{G_k}(A) \theta(U_k) \), and by uniqueness of the \( N_{G_k}(A) \)-component, we conclude that \( u \in U_k \) must be trivial. Hence \( y = z \in Z_{G_k}(A) \), and so \( g_2 = z g_1 x \in Z_{G_k}(A) g_1 H_k \), as claimed.

If in addition we consider that $H$ is an open $k$-subgroup of the the fixed point group $G^{\theta}$, we have to continue doing the following.  Since both \( g_1, g_2 \in \tau^{-1}(N_{G_k}(A)) \), the \( \theta \)-twisted Bruhat decomposition from \cite[Section 6.5(ii)]{HW93}, which guarantees uniqueness of the \( N_{G_k}(A) \)-component, implies that $
g_2 \theta(g_2^{-1}) = (z g_1) \theta((z g_1)^{-1})$,
and hence \( g_2 = z g_1 x_1 \), where \( x_1 \in (G_k)^{\theta} \). We are nearly done, but it is not immediately clear whether \( (G_k)^{\theta}\) is in the subgroup \(H_k \), the $k$-points of the open $H$.

However, from our earlier expression \( g_2 = u z g_1 x \), we also have
\[
u z g_1 x = z g_1 x_1,
\]
which implies
\[
x_1 x^{-1} = g_1^{-1} z^{-1} u z g_1 = (z g_1)^{-1} u z g_1.
\]
This element lies in the conjugate subgroup \( (z g_1)^{-1} U z g_1 \), and is also \( \theta \)-fixed, i.e.,
\[
\theta(x_1 x^{-1}) = x_1 x^{-1}.
\]
We now apply \cite[Lemma 10.1]{HW93}, which states that the set of \( \theta \)-fixed points in a unipotent subgroup -- such as \( (z g_1)^{-1} U z g_1 \), a conjugate of \( U \) -- is connected and defined over \( k \). Therefore, \( x_1 x^{-1} \in (H^0)_k \), and it follows that \( x_1 \in H_k \).  This shows that \( g_2 \in Z_G(A)_k g_1 H_k \), completing the proof of the claim.
\end{proof}

\subsection{Finite number of $H_k$-conjugacy classes of maximal $(\theta,k)$-split tori}

We again fix a minimal parabolic \( k \)-subgroup \( P \) of \( G \). By Lemma \ref{lem::if_correct}, \( P \) admits a \( \theta \)-stable maximal \( k \)-split torus \( A _1 \subset P\) of $G$  that is unique up to conjugation by an element of $(H \cap U)_k$, where $U$ is the unipotent radical of $P$. The torus $A_1$ is fixed for what follows.

In this section, we recall from \cite[Section 6.9]{HW93} how the \( H_k \)-conjugacy classes of \( \theta \)-stable maximal \( k \)-split tori in \( G \) can be described in terms of the double cosets \( H_k v P_k \).

\medskip
Let \( \mathcal{A} \) denote the set of all \( \theta \)-stable maximal \( k \)-split tori of \( G \), on which \( H_k \) acts by conjugation:
\[
(h, A) \in H_k \times \mathcal{A} \mapsto h \cdot A := h A h^{-1}.
\]

For \( A \in \mathcal{A} \), note that $A$ is, in particular, a maximal \( k \)-split torus of \( G \). It follows from \cite[Theorem 4.15.a and Corollary 4.16]{BoTi} that \( A \) is contained in a minimal parabolic \( k \)-subgroup of \( G \). Moreover, any two minimal parabolic \( k \)-subgroups of \( G \) are \( G_k \)-conjugate. Therefore, there exists \( g \in G_k \) such that $A \subseteq g P g^{-1}$, and  $A$ is, in particular, a \( \theta \)-stable maximal \( k \)-split torus  in $gPg^{-1}$.

We now define the map
\[
\pi : H_k \backslash G_k / P_k \longrightarrow H_k \backslash \mathcal{A}
\]
by sending a double coset \( H_k g P_k \), with \( g \in G_k \), to the \( H_k \)-conjugacy class \([A]\) of a \( \theta \)-stable maximal \( k \)-split torus \( A \subset g P g^{-1} \) in \( G \). By Lemma~\ref{lem::if_correct}, this is well-defined: such a torus exists, and any two \( \theta \)-stable maximal \( k \)-split tori contained in \( gPg^{-1} \) are conjugate by an element of \( (H \cap U)_k \).

\medskip
We observe that the map \( \pi \) is surjective. Indeed, since all maximal \( k \)-split tori of \( G \) are \( G_k \)-conjugate, the same holds for the \( \theta \)-stable ones. Moreover, because our fixed minimal parabolic \( k \)-subgroup \( P \) contains the \( \theta \)-stable maximal \( k \)-split torus \( A_1 \), any other element \( A \in \mathcal{A} \) is contained in some conjugate \( g P g^{-1} \), where \( A = g A_1 g^{-1} \) for some \( g \in G_k \).

Thus, given a class \([A] \in H_k \backslash \mathcal{A}\) corresponding to a parabolic subgroup \( g P g^{-1} \), for some \( g \in G_k \) such that \( A = g A_1 g^{-1} \), we aim to compute the fiber \( \pi^{-1}([A]) \) -- that is, the set of all \( G_k \)-conjugates of \( P \) that contain a torus from \([A]\).

To this end, recall that for any \( A \in \mathcal{A} \), the minimal parabolic \( k \)-subgroups of \( G \) containing \( A \) are \( N_{G_k}(A) \)-conjugate. This follows from viewing \( A \) as an apartment in the Bruhat–Tits building of \( G_k \). The ideal boundary \( \partial A \) has finitely many ideal chambers, each corresponding to a minimal parabolic \( k \)-subgroup of \( G \) containing \( A \). These ideal chambers are in bijection with the elements of the Weyl group
\[
W_{G_k}(A) := N_{G_k}(A) / Z_{G_k}(A).
\]
Therefore, the fiber is given by
\[
\pi^{-1}([A]) = H_k \backslash H_k N_{G_k}(A) g P_k / P_k.
\]

Next we want to further understand $\pi^{-1}([A]) =  H_k \backslash H_k N_{G_k}(A)g P_k /P_k$ as a possible disjoint union of double cosets $H_k v P_k$, and eventually relate it to the Weyl group $W_{G_k}(A) :=N_{G_k}(A)/Z_{G_k}(A)$.
Indeed, first recall that $A= gA_1g^{-1}$, for some $g \in G_k$. Then let $n,n_1 \in N_{G_k}(A)$ so that $H_k n_1g P_k = H_kngP_k$. We want to see how $n_1$ is related to $n$. By the assumption we can write
$$n_1 g = hngzu$$
with $h \in H_k$, $z$ in the Levi factor $(Z_G(A_1))_k = Z_{G_k}(A_1)$ and $u \in U_k$ which is the $k$-points of the unipotent radical $U$ of $P$.
Because $A= gA_1g^{-1}$, and $n,n_1 \in  N_{G_k}(A)$, we notice that 
$$ngA_1 g^{-1}n^{-1}= nAn^{-1}=A \; \; \text{ and }\; \; ngzu A_1 (ngzu)^{-1} = h^{-1} n_1 g A_1 (n_1g)^{-1}h = h^{-1}Ah$$
are both $\theta$-stable maximal $k$-split tori of $ng P (ng)^{-1}$, as $zu A_1(zu)^{-1}$ is a maximal $k$-split torus of $P$ but maybe not $\theta$-stable. By \cite[Lemma 2.4]{HW93} $\theta$-stable maximal $k$-split tori are unique in their corresponding parabolic, up to some conjugacy, and so we find $v' \in (H \cap R_u(ngP(ng)^{-1}))_k \leq H_k \leq (G^{\theta})_k$ such that 
$$v'(ngA_1 g^{-1}n^{-1})(v')^{-1} = ngzu A_1 (ngzu)^{-1}.$$

Since \( z \) lies in the Levi factor \( (Z_G(A_1))_k = Z_{G_k}(A_1) \), we have \( z A_1 z^{-1} = A_1 \). With a bit of manipulation, this yields:
\[
(ngz)^{-1} v' ngz \cdot A_1 \cdot \left((ngz)^{-1} v' ngz\right)^{-1} = u A_1 u^{-1},
\]
where \( u \in U \), and \( (ngz)^{-1} v' ngz \in (ngz)^{-1} R_u(ngP(ng)^{-1}) (ngz) = U \). This implies that \( u^{-1} (ngz)^{-1} v' ngz \in U \), but also that it lies in \( N_{G_k}(A_1) \).

However, since \( N_G(A_1) \cap U = \{e\} \), it follows that \( u = (ngz)^{-1} v' ngz \), and hence 
$$ngz u = v' ngz$$
for some \( v' \in (H \cap R_u(ngP(ng)^{-1}))_k \subseteq H_k \).

With that in hand, the we can rewrite $n_1g = hngzu = hv'ngz$ and so  $n_1 = h' ngzg^{-1}$ for some $h'\in H_k$. Since $gzg^{-1} \in Z_{G_k}(A)$, then $n_1 \in H_k n Z_{G_k}(A)$. This implies that the double coset $H_k n_1g P_k = H_kngP_k$ produces a unique element in $W_{H_k}(A) \backslash W_{G_k}(A)$, where $W_{H_k}(A):=N_{H_k}(A)/Z_{H_k}(A)$  and  $W_{G_k}(A):=N_{G_k}(A)/Z_{G_k}(A)$. 

Therefore we have $$ \pi^{-1}([A]) \cong W_{H_k}(A) \backslash W_{G_k}(A) $$
and the following characterization of the double coset decomposition of $G_k$.

\begin{proposition}(Proposition 6.10 in \cite{HW93})
\label{prop::diff_coset_decom_G_k}
Let $\{A_i \; \vert \; i \in I\} $ be the representatives of the $H_k$-conjugacy classes of $\theta$-stable maximal $k$-split tori of $G$. Then 
$$ H_k \backslash G_k /P_k \cong \bigcup_{i \in I}W_{H_k}(A_i) \backslash W_{G_k}(A_i).$$
\end{proposition}

% Take $v := (ng)^{-1} v' (ng)$ which is an element of $R_u(P)_k$ and we have $$ngzu A_1 (ngzu)^{-1} = v'(ngA_1 g^{-1}n^{-1})(v')^{-1} = ngv A_1 (ngv)^{-1} \Rightarrow zu A_1 (zu)^{-1} = vA_1v^{-1}.$$Rewriting the above equality using $z \in (Z_G(A_1))_k$ we get 
%$$ vA_1v^{-1} = zu A_1 (zu)^{-1} =  zuz^{-1}z A_1z^{-1} (zuz^{-1})^{-1} =  zuz^{-1} A_1 (zuz^{-1})^{-1}. $$
%Since $v ,zuz^{-1} \in R_u(P)$ and $N_G(A_1) \cap  R_u(P) = \{e\}$, we must have that $v =zuz^{-1}$ which implies that
%$$ngzu = ng vz = ngv(ng)^{-1} ng z.$$
%This gives that $ngv(ng)^{-1} \in R_u(ngzP (ngz)^{-1})$ as well as in $H_k$ by the above properties of $v' = ngv (ng)^{-1}$. We can then  apply \cite[Lemma 10.1]{HW93}, which basically says that the set of $\theta$-fixed-points  of a unipotent subgroup (in our case $ngzU (ngz)^{-1}$ is a conjugate of the unipotent $U$) is connected and defined over $k$, and obtain $v' = ngv (ng)^{-1} \in (H^{0})_k$. 

Recall the notation:
\[
Q' := \{ g \in G \mid \theta(g) = g^{-1} \}, \quad \text{and} \quad Q := \{ g \theta(g^{-1}) \mid g \in G \}.
\]
Note that \( Q \subset Q' \).

We have fixed a minimal parabolic \( k \)-subgroup \( P \subset G \), and we know that \( P \) admits a \( \theta \)-stable maximal \( k \)-split torus \( A \subset P \) of \( G \), that we fix for what follows.  Let \( U \) denote the unipotent radical of \( P \).

Our interest lies in the variety \( Q \), but by studying the larger variety \( Q' \), we gain more flexibility to understand the structure of \( Q \). So far, we have used Proposition~\ref{prop:unipotent_moves_Q} and the \( (U_k, \theta) \)-action to reduce the study of the \( k \)-points \( Q_k \) of \( Q \) to the \( \theta \)-stable ``apartment" \( N_{G_k}(A) \).

Drawing a parallel with the finite Weyl group and the Bruhat decomposition, the next step is to investigate how many \( (P_k, \theta) \)-twisted orbits -- equivalently, how many Weyl group elements -- are present in \( Q'_k \), and hence in \( Q_k \). This is the focus of the next proposition, where the nature of the field \( k \) plays a crucial role.

\begin{proposition}[Proposition 6.15 in \cite{HW93}]
\label{prop::6.15}
Let \( k \) be a local field of characteristic zero. Then the intersection \( G_k \cap Q' \) contains only finitely many \( (P_k, \theta) \)-twisted orbits.
\end{proposition}

\begin{proof}
The proof proceeds in several steps.

\textbf{Step 1.}  
By definition, the finite Weyl group \( W_{G_k}(A) := N_{G_k}(A) / Z_{G_k}(A) \) associated to \( G_k \) has finitely many elements. Let \( n_1, \dots, n_\ell \in N_{G_k}(A) \) be representatives of these cosets, so that:
\[
N_{G_k}(A) = \bigsqcup_{j=1}^{\ell} Z_{G_k}(A) n_j.
\]
By Proposition~\ref{prop:unipotent_moves_Q}, every \( (U_k, \theta) \)-twisted orbit in \( Q'_k \) intersects the set:
\[
N_{G_k}(A) \cap Q' = \left( \bigsqcup_{j=1}^{\ell} Z_{G_k}(A) n_j \right) \cap Q'.
\]
Since \( P_k = U_k \rtimes Z_{G_k}(A) \), and we are interested in the finiteness of \( (P_k, \theta) \)-twisted orbits, it suffices to show that for each \( j \in \{1, \dots, \ell\} \), the set \( Z_{G_k}(A) n_j \cap Q' \) contains only finitely many \( (Z_{G_k}(A), \theta) \)-twisted orbits.

Moreover, we may assume that \( n_j \in N_{G_k}(A) \cap Q' \) whenever \( Z_{G_k}(A) n_j \cap Q' \neq \emptyset \). After possibly renumbering and discarding irrelevant indices, we can assume that \( n_j \in N_{G_k}(A) \cap Q' \), i.e., \( \theta(n_j) = n_j^{-1} \), for all \( j \in \{1, \dots, \ell\} \).

\textbf{Step 2.} By \textbf{Step 1}, given \( n \in N_{G_k}(A) \cap Q' \), it suffices to study the \((Z_{G_k}(A), \theta)\)-twisted orbits of the set \( Z_{G_k}(A) n \cap Q' \). However, it is more natural and convenient to study the \((Z_{G_k}(A), \theta)\)-twisted orbits of \( Z_{G_k}(A) \cap Q' \), thereby eliminating the explicit appearance of \( n \).  This can be achieved by shifting the perspective from the involution \( \theta \) to a different involution. Indeed, since \( \theta(n) = n^{-1} \), it follows from \cite[Lemma 6.1]{HW93} that the map
\[
\theta_n(x) := n \theta(x) n^{-1}, \quad \text{for all } x \in G,
\]
defines a \( k \)-involution on \( G \).

Moreover, by \cite[Lemma 6.4]{HW93}, and since \( Z_G(A) \) is a \( \theta \)-stable \( k \)-subgroup of \( G \), we can replace the study of \((Z_{G_k}(A), \theta)\)-twisted orbits of \( Z_{G_k}(A) n \cap Q' \) with the study of \((Z_{G_k}(A), \theta_n)\)-twisted orbits of the variety
\[
Q'(n, Z_{G_k}(A)) := \left\{ x \in Z_{G_k}(A) \mid \theta_n(x) = x^{-1} \right\},
\]
thus switching from \( \theta \) to \( \theta_n \) as desired.  From now on, we may assume that
\[
G := Z_G(A) = A \cdot M,
\]
where \( M \) is a $k$-anisotropic subgroup of \( G \) defined over \( k \), and \( A \) is the \( k \)-split central torus of \( G \).

We now consider \( \theta := \theta_n \) as our \( k \)-involution, and define the variety
\[
Q'_k := Q'(n, Z_{G_k}(A)) = \left\{ x \in Z_{G_k}(A) \mid \theta_n(x) = x^{-1} \right\}.
\]
Note that both \( A \) and \( G \) remain \( \theta_n \)-stable.

\textbf{Step 3.}
Let \( G := Z_G(A) = A \cdot M \), where \( A \) is a \( k \)-split and \( \theta \)-stable central torus, and \( M \) is a $k$-anisotropic \( k \)-subgroup. Let \( \theta := \theta_n \), and define
\[
Q' := \left\{ x \in Z_G(A) \mid \theta_n(x) = x^{-1} \right\}.
\]
We aim to show that there are only finitely many \((G_k, \theta)\)-twisted orbits in \( G_k \cap Q' \).  Without loss of generality, we may assume that \( M = G / A \), so \( M \) is a \( \theta \)-stable, $k$-anisotropic \( k \)-subgroup with \( A \cap M = \{ \mathrm{id} \} \).

Let \( x \in Q' \), and we have: $\theta(x) = \theta(a_1 m_1) = \theta(a_1) \theta(m_1) = (a_1 m_1)^{-1} = m_1^{-1} a_1^{-1}$.
Since \( A \) is central in \( G \), we can write $\theta(x) = a_1^{-1} m_1^{-1}$.

Given that both \( A \) and \( M \) are \( \theta \)-stable, and since \( A \cap M = \{ \mathrm{id} \} \), it follows that
\[
\theta(a_1) = a_1^{-1}, \quad \theta(m_1) = m_1^{-1}.
\]
Therefore, we can decompose: $Q' = Q'_A \cdot Q'_M$, where for \( B \in \{A, M\} \), we define $Q'_B := \left\{ b \in B \mid \theta(b) = b^{-1} \right\}$.

Since \( A_k \) is central in \( G_k \), we also have $G_k \cap Q' = (A_k \cap Q'_A) \cdot (M_k \cap Q'_M)$.

Moreover, the \((G_k, \theta)\)-twisted action splits accordingly. For any \( x = a_1 m_1 \in G_k \cap Q' \) and \( g = a m \in G_k = A_k \cdot M_k \), we compute:
\[
g \cdot x \cdot \theta(g)^{-1} = a m \cdot a_1 m_1 \cdot \theta(a m)^{-1} = a a_1  \theta(a)^{-1} \cdot m m_1 \theta(m)^{-1} = a \theta(a)^{-1} a_1  \cdot m m_1 \theta(m)^{-1}.
\]

Hence, it is natural to count separately the number of \((A_k, \theta)\)-twisted orbits in \( A_k \cap Q'_A \), and  the number of \((M_k, \theta)\)-twisted orbits in \( M_k \cap Q'_M \).

\textbf{Step 4.} 
Assume \( G = A \), i.e., \( G \) is a \( k \)-split torus that is \( \theta \)-stable. Since \( A = A_+ \cdot A_- \), we have:
\[
Q'_A = A_- := \left\{ a \in A \mid \theta(a) = a^{-1} \right\},  \quad \tau(A) := \left\{ a \theta(a)^{-1} \mid a \in A \right\} = \left\{ a^2 \mid a \in A_- \right\}.
\]
This implies that the number of \((A_k, \theta)\)-twisted orbits in \( A_k \cap Q'_A \) is given by the coset space
\[
(Q'_A)_k \backslash \tau(A)_k,
\]
since \( A_k \) is abelian and the $(A_k,\theta)$-twisted conjugation reduces to multiplication by \( \tau(a) \).

Now, since \( k \) is a local field, the group \( k^\times / (k^\times)^2 \) is finite. Therefore, the image of \( \tau(A)_k \) in \( A_k \) has finite index.  Hence, the number of \((A_k, \theta)\)-twisted orbits in \( A_k \cap Q'_A \) is finite. This concludes Step 4, using the fact that \( A_k \) is an abelian group of finite \( k^\times \)-dimension.

\textbf{Step 5.} 
Let \( G = M \) be a \( k \)-anisotropic subgroup defined over a local field \( k \), equipped with a \( k \)-involution \( \theta \). We assume that the local field \( k \) is endowed with a discrete valuation \( \nu \), and that Hensel's Lemma holds for the pair \( (k, \nu) \); that is, \( (k, \nu) \) is Henselian. According to \cite[Example 2.1.4]{KaPra} or \cite[Proposition 2.A.5]{Ach}, any field \( k \) equipped with a valuation \( \nu \) and complete with respect to the topology induced by \( \nu \) is Henselian. In particular, since we are working with non-Archimedean local fields of characteristic zero -- which are complete by definition -- this setting ensures that we are in the setting of Henselian fields.

By a theorem of Bruhat and Tits, and also of Guy Rousseau (see \cite{Rou77} or \cite{Pra82}), we know that \( M_k \) is compact with respect to the topology induced by the local field \( k \). This compactness follows from the fact that $k$-anisotropic reductive groups over local fields have bounded \( k \)-points, which is equivalent to compactness in the valuation topology.  Furthermore, by \cite[Section 9]{Richar}, there are only finitely many \((M, \theta)\)-twisted conjugacy classes in the variety
\[
Q'_M := \left\{ m \in M \mid \theta(m) = m^{-1} \right\},
\]
say \( X_1, \dots, X_r \), where each \( X_i \) is a closed \( k \)-subvariety of \( M \).

By the definition of the topology on the set of \( k \)-points of an algebraic variety defined over \( k \), and using results from Section~\ref{sec::inv_fct_thm}, we conclude that:
\[
(Q'_M)_k = Q'_M \cap M_k = \bigcup_{i=1}^{r} (X_i)_k
\]
is a finite union of closed subsets \( (X_i)_k \subset M_k \), each closed with respect to the topology induced from \( k \). Since \( M_k \) is compact, this implies that \( (Q'_M)_k \) is a compact subset of \( M_k \), partitioned into finitely many $(M_k,\theta)$-twisted conjugacy classes.

\medskip
Since $M_k$ is compact, for every closed $X_i$, with $i \in \{1,...,r\}$, the intersection $(X_i)_k = X_i \cap M_k$ is also compact in $M_k$. 

Given $X_i$, suppose we can show that for every $x \in X_i \cap M_k$, the $(M_k, \theta)$-twisted orbit of $x$ in $X_i \cap M_k$ is open. Then, since $(X_i)_k$ is compact we can select a finite number of $(M_k, \theta)$-twisted orbit in $(X_i)_k=X_i \cap M_k$ to cover the compact set  $(X_i)_k$. As there are a finite number of $X_i$ the proof of Step 4 is finished.

It remains to show that given an $(M,\theta)$-twisted orbit $X_i$ in $Q'_M$, that is a closed subvariety of $M$, and a point  $x$ in $X_i \cap M_k$, the $(M_k, \theta)$-twisted orbit $M_k\cdot_\theta x$ of $x$ in $X_i \cap M_k$ is open. Indeed, consider the map $\mu_x : M \to X_i$ given by $m \in M \mapsto \mu_x(m):= m x \theta(m^{-1})= m\cdot_\theta x$. Since $X_i$ is already an $(M,\theta)$-twisted orbit and $x \in X_i$, the map $\mu_x$ is surjective, and dominant, and so its differential $(d\mu_x)_{\id}: T_{\id} (G) \to T_x (X_i)$ is surjective. By \cite[Chapter 3, Proposition 3.3]{PlaRa} (\textbf{where we need that $char(k)=0$}) the orbit  $M_k\cdot_\theta x$ is open in $(M\cdot_\theta x)_k= (X_i)_k$, and we are done.
 \end{proof}

Recall that we have fixed a minimal parabolic \( k \)-subgroup \( P \subset G \), and that \( P \) admits a \( \theta \)-stable maximal \( k \)-split torus \( A \subset P \) of \( G \), which we fix for what follows. Let \( U \) denote the unipotent radical of \( P \), and let \( V \) be a set of representatives for the \( Z_{G_k}(A) \times H_k \)-orbits from (\ref{equ::action_Z_H}) in
\[
\tau^{-1}(N_{G_k}(A)) = \left\{ g \in G_k \mid \tau(g) = g \theta(g^{-1}) \in N_{G_k}(A) \right\}.
\]

By Proposition~\ref{prop::HvP_orbits}, we know that \( G_k \) is the disjoint union of the double cosets \( P_k v H_k \), with \( v \in V \). Any element \( g \in P_k v H_k \) produces the element \( g \theta(g)^{-1} \), which lies in one of the \( (P_k, \theta) \)-orbits in \( Q' \), and in fact in \( Q \). By Proposition~\ref{prop::6.15}, the intersection \( G_k \cap Q' \) contains only finitely many \( (P_k, \theta) \)-twisted orbits. In particular, this also holds for \( G_k \cap Q \).

Combining these two propositions, we obtain the following:

\begin{corollary}[Proposition 6.16 in \cite{HW93}]
\label{cor::6.16}
Let \( k \) be a local field of characteristic zero. Then the double coset space \( P_k \backslash G_k / H_k \) is finite. In particular, the set \( \{A_i \mid i \in I\} \) of representatives of the \( H_k \)-conjugacy classes of \( \theta \)-stable maximal \( k \)-split tori in \( G \) is finite (see Proposition~\ref{prop::diff_coset_decom_G_k}). Moreover, the set of \( H_k \)-conjugacy classes of maximal \( (\theta, k) \)-split tori in \( G \) is also finite.
\end{corollary}

\subsection{Visualizing the decomposition \( K \mathcal{B} H \) in the Bruhat--Tits building}
Finally, we present an intuitive and geometric interpretation of the proof of the polar decomposition \( K \mathcal{B}_k H_k \) of \( G_k \) given in \cite[Theorem 1.1]{BeOh}. Recall that \( \Delta_{G_k} \) denotes the Bruhat--Tits building of \( G_k \), and that the involution \( \theta \) induces an involutive automorphism of \( \Delta_{G_k} \), which we also denote by \( \theta \).

Let \( \Delta_{G_k}^{\theta} := \{ x \in \Delta_{G_k} \mid \theta(x) = x \} \) be the set of \( \theta \)-fixed points in \( \Delta_{G_k} \). The group \( H_k \) stabilizes this geodesically convex and closed subset, and its action on \( \Delta_{G_k}^{\theta} \) admits (in theory) a compact fundamental domain, which we denote by \( F \).

Fix a vertex \( x \in F \), and define
\[
K := \{ g \in G_k \mid g(x) \in F \}.
\]
Note that \( K \) is a compact subset of \( G_k \), although it is not necessarily a subgroup of $G_k$

Now, let \( g \in G_k \). If \( g(x) \in F \), then clearly \( g \in K \). If instead \( g(x) \notin F \), we consider the projection \( y \) of \( g(x) \) onto the geodesically convex and closed subset \( \Delta_{G_k}^{\theta} \). Using an element \( h \in H_k \), we can move the point \( y \in \Delta_{G_k}^{\theta} \) into the fundamental domain \( F \). Consequently, \( h(g(x)) \) projects to the point \( h(y) \in \Delta_{G_k}^{\theta} \).

Now, recall that there are only finitely many \( H_k \)-conjugacy classes of maximal \( (\theta, k) \)-split tori in \( G \). Geometrically, these can be viewed as a finite collection of flats in \( \Delta_{G_k} \) that are ``\( \theta \)-perpendicular'' to \( \Delta_{G_k}^{\theta} \). Since we are working modulo \( H_k \)-conjugacy, we may assume that these flats emanate ``perpendicularly'' from the fundamental domain \( F \subset \Delta_{G_k}^{\theta} \).

Thus, the point \( h(g(x)) \) lies in one of these flats associated with a maximal \( (\theta, k) \)-split torus of $G$. Each such torus has a compact fundamental domain on its corresponding flat that intersects $\Delta_{G_k}^{\theta}$, so we can use an element \( a \) from the torus to move \( h(g(x)) \) back into \( F \). This implies that \( a h g(x) \in F \), and hence \( a h g \in K \), which shows that \( g \in H_k \mathcal{B}_k K \). This gives the desired decomposition.

\section{Appendix B: On the Inverse Function Theorem and its consequences }
\label{appen:B}

\label{sec::inv_fct_thm}

We revisit the Inverse Function Theorem within the framework of analytic manifolds defined over non-Archimedean fields, highlighting its implications for algebraic group actions. These foundational results play a role in establishing openness and transitivity properties, which are essential to the study of Chabauty limits of fixed-point groups arising from involutive automorphisms.

In the context of results referenced from \cite{Borel}, the group \( G \) is assumed to be a connected reductive affine algebraic group, which is equivalent to a connected reductive linear algebraic group. This equivalence implies that \( G \) can be embedded as a subgroup of a general linear group, and its points correspond to the solutions of a finite system of polynomial equations over a fixed algebraically closed field (with unspecified characteristic).

% \begin{definition}[See Chapter IV, Definition 11.13 from \cite{Borel}]
% \label{def::Cartan_subgroup}
% A Cartan subgroup of \( G \) is the centralizer of a maximal torus in \( G \).
% \end{definition}

The notions of the radical and the unipotent radical of a reductive or semisimple group \( G \) are defined in \cite[Chapter IV, Section 11.21, p.157]{Borel}.

\begin{proposition}[See Chapter IV, Corollary 14.19 from \cite{Borel}]
\label{prop::Levi_centralizers}
Let \( P \) be a parabolic subgroup of \( G \). Then \( P \) contains Levi subgroups. The maximal tori of the radical \( \mathcal{R}P \) of \( P \) are also Cartan subgroups of \( \mathcal{R}P \). The Levi subgroups of \( P \) are precisely the centralizers of the maximal tori of \( \mathcal{R}P \). Moreover, any two Levi subgroups of \( P \) are conjugate by a unique element of the unipotent radical \( \mathcal{R}_u P \) of \( P \).
\end{proposition}

\begin{proposition}[See Chapter IV, Proposition 14.21 from \cite{Borel}]
\label{prop::opposite_parabolics} 
Let \( P \) be a parabolic subgroup of \( G \), and let \( L \) be a Levi subgroup of \( P \).
\begin{enumerate}
\item
There exists a unique parabolic subgroup \( P' \) of \( G \) that is opposite to \( P \) and contains \( L \). Moreover, any two parabolic subgroups of \( G \) that are opposite to \( P \) are conjugate by a unique element of the unipotent radical \( \mathcal{R}_u P \) of \( P \).
\item
Two parabolic subgroups \( P \) and \( Q \) of \( G \) contain opposite Borel subgroups if and only if their unipotent radicals intersect trivially, i.e., only at the identity element.
\item
If \( P' \) is opposite to \( P \), then the product map \( \mu: (x, y) \mapsto x \cdot y \) induces an isomorphism of varieties from \( \mathcal{R}_u P' \times P \) onto an open subset of \( G \), equal to \( P' \cdot P \).
\end{enumerate}
\end{proposition}

We emphasize that all algebraic groups considered here are affine, and hence linear. Since simple points always exist in a variety (see Chapter 2, Section 2.4.3 of \cite{PlaRa}), and since algebraic groups are homogeneous spaces, it follows that every affine algebraic group is smooth -- that is, each of its points is simple (see p.\,98 of \cite{PlaRa}). Furthermore, by Chapter AG, Section 18.1 of \cite{Borel}, every simple point in a variety is normal. Consequently, every affine algebraic group is both smooth and normal. Additionally, the direct product of two normal varieties is itself normal.

\begin{theorem}[Theorem 17.3 from \cite{Borel}, Chapter AG]
Let \( \alpha : V \to W \) be a morphism of varieties. Then the following conditions are equivalent:
\begin{enumerate}
\item
The morphism \( \alpha \) is dominant (i.e., \( \alpha(V) \) is dense in \( W \) with respect to the Zariski topology) and separable.
\item
There exists a dense open subvariety \( V_0 \subset V \) such that the differential \( d_x \alpha \) is surjective for all \( x \in V_0 \).
\item
In each irreducible component of \( V \), there exists a simple point \( x \in V \) such that \( \alpha(x) \) is a simple point of \( W \), and \( d_x \alpha \) is surjective.
\end{enumerate}
\end{theorem}

We also cite a result from \cite{Vust}, where the groups under consideration are connected reductive algebraic groups \( G \), and the base field \( k \) is algebraically closed of characteristic zero. Let \( \theta \) be an involutive automorphism of \( G \), distinct from the identity, and let \( H \subset G \) denote the subgroup of fixed points under \( \theta \). 

We say that the parabolic subgroups \( P \) and \( \theta(P) \) of \( G \) are \textbf{\( \theta \)-split} if they are opposite. In \cite{Vust}, the term \textbf{\( \theta \)-anisotropic} is used to refer to this notion of \( \theta \)-splitness.

\begin{theorem}[1.3 Theorem 1 from \cite{Vust}]
\label{thm::Vust_open_split}
Let \( G \) be a connected reductive affine algebraic group, and the base field \( k \) is algebraically closed of characteristic zero.  
If \( P \) is a parabolic subgroup of \( G \) and \( P \) is \( \theta \)-split, then the product \( H^0 P \) is open in \( G \), where \( H^0 \) denotes the connected component of the identity in \( H \).
\end{theorem}

\begin{proof}
We recall the proof here. Let \( \mathfrak{h}, \mathfrak{p}, \mathfrak{g} \) denote the Lie algebras of \( H^0 \), \( P \), and \( G \), respectively.  
Consider the natural action of \( H^0 \times P \) on the variety \( G \), defined by
\[
(h, p, g) \in H^0 \times P \times G \mapsto \mu(h, p, g) := hgp \in G.
\]

Taking the identity element \( e \in G \), the Closed Orbit Lemma (see \cite[Chapter AG, Corollary 10.2]{Borel} or \cite[Proposition 2.23]{PlaRa}) implies that the orbit of \( e \) under \( H^0 \times P \) is a smooth variety that is open in its closure in \( G \).

Define the map $\mu_e : H^0 \times P \to G, \quad \mu_e(h, p) := hep = hp.$

If we can show that \( \mathfrak{g} = \mathfrak{h} + \mathfrak{p} \), then the differential \( d_e \mu_e \) is surjective. By the homogeneity of the group action, this implies that \( d_v \mu_e \) is surjective for every \( v \in H^0 \times P \). Consequently, by Theorem 17.3 from \cite[Chapter AG]{Borel}, the map \( \mu_e \) is dominant, meaning that \( H^0 e P = H^0 P \) is dense in \( G \) with respect to the Zariski topology. Since \( H^0 e P \) is also open in its closure by the Closed Orbit Lemma, it follows that \( H^0 P \) is open in \( G \).

To establish that \( \mathfrak{g} = \mathfrak{h} + \mathfrak{p} \), we use the fact that \( P \) and \( \theta(P) \) are opposite parabolic subgroups. By Proposition \ref{prop::opposite_parabolics}, there exists an isomorphism of varieties
\[
\mathcal{R}_u \theta(P) \times P \cong U \subset G,
\]
where \( U \) is an open subset of \( G \). This implies that every element \( a \in \mathfrak{g} \) can be expressed as
$a = b + d_e \theta(c)$ or some \( b, c \in \mathfrak{p} \). Rewriting this as $a = (b - c) + (c + d_e \theta(c))$, we observe that \( b - c \in \mathfrak{p} \) and \( c + d_e \theta(c) \in \mathfrak{h} \), which confirms the surjectivity of \( d_e \mu_e \).
\end{proof}

Although an Inverse Function Theorem does not generally apply in the context of algebraic varieties, there exists a well-defined version for analytic manifolds over a field \( k \), as presented in \cite[Part II, Chapter 3]{Serre}. The setup in \cite{Serre} assumes that \( k \) is a field complete with respect to a non-trivial absolute value -- this includes the real numbers \( \mathbb{R} \), complex numbers \( \mathbb{C} \), or a \( p \)-adic field.

Analytic manifolds with charts defined over such a field \( k \) are introduced in Part II, Chapter 3, Section 3 of \cite{Serre}. However, when studying the analytic structure of the \( k \)-points \( V_k \) of an affine algebraic variety \( V \) defined over \( k \), the assumption that the characteristic of \( k \) is zero appears to be essential (see \cite{PlaRa}, middle of page 110).

Given two analytic manifolds \( X \) and \( Y \) over the field \( k \), a map \( f : X \to Y \) is called a \textit{morphism} (or \textit{analytic function}) if it is continuous and locally given by analytic functions (see Part II, Chapter 3, Section 5 of \cite{Serre}). Furthermore, if there exist maps \( f : X \to Y \) and \( g : Y \to X \) such that \( g \circ f = \mathrm{id}_X \) and \( f \circ g = \mathrm{id}_Y \), then \( f \) is called an \textit{isomorphism} if and only if both \( f \) and \( g \) are morphisms.

For the definition of the tangent space at a point of an analytic variety, see Part II, Chapter 3, Section 8 of \cite{Serre}. The Inverse Function Theorem is stated as follows:

\begin{theorem}[See Part II, Chapter 3, Section 9 of \cite{Serre}]
\label{thm::inverse_thm}
Let \( X \) and \( Y \) be analytic manifolds over the field \( k \) of characteristic zero, and let \( x \in X \), \( y \in Y \), with \( \phi : X \to Y \) a morphism such that \( \phi(x) = y \). Then the following conditions are equivalent:
\begin{enumerate}
\item
\( \phi \) is a local isomorphism at \( x \); that is, there exist open neighborhoods \( U \subset X \) of \( x \) and \( V \subset Y \) of \( y \) such that \( \phi : U \to V \) is an isomorphism.
\item
The tangent map \( T_x \phi : T_x X \to T_y Y \) is an isomorphism.
\item
The cotangent map \( T_x^* \phi : T_y^* Y \to T_x^* X \) is an isomorphism.
\end{enumerate}
\end{theorem}

The Immersion Function Theorem provides a criterion for when a morphism between analytic manifolds behaves locally like an embedding. It is stated as follows:

\begin{theorem}[See Part II, Chapter 3, Section 10 of \cite{Serre}]
\label{thm::immersion_thm}
Let \( X \) and \( Y \) be analytic manifolds over a field \( k \) of characteristic zero, with dimensions \( n \) and \( m \), respectively. Let \( x \in X \), \( y \in Y \), and let \( \phi : X \to Y \) be a morphism such that \( \phi(x) = y \). Then the following conditions are equivalent:
\begin{enumerate}
\item
The tangent map \( T_x \phi : T_x X \to T_y Y \) is injective.
\item
There exist open neighborhoods \( U \subset X \) of \( x \), \( V \subset Y \) of \( y \), and \( W \subset k^{n-m} \) containing \( 0 \), along with an isomorphism \( \psi : V \to U \times W \), such that:
\begin{enumerate}
\item[a)]
\( \phi(U) \subset V \),
\item[b)]
If \( i \) denotes the inclusion \( U \hookrightarrow U \times \{0\} \subset U \times W \), then the following diagram commutes:
\begin{center}
\begin{tikzcd}
U \arrow{r}{\phi}  \arrow{rd}{i} 
  &V  \arrow{d}{\psi} \\
    &  U \times W,
\end{tikzcd}
\end{center}
\end{enumerate}
\item
There exist open neighborhoods \( U \subset X \) of \( x \), \( V \subset Y \) of \( y \), and a morphism \( \sigma : V \to U \) such that \( \phi(U) \subset V \) and \( \sigma \circ \phi = \mathrm{id}_U \).
\end{enumerate}
\end{theorem}

In the same analytic setting, there exists a Submersion Function Theorem that characterizes when a morphism locally resembles a projection. The statement is as follows:

\begin{theorem}[See Part II, Chapter 3, Section 10 of \cite{Serre}]
\label{thm::submersion_thm}
Let \( X \) and \( Y \) be analytic manifolds over a field \( k \) of characteristic zero, with dimensions \( n \) and \( m \), respectively. Let \( x \in X \), \( y \in Y \), and let \( \phi : X \to Y \) be a morphism such that \( \phi(x) = y \). Then the following conditions are equivalent:
\begin{enumerate}
\item
The tangent map \( T_x \phi : T_x X \to T_y Y \) is surjective.
\item
There exist open neighborhoods \( U \subset X \) of \( x \), \( V \subset Y \) of \( y \), and \( W \subset k^{m-n} \) containing \( 0 \), along with an isomorphism \( \psi : U \to V \times W \), such that:
\begin{enumerate}
\item[a)]
\( \phi(U) = V \),
\item[b)]
If \( p \) denotes the projection \( V \times W \to V \), then the following diagram commutes:
\begin{center}
\begin{tikzcd}
U \arrow{r}{\phi}   \arrow{rd}{\psi}
  &V   \\
    &  V \times W  \arrow{u}{p},
\end{tikzcd}
\end{center}
\end{enumerate}
\item
There exist open neighborhoods \( U \subset X \) of \( x \), \( V \subset Y \) of \( y \), and a morphism \( \sigma : V \to U \) such that \( \phi(U) \subset V \) and \( \phi \circ \sigma = \mathrm{id}_V \).
\end{enumerate}
\end{theorem}

From Theorem \ref{thm::submersion_thm}, one deduces the following result. Recall that projection maps are open in the analytic category.

\begin{proposition}[See Proposition 3.2 from \cite{PlaRa}, Chapter 3]
\label{prop::open_map_submersion}
Let \( k \) be a non-discrete, locally compact field of characteristic zero. Let \( X \) and \( Y \) be analytic manifolds over \( k \), and let \( x \in X \), \( y \in Y \), with \( \phi : X \to Y \) a morphism such that \( \phi(x) = y \) (in the sense of \cite[Part II, Chapter 3]{Serre}). If the tangent map \( T_x \phi : T_x X \to T_y Y \) is surjective, then \( \phi \) is an open map at \( x \); that is, there exist open neighborhoods \( U \subset X \) of \( x \) and \( V \subset Y \) of \( y \) such that \( \phi : U \to V \) is an open map and \( \phi(U) = V \).
\end{proposition}

Let us now examine the framework established in Chapter 3 of \cite{PlaRa}. The authors work under the assumption that \( k \) is a non-discrete, locally compact field of characteristic zero. Using the topology on \( k \), one can define a natural topology on the set of \( k \)-points \( V_k \) of an algebraic variety \( V \) defined over \( k \).

This is constructed as follows: for any Zariski-open \( k \)-subset \( U \subset V \) and a finite collection of regular \( k \)-functions \( f_1, \ldots, f_r \) on \( U \), one considers the subset
\[
V(f_1, \ldots, f_r; \epsilon) := \left\{ x \in U_k \;\middle|\; |f_i(x)|_v < \epsilon \text{ for all } i = 1, \ldots, r \right\},
\]
where \( |\cdot|_v \) denotes the norm on \( k \) induced by a valuation \( v \). The collection of all such subsets forms a basis for a topology on \( V_k \), known as the \( v \)-adic topology. This topology is strictly finer than the Zariski topology and enjoys better topological properties.

For instance:
\begin{itemize}
\item If \( V = V_1 \times V_2 \) is a product of two \( k \)-varieties, then \( V_k \) is canonically homeomorphic to \( (V_1)_k \times (V_2)_k \), equipped with the product topology.
\item If \( W \subset V \) is a Zariski-open (resp. closed) \( k \)-subvariety, then \( W_k \subset V_k \) is open (resp. closed) in the \( v \)-adic topology.
\end{itemize}

As stated in \cite[middle of p.\,111]{PlaRa} -- without proof and citing \cite{Serre} -- when \( V \) is an affine \( k \)-variety, the set of simple points of \( V_k \) naturally carries the structure of an analytic manifold over \( k \). In particular, if \( G \) is an affine algebraic group defined over \( k \), then its group of \( k \)-points \( G_k \) (which are all simple) inherits an analytic manifold structure. This allows one to apply the Inverse Function Theorem and its consequences to \( G_k \).

Finally, recall that for affine (i.e., linear) algebraic groups, connectedness is equivalent to irreducibility. Thus, when citing Proposition 3.2 from \cite{PlaRa}, we do so in the context of connected linear algebraic groups defined over \( k \), where the analytic and algebraic structures on the \( k \)-points coincide.

\end{document}